\documentclass[preprint,11pt]{elsarticle}
\biboptions{sort,compress} 

\usepackage[a4paper]{geometry}
\usepackage{lmodern} 

\usepackage{amssymb}
\usepackage{amsfonts}
\usepackage{amsmath}
\usepackage{mathtools} 
\usepackage{xfrac} 
\usepackage{gensymb} 
\usepackage{booktabs} 
\usepackage{siunitx} 

\usepackage{xcolor} 
\usepackage{graphicx}

\usepackage{caption}

\usepackage{floatrow}
\usepackage{pgf}
\usepackage{subfig}

\usepackage{multirow}
\usepackage{dashrule}

\usepackage{acro} 
\usepackage{lineno} 
\usepackage[ruled, vlined, english, onelanguage, algo2e]{algorithm2e}

\usepackage{mdframed}

\definecolor{julia-purple}{HTML}{9558B2}
\definecolor{julia-green}{HTML}{389826}
\definecolor{julia-blue}{HTML}{4063D8}
\definecolor{julia-red}{HTML}{CB3C33}

\definecolor{gray03}{gray}{0.3}
\definecolor{gray04}{gray}{0.4}
\definecolor{gray08}{gray}{0.8}
\definecolor{gray09}{gray}{0.9}

\definecolor{minted-red}{HTML}{b00040}

\definecolor{gray08}{gray}{0.8}

\usepackage{makecell} 
\usepackage{hyperref}
\usepackage[capitalise]{cleveref}
\Crefname{equation}{}{}

\newcommand{\orcid}[1]{\href{https://orcid.org/#1}{\includegraphics[width=10pt]{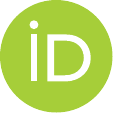}}}

\numberwithin{equation}{section}

\usepackage{fancyhdr}
\usepackage{arydshln} 

\newcommand*\eps{\varepsilon} 
\newcolumntype{?}[1]{!{\vrule width #1pt}} 

\Crefname{algocf}{Algorithm}{Algorithms}

\DeclareAcronym{ODE}{
	short = ODE,
	long  = ordinary differential equation
}

\DeclareAcronym{IVP}{
	short = IVP,
	long  = initial value problem
}

\DeclareAcronym{RHS}{
	short = RHS,
	long = right-hand side
}

\DeclareAcronym{PDE}{
	short = PDE,
	long  = partial differential equation
}

\DeclareAcronym{MoL}{
	short = MoL,
	long  = method of lines
}

\DeclareAcronym{SSP}{
	short = SSP,
	long  = strong stability preserving
}

\DeclareAcronym{TVD}{
	short = TVD,
	long  = total variation diminishing
}

\DeclareAcronym{TVB}{
	short = TVB,
	long  = total variation bounded
}

\DeclareAcronym{DG}{
	short = DG,
	long  = Discontinuous Galerkin
}

\DeclareAcronym{FV}{
	short = FV,
	long  = finite volume
}

\DeclareAcronym{SBP}{
	short = SBP,
	long  = summation--by--parts
}

\DeclareAcronym{WF}{
	short = WF,
	long  = weak form
}

\DeclareAcronym{FD}{
	short = FD,
	long  = flux--differencing
}

\DeclareAcronym{VT}{
	short = VT,
	long  = volume term
}

\DeclareAcronym{VI}{
	short = VI,
	long  = volume integral
}

\DeclareAcronym{VT/VI}{
	short = VT/VI,
	long  = volume term/volume integral
}

\DeclareAcronym{RKDG}{
	short = RKDG,
	long  = Runge-Kutta Discontinuous Galerkin
}

\DeclareAcronym{DGSEM}{
	short = DGSEM,
	long  = discontinuous Galerkin spectral element method
}

\DeclareAcronym{LLF}{
	short = LLF,
	long  = local Lax--Friedrichs
}

\DeclareAcronym{HLLC}{
	short = HLLC,
	long  = Harten--Lax--Van Leer--Contact
}

\DeclareAcronym{HLLE}{
	short = HLLE,
	long  = Harten--Lax--Van Leer--Einfeldt
}

\DeclareAcronym{HLL}{
	short = HLL,
	long  = Harten--Lax--Van Leer
}

\DeclareAcronym{MHD}{
	short = MHD,
	long  = magnetohydrodynamics
}

\DeclareAcronym{vrMHD}{
	short = VRMHD,
	long  = visco-resistive Magnetohydrodynamics
}

\DeclareAcronym{GLM}{
	short = GLM,
	long  = generalized Lagrangian multiplier
}

\DeclareAcronym{glm-mhd}{
	short = GLM-MHD,
	long  = generalized Lagrangian multiplier magnetohydrodynamics
}

\DeclareAcronym{PRKM}{
	short = PRKM,
	long = Partitioned Runge-Kutta method
}
\DeclareAcronym{ARKMs}{
	short = ARKMs,
	long = Additive Runge-Kutta methods
}
\DeclareAcronym{PERK}{
	short = P-ERK,
	long = Paired Explicit Runge-Kutta
}

\DeclareAcronym{PERRK}{
	short = P-ERRK,
	long = Paired Explicit Relaxation Runge-Kutta
}

\DeclareAcronym{RRKM}{
	short = RRKM,
	long = Relaxation Runge-Kutta method
}

\DeclareAcronym{IMEX}{
	short = IMEX,
	long = implicit-explicit
}

\DeclareAcronym{AMR}{
	short = AMR,
	long = adaptive mesh refinement
}

\DeclareAcronym{CFL}{
	short = CFL,
	long = Courant--Friedrichs--Lewy
}

\DeclareAcronym{DoF}{
	short = DoF,
	long = degree of freedom
}

\DeclareAcronym{IDT}{
	short = IDT,
	long = incremental direction technique
}

\DeclareAcronym{BR1}{
	short = BR1,
	long = Bassi--Rebay 1
}

\DeclareAcronym{EC/ES}{
	short = EC/ES,
	long = entropy--conservative/entropy--stable
}

\DeclareAcronym{EC}{
	short = EC,
	long = entropy--conservative
}

\DeclareAcronym{ES}{
	short = ES,
	long = entropy--stable
}

\DeclareAcronym{CRM}{
	short = CRM,
	long = Common Research Model
}

\DeclareAcronym{RANS}{
	short = RANS,
	long = Reynolds--averaged Navier--Stokes
}

\DeclareAcronym{MOOD}{
	short = MOOD,
	long = multi-dimensional optimal order detection
}

\journal{??}

\begin{document}

\mdfsetup{
	backgroundcolor=gray08,
	roundcorner=10pt,
	topline=false,
	rightline=false,
	bottomline=false,
	leftline=false}
	
	\begin{frontmatter}
		
		
		
		\title{Volume Term Adaptivity for Discontinuous Galerkin Schemes}
		
		
		\cortext[cor1]{Corresponding author}
		
		\author[1]{Daniel Doehring\corref{cor1} \orcid{0009-0005-4260-0332}}
		\author[2]{Jesse Chan \orcid{0000-0003-2077-3636}}
		\author[3]{Hendrik Ranocha \orcid{0000-0002-3456-2277}}
		\author[4]{\\Michael Schlottke--Lakemper \orcid{0000-0002-3195-2536}}
		\author[1]{Manuel Torrilhon \orcid{0000-0003-0008-2061}}
		\author[5]{and Gregor Gassner \orcid{0000-0002-1752-1158}}
		
		\affiliation[1]{
			organization={Applied~and~Computational~Mathematics,~RWTH~Aachen~University},
			country={Germany}.
		}

		\affiliation[2]{
			organization={Oden~Institute~for~Computational~Science~and~Engineering, Aerospace Engineering~\&~Engineering~Mechanics, ~University~of~Texas~at~Austin},
			country={USA}.
		}

		\affiliation[3]{
			organization={Institute~of~Mathematics,~Johannes~Gutenberg~University~Mainz},
			country={Germany.}
		}
		
		\affiliation[4]{
			organization={High-Performance Computing, Center for Advanced Analytics and Predictive Sciences, \\University of Augsburg},
			country={Germany.}
		}
		
		\affiliation[5]{
			organization={Department of Mathematics and Computer Science, Center for Data and Simulation Science, \\University of Cologne},
			country={Germany.}
		}
		
		\begin{abstract}
			We introduce the concept of volume term adaptivity for high--order discontinuous Galerkin (DG) schemes solving time--dependent partial differential equations.
			Termed \textit{v}--adaptivity, we present a novel general approach that exchanges the discretization of the volume contribution of the DG scheme at every Runge--Kutta stage based on suitable indicators.
			Depending on whether robustness or efficiency is the main concern, different adaptation strategies can be chosen.
			Precisely, the weak form volume term discretization is used instead of the entropy-conserving flux--differencing volume integral whenever the former produces more entropy than the latter, resulting in an entropy--stable scheme.
			Conversely, if increasing the efficiency is the main objective, the weak form volume integral may be employed as long as it does not increase entropy beyond a certain threshold or cause instabilities.
			Thus, depending on the choice of the indicator, the \textit{v}--adaptive DG scheme improves robustness, efficiency and approximation quality compared to schemes with a uniform volume term discretization.
			We thoroughly verify the accuracy, linear stability, and entropy--admissibility of the \textit{v}--adaptive DG scheme before applying it to various compressible flow problems in two and three dimensions.
		\end{abstract}
		
		%
		%
		
		%
		%
		
		
	\end{frontmatter}
	
	
	
	
	%
	\section{Introduction}
	Adaptivity is a necessity for the feasible simulation of complex physical phenomena with high--fidelity, high--order methods.
	When referring to adaptivity, we mean the ability of a numerical method to automatically adjust the discretization technique to the current local solution approximation over the course of the simulation.
	Adaptivity is a well--known concept in the numerical treatment of \acp{PDE} and has been studied extensively in the past decades.
	Established terms are \textit{h--adaptivity}, often also referred to as \ac{AMR} \cite{babuvvska1978error, osher1983numerical, berger1984adaptive, berger1989local, leveque2002finite}, where the cell size is locally adjusted to capture small--scale features of the solution, for instance shock fronts and boundary layers.
	Next, \textit{p--adaptivity} \cite{babuska1981p, biswas1994parallel, houston2001hp, schwab1998p, balsara2016efficient} refers to the local adjustment of an element's polynomial degree $p$ which is particularly useful to capture smooth solution features with high accuracy, such as vortical structures in turbulent flows.
	Third, \textit{r--adaptivity} \cite{miller1981moving1, miller1981moving2, huang1994moving, huang2010adaptive} denotes methods where the cells are moved and deformed to better align them with solution features such as discontinuities.
	Adaptivity of the approximating scheme similar to the this work has also been explored.
	Some examples related to the present work are order adaptation \cite{clain2011high, diot2012improved, diot2013multidimensional, dumbser2014posteriori, zanotti2015space, dumbser2016simple, hennemann2021provably, rueda2021entropy}, order--preserving stabilization \cite{FISHER2013high, carpenter2014entropy, abgrall2018general, abgrall2022reinterpretation, bilocq2025comparison, chan2025artificial} and explicit shock--capturing strategies \cite{hennemann2021provably, rueda2021entropy, lin2024high}.
	Furthermore, for time--dependent problems solved with explicit time integrators, adaptive strategies to alleviate the time step restriction imposed by the \ac{CFL} condition have been developed.
	Classic examples thereof are local time stepping \cite{osher1983numerical, berger1984adaptive, dawson2001high, krivodonova2010efficient} which reduce the timestep locally per cell, and multirate time integration methods \cite{gunther2001multirate, constantinescu2007multirate, seny2013multirate, vermeire2019paired, doehring2024multirate} which employ methods with enlarged regions of absolute stability for cells with more restricted timestep sizes.
	Finally, we would like to remark that adaptivity is not limited to discretization techniques, but is also extended to the modeling itself, i.e., different physical models are employed in dynamically adapted regions of the domain, see for instance \cite{tiwari1998adaptive, kolobov2012towards, mathis2015dynamic, shou2021magnetohydrodynamic, allmann2024muphyii} and references therein.

	In this work, we focus on \ac{DG} methods \cite{reed1973triangular, cockburn2001runge, hesthaven2007nodal} and their application to time--dependent hyperbolic--parabolic \acp{PDE}, most notably the compressible Euler and Navier--Stokes--Fourier equations.
	A standard weak/variational formulation of a \ac{PDE} with a \ac{DG} discretization consists of two main parts: The surface term and the volume term, which arise from the application of integration by parts to the integrated \ac{PDE}.
	The treatment of the surface term boils down to the numerical solution of Riemann problems at element interfaces, for which a wealth of approximate Riemann solvers/numerical flux functions is available \cite{leveque2002finite, RUSANOV1962304, harten1983upstream, einfeldt1988godunov, toro1994restoration, toro2009riemann}.
	The treatment of the volume term, however, leaves more freedom in the design of the numerical scheme.
	This flexibility is crucial to obtain \ac{EC/ES} \cite{tadmor1986entropy, tadmor1987numerical, tadmor2003entropy} 
	\ac{DG} schemes
	which employ a special type of volume term discretization, the so--called \textit{\ac{FD}} form \cite{FISHER2013high, gassner2013skew,  carpenter2014entropy, gassner2016split, CHEN2017427, chan2018discretely}, which greatly enhances the stability of these methods \cite{SJOGREEN2018153, MANZANERO2020109241, rojas2021robustness, chan2022entropy, bilocq2025comparison}.

	In general, the usage of the (in general) more robust \ac{EC/ES} \ac{DG} schemes comes at the cost of a more involved volume term discretization.
	In particular, for the \ac{FD} form, the volume term evaluation requires the computation of special two--point flux functions at pairs of the quadrature nodes of the element.
	This is especially troubling for element types which do not naturally possess a tensor--product/sum factorization structure, such as triangles, tetrahedra, pyramids, and prisms unless specialized techniques are employed, see for instance \cite{montoya2024efficient} and references therein.
	But even for tensor--product elements such as quadrilaterals and hexahedra, the \ac{FD} form is for higher polynomial degrees more expensive than a standard \ac{WF} volume integral, which only involves a single analytical flux evaluation per quadrature node.
	This observation serves as the main motivation for this work:
	We seek to use the cheaper weak form volume term as often as possible, with the option to switch to the more expensive \ac{FD} form only when necessary to retain stability of the simulation.
	The choice of the volume term discretization is made locally per element and Runge--Kutta stage based on a suitable indicator, which may be of a--priori or a--posteriori nature.
	This gives rise to a new type of adaptivity for \ac{DG} schemes, which we refer to as \textit{volume term--} or in short \textit{v--adaptivity}.
	In some sense, this idea is a generalization of the approach taken in \cite{hennemann2021provably, rueda2021entropy} wherein the entropy--conservative \ac{FD} volume integral is blended with a lower--order finite volume scheme defined on the \ac{DG} subcells to provide stabilization.
	A blueprint for the a--posteriori adaptivity are the works on \ac{MOOD} \cite{clain2011high, diot2012improved, diot2013multidimensional} and its extension to \ac{DG} schemes \cite{dumbser2014posteriori, zanotti2015space, dumbser2016simple}, where the solution is checked after computation for admissibility given by certain indicators and adjusted otherwise.
	Another related approach is proposed in \cite{FISHER2013high, carpenter2014entropy} where the finite--difference scheme is equipped with a more stable finite--volume discretization to stabilize the simulation.
	In a recent work \cite{bilocq2025comparison} the authors proposed combining the standard \ac{DG} weak form with the entropy stable \ac{DG} scheme from \cite{chan2019skew} on Gauss--Legendre points based on the indicator from \cite{persson2006sub, hennemann2021provably}.
	We generalize this approach to general quadrature rules and general indicators.

	Besides enhancing efficiency, it is also possible to improve the robustness of the discretization by using the \ac{WF} update for an element whenever it dissipates more entropy than the \ac{FD} \ac{VT/VI} on that element, which is a rigorous criterion for entropy--stability/diffusion.
	The key to efficiency for this approach is that one can compute the analytical element--wise entropy production, which is identical to the entropy production of the \ac{FD} \ac{VT/VI}, without the need to compute the \ac{FD} volume term itself.
	Instead, only the element--wise surface integral of the entropy potential is required \cite{CHEN2017427, chan2018discretely, lin2024high, chan2025artificial}.

	The paper is structured as follows:
	We begin by reviewing the classic \ac{WF} and \ac{EC/ES} \ac{FD} volume term discretizations of \ac{DG} schemes in \cref{sec:vol_term_formulations_DG}.
	Next, in \cref{sec:indicator} we present an indicator based on the element--wise entropy production rate to determine the appropriate volume term.
	In particular, we devise a simple rigorous entropy--stable criterion that increases robustness of the simulation, and a heuristic criterion that aims to reduce computational cost without impairing stability significantly.
	Additionally, we briefly revisit the indicator from \cite{persson2006sub, hennemann2021provably} as this will also be used to demarcate between the \ac{WF} and the \ac{FV}--stabilized \ac{FD} \ac{VT}.
	In \cref{sec:Verification} we verify the newly obtained entropy--dissipative adaptive volume integral in terms of order of accuracy, linear stability, and entropy--stability.
	Examples for coupling the \ac{WF} with pure \ac{FD} and stabilized \ac{FD}--\ac{FV} are presented in \cref{sec:Examples}, before we conclude the work in \cref{sec:Conclusions}.
	\section{Volume Term Formulations in Nodal Discontinuous Galerkin Schemes}
	\label{sec:vol_term_formulations_DG}
	In this section, we briefly review the major steps in the construction of the \ac{EC/ES} \ac{DG} schemes based on the \ac{FD} form.
	We restrict ourselves here for simplicity to the \ac{DGSEM} \cite{black1999conservative, hesthaven2007nodal, kopriva2009implementing} with collocated Gauss--Legendre--Lobatto interpolation/quadrature nodes on tensor--product elements.
	The interested reader is referred to \cite{kopriva2010quadrature, chan2018discretely, chan2019efficient, chan2019discretely, chan2022entropy} for more general choices of the quadrature rule and to \cite{CHEN2017427, chan2019skew, chen2020review} for a discussion on modal and non tensor--product elements.
	\subsection{Weak Formulation of the 1D Discontinuous Galerkin Spectral Element Method}
	Consider the vector--valued conservation law in one spatial dimension
	\begin{equation}
		\partial_t \boldsymbol u(t, x) + \partial_x \boldsymbol f(\boldsymbol u) = \boldsymbol 0,	
	\end{equation}
	which can be transformed from a physical element $\Omega_i = [x_{i-\sfrac{1}{2}}, x_{i+\sfrac{1}{2}}]$ to the reference element $\widehat{\Omega} \coloneqq [-1, 1]$ via the affine mapping
	\begin{equation}
		\xi(x) \coloneqq 2\frac{x - x_{i-\sfrac{1}{2}}}{x_{i+\sfrac{1}{2}} - x_{i-\sfrac{1}{2}}} - 1, \quad x \in \Omega_i
	\end{equation}
	and inverse
	\begin{equation}
		x_i(\xi) \coloneqq x_{i - \sfrac{1}{2}} + \frac{x_{i+\sfrac{1}{2}} - x_{i-\sfrac{1}{2}}}{2} (\xi + 1), \quad \xi \in \widehat{\Omega}
	\end{equation}
	with the Jacobian of the transformation
	\begin{equation}
		J_i \coloneqq \frac{\partial x_i(\xi)}{\partial \xi} = \frac{x_{i+\sfrac{1}{2}} - x_{i-\sfrac{1}{2}}}{2} \: .
	\end{equation}
	Thus, the conservation law in reference space reads
	\begin{subequations}
		\label{eq:1D_conservation_law_reference}
		\begin{align}
			\partial_t \boldsymbol u\big(t, x_i(\xi)\big) + \frac{\partial \xi}{\partial \xi} \frac{\partial}{\partial x} \boldsymbol f\big(\boldsymbol u(t, x_i(\xi))\big) &= \boldsymbol 0 \\
			\Leftrightarrow \partial_t \boldsymbol u\big(t, x_i(\xi)\big) + J_i^{-1} \frac{\partial}{\partial \xi} \boldsymbol f\big(\boldsymbol u(t, x_i(\xi))\big) &= \boldsymbol 0 \\
			\label{eq:1D_conservation_law_reference3}
			\Leftrightarrow J_i \partial_t \boldsymbol u\big(t, x_i(\xi)\big) + \partial_\xi \boldsymbol f\big(\boldsymbol u(t, x_i(\xi))\big) &= \boldsymbol 0 \: .
		\end{align}
	\end{subequations}
	We expand the solution $\boldsymbol u$ in terms of Lagrange polynomials $\ell_j(\xi)$ of degree $p$ defined on the Gauss--Lobatto nodes $\xi_0, \ldots, \xi_p \in \widehat{\Omega}$ as
	\begin{equation}
		\label{eq:DG_solution_expansion}
		\boldsymbol u(t, x_i(\xi)) \coloneqq \sum_{j=0}^p \boldsymbol u_{i,j}(t) \ell_j(\xi) \: ,	
	\end{equation}
	with time--dependent degrees of freedom $\boldsymbol u_{i,j}(t)$ which are defined for every element $\Omega_i$ and node $j = 0, \ldots, p$.
	We insert this expansion into \cref{eq:1D_conservation_law_reference3} and make the traditional Galerkin ansatz by multiplying with test functions identical to ansatz functions $\ell_k(\xi)$ and integrating over the reference element $\widehat{\Omega}$ to obtain
	\begin{subequations}
		\label{eq:DG_weak_form}
		\begin{align}
			\boldsymbol 0 &= \int_{-1}^1 \left[ J_i \partial_t \boldsymbol u\big(t, x_i(\xi)\big) + \partial_\xi \boldsymbol f(\boldsymbol u ) \right] \ell_k(\xi)d\xi \: , \quad k = 0, \ldots, p \\
			\label{eq:DG_weak_form2}
   		& = J_i \int_{-1}^1 \partial_t \boldsymbol u \big(t, x_i(\xi)\big) \ell_k(\xi)d\xi + \int_{-1}^1 \partial_\xi \boldsymbol f(\boldsymbol u) \ell_k(\xi)d\xi \: , \quad k = 0, \ldots, p \: .
		\end{align}
	\end{subequations}
	The first term in \cref{eq:DG_weak_form2} can be simplified by inserting \cref{eq:DG_solution_expansion}
	\begin{subequations}
		\label{eq:DG_weak_form_first_term}
		\begin{align}
			J_i \int_{-1}^1 \partial_t \boldsymbol u \big(t, x_i(\xi)\big) \ell_k(\xi)d\xi &= J_i \int_{-1}^1 \sum_{j=0}^p \partial_t \boldsymbol u_{i,j}(t) \ell_j(\xi) \ell_k(\xi) d\xi \\
			\label{eq:DG_weak_form_first_term2}
			&= J_i \sum_{l=0}^p \omega_l \sum_{j=0}^p \partial_t \boldsymbol u_{i,j}(t) \ell_j(\xi_l) \ell_k(\xi_l)
		\end{align}
	\end{subequations}
	with $\omega_l$ denoting the quadrature weights of the Gauss--Lobatto quadrature rule.
	Due to the collocation property of the Lagrange polynomials, i.e., $\ell_j(\xi_l) = \delta_{jk}$, \cref{eq:DG_weak_form_first_term2} simplifies to
	\begin{equation}
		J_i \int_{-1}^1 \partial_t \boldsymbol u \big(t, x_i(\xi)\big) \ell_k(\xi)d\xi = J_i \omega_k \partial_t \boldsymbol u_{i,k}(t) = J_i \omega_k \dot{\boldsymbol u}_{i,k}(t)
	\end{equation}
	which is can be expressed in matrix notation with diagonal \textit{mass matrix}
	\begin{equation}
		\label{eq:mass_matrix}
		\mathbf{M}_{ij} \coloneqq \delta_{ij} \omega_i, \quad i, j = 0, \ldots, p
	\end{equation}
	as
	\begin{equation}
		J_i \int_{-1}^1 \partial_t \boldsymbol u \big(t, x_i(\xi)\big) \ell_k(\xi)d\xi = J_i \mathbf{M} \underline{\dot{\boldsymbol u}}_i(t)
	\end{equation}
	where we assemble the coefficients $\boldsymbol u_{i,j}(t)$ at quadrature nodes $j$ into
	\begin{equation}
		\underline{\boldsymbol u}_i(t) \coloneqq \begin{pmatrix}
			\boldsymbol u_{i,0}(t) \\
			\vdots \\
			\boldsymbol u_{i,p}(t)
		\end{pmatrix}.
	\end{equation}
	Next, we consider the second summand in \cref{eq:DG_weak_form2} and apply integration by parts to obtain
	\begin{equation}
		\label{eq:SurfaceVolumeTermSplit}
		\int_{-1}^1 \partial_\xi \boldsymbol f(\boldsymbol u) \ell_k(\xi)d\xi
		= \left[ \boldsymbol f \big( \boldsymbol u (t, x_i(\xi))\big) \ell_k(\xi) \right]_{-1}^1 - \int_{-1}^1 \boldsymbol f ( \boldsymbol u ) \ell_k'(\xi) d\xi.
	\end{equation}
	The first term on the right--hand side of \cref{eq:SurfaceVolumeTermSplit} represents the interface or \textit{surface term} of the \ac{DG} scheme which simplifies with $\ell_k(-1) = \delta_{k0}$ and $\ell_k(1) = \delta_{kp}$ to
	\begin{equation}
		\label{eq:DG_surface_term}
		\left[ \boldsymbol f \big( \boldsymbol u (t, x_i(\xi))\big) \ell_k(\xi) \right]_{-1}^1
		= \boldsymbol f\big( \boldsymbol u(t, x_{i+\sfrac{1}{2}}) \big) \delta_{kp} - \boldsymbol f\big( \boldsymbol u(t, x_{i-\sfrac{1}{2}}) \big) \delta_{k0}.
	\end{equation}
	Analogous to \ac{FV} schemes, the physical flux evaluations at element interfaces with potentially discontinuous states $\boldsymbol u\left(t, x^{L}_{i \pm \sfrac{1}{2}}\right) \neq \boldsymbol u\left(t, x^{R}_{i \pm \sfrac{1}{2}}\right)$ are replaced by numerical flux functions $\boldsymbol f^{\star}$ which approximate the solution of the Riemann problem.
	By introducing the \textit{boundary matrix}
	\begin{equation}
		\label{eq:boundary_matrix}
		\mathbf{B} \coloneqq \begin{pmatrix}
			-1 & 0 & \cdots & 0 \\
			0 & 0 & \cdots & 0 \\
			\vdots & \vdots & \ddots & \vdots \\
			0 & 0 & \cdots & 1
		\end{pmatrix}	
	\end{equation}
	we can further simplify \cref{eq:DG_surface_term} to
	\begin{equation}
		\left[ \boldsymbol f \big( \boldsymbol u (t, x_i(\xi))\big) \ell_k(\xi) \right]_{-1}^1
		= \mathbf{B} \begin{pmatrix} \boldsymbol f^\star\vert_{-1} \\ 0 \\ \vdots \\ 0 \\ \boldsymbol f^\star\vert_{1} \end{pmatrix} \eqqcolon \mathbf{B} \underline{\boldsymbol f}_i^\star \: .
	\end{equation}
	Finally, we come to the \textit{volume term} in \cref{eq:SurfaceVolumeTermSplit} for which we apply the quadrature rule to obtain
	\begin{equation}
		- \int_{-1}^1 \boldsymbol f ( \boldsymbol u ) \ell_k'(\xi) d\xi = - \sum_{l=0}^p \omega_l \boldsymbol f\big( \boldsymbol u(t, x_i(\xi_l)) \big) \ell_k'(\xi_l) \: .
	\end{equation}
	With the \textit{derivative matrix}
	\begin{equation}
		\label{eq:derivative_matrix}
		\mathbf{D}_{ij} \coloneqq \ell_j'(\xi_i), \quad i, j = 0, \ldots, p	
	\end{equation}
	and mass matrix from \cref{eq:mass_matrix}, we can express the volume term in matrix notation as
	\begin{equation}
		- \int_{-1}^1 \boldsymbol f ( \boldsymbol u ) \ell_k'(\xi) d\xi = - \mathbf{D}^T \mathbf{M} \underline{\boldsymbol f}_i \quad k = 0, \ldots, p 
	\end{equation}
	which gives the complete weak form \ac{DG} semidiscretization in reference space on element $\Omega_i$ as
	\begin{subequations}
		\label{eq:DG_weak_form_ref_space}
		\begin{align}
			& \, J_i \mathbf{M} \underline{\dot{\boldsymbol u}}_i(t) + \mathbf{B} \underline{\boldsymbol f}_i^\star - \mathbf{D}^T \mathbf{M} \underline{\boldsymbol f}_i = \boldsymbol 0 \\
			\label{eq:DG_weak_form_ref_space2}
			\Leftrightarrow 
			& \, J_i \underline{\dot{\boldsymbol u}}_i(t) = - \mathbf{M}^{-1} \mathbf{B} \underline{\boldsymbol f}_i^\star + \mathbf{M}^{-1} \mathbf{D}^T \mathbf{M} \underline{\boldsymbol f}_i \: .
		\end{align}
	\end{subequations}
	\subsection{Flux--Differencing Formulation of the 1D Discontinuous Galerkin Spectral Element Method}
	For the sake of notational compactness we drop the element index $i$ in the following.
	A crucial observation for the construction of \ac{EC/ES} \ac{DG} schemes is that the mass, boundary, and derivative matrices fulfill the \textit{\ac{SBP} property} \cite{kreiss1974finite, strand1994summation, olsson1995summation1, olsson1995summation2, carpenter1996spectral, gassner2013skew}
	\begin{equation}
		\label{eq:SBP_property}
		\underbrace{\mathbf{M} \mathbf{D}}_{\eqqcolon \mathbf{Q}} + (\mathbf{M} \mathbf{D})^T = \mathbf{B} \: .
	\end{equation}
	Employing this property, one can rewrite $\mathbf{D}^T \mathbf{M} = \mathbf{D}^T \mathbf{M}^T = (\mathbf{M} \mathbf{D})^T = \mathbf{B} - \mathbf{M} \mathbf{D}$ and thus express the volume term as
	\begin{equation}
		\mathbf{M}^{-1} \mathbf{D}^T \mathbf{M} \underline{\boldsymbol f} = \mathbf{M}^{-1} \left( \mathbf{B} - \mathbf{M} \mathbf{D} \right) \underline{\boldsymbol f}
	\end{equation}
	which can be grouped into boundary and volume contributions as
	\begin{equation}
		\label{eq:flux_differencing_grouped}
		\mathbf{M}^{-1} \mathbf{D}^T \mathbf{M} \underline{\boldsymbol f} = \mathbf{M}^{-1} \mathbf{B} \underline{\boldsymbol f} - \mathbf{D} \underline{\boldsymbol f} \: .
	\end{equation}
	The key insight from \cite{FISHER2013high} is that if the derivative operator $\mathbf{D}$ fulfills a telescoping property, i.e., every row sums to zero, one can rewrite
	\begin{equation}
		\label{eq:flux_differencing_volume_term}
		\left(\mathbf{D} \underline{\boldsymbol f} \right)_j = \sum_{k=0}^p \mathbf{D}_{j,k} \boldsymbol f_k = \sum_{k=0}^p \mathbf{D}_{j,k} \boldsymbol f_k + \boldsymbol f_j \underbrace{\sum_{k=0}^p \mathbf{D}_{j,k}}_{= \boldsymbol 0} = 2 \sum_{k=0}^p \frac{1}{2} \mathbf{D}_{j,k} \left( \boldsymbol f_j + \boldsymbol f_k \right)
	\end{equation}
	which is the \ac{FD} form of the volume term.
	Note that the \ac{WF} volume integral is equivalent to the \ac{FD} volume integral with the central two--point flux employed as the \textit{volume flux}.
	Replacing the central two--point flux
	\begin{equation}
		\boldsymbol f^\text{Central} \left(\boldsymbol f_j, \boldsymbol f_k \right) = \frac{1}{2} \left( \boldsymbol f_j + \boldsymbol f_k \right)
	\end{equation}
	with an appropriate entropy--conservative and consistent two--point flux function $\boldsymbol f^{\text{EC}} \left( \boldsymbol u_j, \boldsymbol u_k \right)$ \cite{ismail2009affordable, Chandrashekar_2013, winters2016affordable, ranocha2018comparison, Ranocha2020Entropy} leads to an entropy--conservative volume term discretization \cite{gassner2016split}.
	We thus replace $-\mathbf{D} \underline{\boldsymbol f}$ in \cref{eq:flux_differencing_grouped} by $-2 \mathbf{D} \underline{\boldsymbol f}^\text{EC}$ and reformulate the \ac{DG} scheme \cref{eq:DG_weak_form_ref_space2} using the above manipulations as
	\begin{subequations}
		\label{eq:DG_flux_differencing_form_split_matrix}
		\begin{align}
			J \underline{\dot{\boldsymbol u}}(t) &= - \mathbf{M}^{-1} \mathbf{B} \underline{\boldsymbol f}^\star + \mathbf{M}^{-1} \mathbf{B} \underline{\boldsymbol f}^\text{EC} - \mathbf{D} \underline{\boldsymbol f}^\text{EC}	\\
			&= - \mathbf{M}^{-1} \mathbf{B} \left( \underline{\boldsymbol f}^\star - \underline{\boldsymbol f}^\text{EC} \right) - 2 \mathbf{D} \underline{\boldsymbol f}^\text{EC} \\
			\label{eq:DG_flux_differencing_form_split_matrix3}
			&= - \mathbf{M}^{-1} \mathbf{B} \underline{\boldsymbol f}^\star - \left( 2 \mathbf{D} - \mathbf{M}^{-1} \mathbf{B} \right) \underline{\boldsymbol f}^\text{EC} \: .
		\end{align}
	\end{subequations}
	The matrix $\mathbf{D}_\text{split} \coloneqq 2 \mathbf{D} - \mathbf{M}^{-1} \mathbf{B}$ is typically referred to as \textit{derivative--split matrix} \cite{gassner2016split}.

	A closer inspection of \cref{eq:flux_differencing_volume_term} reveals why the \ac{FD} approach is computationally more expensive than the \ac{WF} volume term.
	The \ac{WF} volume term requires $d \cdot (p + 1)^d$ \underline{analytical/physical flux} evaluations to construct $\underline{\boldsymbol f}$ in $d$ spatial dimensions.
	In contrast, the \ac{FD} discretization requires $\sfrac{d p}{2} \cdot (p+1)^d$ \underline{two--point numerical} flux evaluations per element, i.e., there is quadratic growth in the number of \ac{EC} flux evaluations with increasing polynomial degree $p$.
	The ratio of flux evaluations between the \ac{FD} and \ac{WF} volume terms is thus for all $d$ given by $\sfrac{p}{2}$.
	In other words, if the \ac{WF} \ac{VT} would be naively implemented in \ac{FD} form, the simulation would be about $\sfrac{p}{2}$ times more expensive and take correspondingly longer solely due to the increased number of flux evaluations.

	For non--tensor product/sum factorization elements, though, the situation is dramatically worse, as now two--point fluxes need to  be evaluated between all pairs of nodes within an element.
	For an element with $N$ nodes, this results in
	\begin{equation}
		N_f = \sum_{j = 1}^{N-1} j= \frac{1}{2} N(N-1)
	\end{equation}
	flux evaluations per element.
	Assuming now that the number of nodes $N$ scales with the polynomial degree $p$ as $N \sim p^d$, we obtain a scaling of $N_f \sim p^{2d}$ which is much worse than the $p^{d+1}$ scaling for tensor--product elements.

	Additionally, the computation of the entropy--conserving two--point fluxes is more expensive than the analytical flux evaluation.
	While the analytical/physical flux can be computed using the most basic arithmetic operations, the computation of the \ac{EC} fluxes involves usually some logarithmic means \cite{ismail2009affordable, Chandrashekar_2013, ranocha2018comparison, winters2016affordable} which considerably increase the computational cost per flux evaluation.
	To illustrate this, we present runtimes for the evaluation of the analytical flux and the \ac{EC} flux by Chandrashekar \cite{Chandrashekar_2013} and Ranocha \cite{ranocha2018comparison} as implemented in \texttt{Trixi.jl} \cite{trixi1, trixi2, trixi3} for the compressible Euler equations in \cref{tab:CostFluxEval}.
	The measurements are obtained for averaging over all directions $i = 1, \ldots, d$ and $10^7$ evaluations for blast--alike initial states $\rho_l = \rho_r = 1$, $\boldsymbol v_l = \boldsymbol v_r = \boldsymbol 0$, $p_l = 10^4$, and $p_r = 1$.
	For the analytical flux, the tabulated value is the average of $\boldsymbol f(\boldsymbol u_l), \boldsymbol f(\boldsymbol u_r)$.
	Clearly, the \ac{EC} fluxes are about $1.4$ to $1.6$ times more expensive than the analytical flux evaluation.
	\begin{table}
		\def\arraystretch{1.3}
		\centering
		\begin{tabular}{c?{2}c|c|c}
			 & \multicolumn{3}{c}{$\tau$ [ns]} \\
			 & 1D & 2D & 3D \\
			\Xhline{5\arrayrulewidth}
			Analytical                              & $2.1$ & $2.3$ & $2.3$ \\
			\hdashline
			Chandrashekar \cite{Chandrashekar_2013} & $11.8$ & $13.5$ & $14.8$ \\
			Ranocha \cite{ranocha2018comparison}    & $10.1$ & $12.0$ & $12.3$
		\end{tabular}
		\caption[]
		{Mean runtimes $\tau$ for the evaluation of the analytical compressible Euler flux and the two--point entropy--conservative fluxes by Chandrashekar \cite{Chandrashekar_2013} and Ranocha \cite{ranocha2018comparison} in one, two, and three spatial dimensions.}
		\label{tab:CostFluxEval}
	\end{table}
	\section{Indicator}
	\label{sec:indicator}
	To decide which volume term discretization to use in the $v$--adaptive framework, a suitable indicator is required.
	For the combination of the weak form and flux--differencing form volume terms we propose an indicator based on the cell--wise entropy production of the numerical scheme.
	In particular, one can compute the change in mathematical entropy $S(\boldsymbol u)$ as
	\begin{equation}
		\dot{S} = \partial_t S \big( \boldsymbol u(t) \big) = \frac{\partial S}{\partial \boldsymbol u} \partial_t \boldsymbol u(t)= \boldsymbol w^T \dot{\boldsymbol u}
	\end{equation}
	where we introduced the entropy variables $\boldsymbol w \coloneqq \frac{\partial S}{\partial \boldsymbol u}$ \cite{tadmor2003entropy}.
	The entropy production due to the \ac{VT} can thus be computed as
	\begin{subequations}
		\label{eq:entropy_production_VT}
		\begin{align}
			\dot S\vert_{\text{WF}} &= \boldsymbol w^T \cdot \dot{\boldsymbol u} \vert_\text{WF} \overset{\cref{eq:DG_weak_form_ref_space2}}{=} - J^{-1} \mathbf{M}^{-1} \mathbf{D}^T \mathbf{M} \underline{\boldsymbol f} \\
			\dot S\vert_{\text{FD}} &= \boldsymbol w^T \cdot \dot{\boldsymbol u} \vert_\text{FD} \overset{\cref{eq:DG_flux_differencing_form_split_matrix3}}{=} - J^{-1} \left( 2 \mathbf{D} - \mathbf{M}^{-1} \mathbf{B} \right) \underline{\boldsymbol f}^\text{EC} \: .
		\end{align}
	\end{subequations}
	The discrete scheme satisfies a cell entropy inequality if we can bound the entropy production of the volume term by the surface integral of the \textit{entropy potential} $\psi$ \cite{CHEN2017427, chan2018discretely, lin2024high, chan2025artificial} 
	\begin{equation}
		\label{eq:entropy_production_VT_bound}
		\dot S\vert_{\text{VT}} \leq \int_{\partial \Omega} \psi(\boldsymbol u) \, \mathrm{d}S = \left[ \psi(\boldsymbol u)\right]_{-1}^{+1}
	\end{equation}
	and entropy--dissipative surface fluxes $\boldsymbol f^\star$ are employed.
	The entropy potential $\psi$ \cite{tadmor2003entropy} is defined via the relation
	\begin{equation}
		\label{eq:Def_entropy_pot}
		\psi\big(\boldsymbol w(\boldsymbol u)\big) \coloneqq \boldsymbol w^T \boldsymbol f(\boldsymbol u) - F(\boldsymbol u)
	\end{equation}
	where $F(\boldsymbol u)$ is the \textit{entropy flux} defined via \cite{tadmor2003entropy}
	\begin{equation}
		\label{eq:Def_entropy_flux}
		\left( \frac{\partial F}{\partial \boldsymbol u} \right)^T = \boldsymbol w^T \left( \frac{\partial \boldsymbol f}{\partial \boldsymbol u} \right) \: .
	\end{equation}
	The (component--wise) entropy potential for the compressible Euler equations is given simply by the momentum components, i.e., in one dimension $\psi(\boldsymbol u) = \rho v$ \cite{tadmor2003entropy}.
	\subsection{Rigorous Indicator for more Robust Entropy--Stable Schemes}
	\label{subsec:RigorousIndicator}
	We now turn our attention to designing a numerical scheme that is globally entropy dissipative by combining the \ac{WF} and \ac{FD} \ac{VT} discretizations.
	In particular, we seek to construct a more stable adaptive scheme by using the \ac{WF} volume term if it introduces more decay in the mathematical entropy than the \ac{FD} form, i.e., if
	\begin{equation}
		\label{eq:entropy_production_comparison_naive}
		\Delta S_i \coloneqq \dot{S}_i \vert_{\text{WF}} - \dot{S}_{i} \vert_{\text{FD}} < 0
	\end{equation}
	for element $\Omega_i$.
	This way, we obtain a globally entropy--stable scheme (assuming the surface terms dissipate entropy) as we only use the \ac{WF} volume term when it is more dissipative than the \ac{FD} form.
	A related approach for the entire semidiscretization, i.e., not only the volume integral, was already suggested in \cite{carpenter2014entropy} where also an adaptive scheme based on local entropy production was presented.
	In particular, the authors suggest equipping an entropy--conservative scheme with an entropy--diffusive \textit{companion} scheme and blending them based on local entropy production.
	In \cite{carpenter2014entropy} WENO and other finite volume based schemes are suggested as companion schemes, with the former being employed in \cite{FISHER2013high}.
	Thus, in the terminology of \cite{carpenter2014entropy}, the herein presented approach uses the \ac{WF}/central--flux \ac{FD} volume integral as the diffusive \textit{companion} scheme for the entropy--conserving \ac{FD} volume integral.
	Other closely related works are the a--posteriori subcell limiting \ac{DG} schemes developed in \cite{dumbser2014posteriori, zanotti2015space, dumbser2016simple} which are based on the \ac{MOOD} paradigm \cite{clain2011high, diot2012improved, diot2013multidimensional}.
	Adaptive or switching schemes, albeit not necessarily for entropy--stability, have also been studied in \cite{hou1994nonconservative, sjogreen2023construction, worku2024entropy}, for instance.

	The indicator as presented in \cref{eq:entropy_production_comparison_naive} is unfortunately very expensive, as the \ac{FD} volume term is always computed, even if it is not used in the final update.
	As pointed out in \cite{lin2024high, chan2025artificial}, this can be avoided by replacing $\dot{S}_{i} \vert_\text{FD}$ by
	\begin{equation}
		\dot{S}_{i} \vert_\text{FD} = \psi\big( \boldsymbol{u}_{i,p} \big) - \psi\big( \boldsymbol{u}_{i,0} \big)\: .
	\end{equation}
	Employing this in the indicator equation \cref{eq:entropy_production_comparison_naive} thus gives
	\begin{equation}
		\label{eq:entropy_production_comparison}
		\Delta S_i \coloneqq \dot{S}_{i} \vert_{\text{WF}} - \Big( \psi \big( \boldsymbol{u}_{i,p} \big) - \psi\big( \boldsymbol{u}_{i,0} \big) \Big) < 0 \: .
	\end{equation}
	This way, we only need to compute the entropy production of the \ac{WF} volume term and the surface integral of the entropy potential to evaluate the indicator.
	The \ac{FD} volume term is in this case only computed if entropy is indeed increasing.
	We remark that in the case where for a majority of cells the \ac{WF} \ac{VT} may be used, this is not only a more robust scheme, but also more efficient scheme for polynomial degrees $p > 2$ as the \ac{FD} volume term is only computed on a subset of elements.

	The indicator \cref{eq:entropy_production_comparison} is computed at every Runge--Kutta stage and element in the $v$--adaptive framework.
	Note that this indicator is of \underline{a--posteriori} nature, as the \ac{WF} \ac{VT} is computed in any case, but potentially discarded.
	We present numerical results for this scheme showcasing increased robustness over the pure \ac{FD} form in \cref{sec:Examples}.
	\subsection{Heuristic Indicator for more Efficient Schemes}
	If efficiency is the main concern for a simulation, one can relax the indicator by introducing an upper bound $\sigma_i$ for the tolerated entropy production of the \ac{WF} volume term.
	In particular, the \ac{WF} volume term is used if
	\begin{equation}
		\label{eq:EntropyTargetDecay}
		\dot{S}_{i, \text{WF}} < \sigma_i \: .
	\end{equation}
	The \ac{WF} volume integral serves in this case as a predictor, and correction via the \ac{FD} volume integral is only applied if deemed necessary by the indicator.

	The difficulty lies now in finding a reasonable upper bound $\sigma_i$ for the admissible entropy production of the \ac{WF} volume term.
	Is this bound too large, the \ac{WF} volume integral is potentially applied too often, leading to instabilities.
	On the other hand, if the bound is too small, the \ac{FD} volume term computation is used most of the time, leading to ultimately higher computational cost than a pure \ac{FD} scheme.
	We emphasize here that $\sigma$ applies to the entropy production in reference space, which eliminates issues with magnitude across elements of different sizes and enables better reusability of the parameter across different meshes.
	We experimented with different fixed and dynamic bounds, the latter computed from the erronous total entropy production $\dot S > 0$ observed for a complete Runge--Kutta step.
	In practice, the dynamic approach did not yield improvements over a fixed bound, and thus we present results only for the latter case in \cref{sec:Examples}.
	Similar to identifying a stable CFL number, we determined a stable $\sigma$ by performing bisection starting with e.g. $\sigma_\text{max} = 10^{-2}$ and $\sigma_\text{min} = 0.0$.
	More sophisticated approaches are left for future work.
	\subsection{Shock--Capturing Indicator}
	\label{sec:shock_capturing_indicator}
	In the case where the weak form volume term is combined with the stabilizing blended \ac{DG}--\ac{FV} we typically employ the a--priori troubled--cell indicator from \cite{persson2006sub, hennemann2021provably}.
	The indicator converts the nodal representation in Lagrange polynomials into a modal representation in terms of orthogonal Legendre polynomials via 
	\begin{equation}
		\label{eq:nodal_to_modal}
		\underline{\boldsymbol m} = \mathbf{V}(\boldsymbol \xi)^{-1} \underline{\boldsymbol v}
	\end{equation}
	with $\mathbf V(\boldsymbol \xi)$ denoting the Vandermonde matrix of the Legendre polynomials evaluated at the interpolation nodes $\boldsymbol \xi = \left\{ \xi_0, \dots, \xi_p \right\}$.
	Note that $\underline{\boldsymbol v} = \{ v_j \}$ in \cref{eq:nodal_to_modal} is the nodal representation of an indicator quantity such as density $\rho$, pressure $p$, or the product of these two $\rho \cdot p$.
	From the modal coefficients $\boldsymbol m$ the energies
	\begin{equation}
		\mathcal{E}_n \coloneqq \sum_{j = 0}^{p - n} m_j^2, \quad n = 0, 1, 2
	\end{equation}
	are computed.
	The influence of the higher modes, i.e., potential oscillations, is then measured by the fractions
	\begin{equation}
		\eps_0 \coloneqq \frac{\mathcal{E}_0 - \mathcal{E}_1}{\mathcal{E}_0}, \quad \eps_1 \coloneqq \frac{\mathcal{E}_1 - \mathcal{E}_2}{\mathcal{E}_1}
	\end{equation}
	from which the maximum $\eps \coloneqq \max(\eps_0, \eps_1)$ is computed.
	Then, the blending factor $\beta$ is computed from the standard logistic function as
	\begin{equation}
		\label{eq:blending_factor_beta}
		\beta \coloneqq \frac{1}{1 + \exp \left( -\kappa \frac{\eps - \mathcal{T}}{\mathcal{T}} \right)}	
	\end{equation}
	with hyperparameters 
	\begin{equation}
		\kappa \coloneqq \log \left( \frac{1 - 10^{-4}}{10^{-4}} \right) \approx 9.21024, \quad \mathcal{T} \coloneqq 0.5 \cdot 10^{1.8 \cdot p^{0.25}} 
		\: .
	\end{equation}
	To avoid unnecessary blending in the case of a very small $\beta$, a threshold $\beta_\text{min} \geq 0$ is supplied, below which the indicator value is clipped to zero.
	Similarly, one can limit the maximum blending factor to $\beta_\text{max} \leq 1$ to avoid excessive blending in the case of a very large $\beta$:
	\begin{equation}
		\beta = 
		\begin{cases}
			0 & \beta < \beta_\text{min} \\
			\beta & \beta_\text{min} \leq \beta \leq \beta_\text{max} \\
			\beta_\text{max} & \beta > \beta_\text{max}
		\end{cases}
		\: .
	\end{equation}
	The indicator is computed at every Runge--Kutta stage and element in the blended \ac{DG}--\ac{FV} framework with indicator variable $\boldsymbol v$ taken from the previous Runge--Kutta stage.
	Thus, the blending is determined a--priori any volume term computations of the current stage.
	\subsubsection{Shock--Capturing via Blended DG--FV}
	The blended \ac{DG}--\ac{FV} scheme uses a convex combination of a high--order \ac{DG} discretization and a non--oscillatory first or second--order \ac{FV} discretization in troubled cells identified by the shock indicator $\beta$.
	The volume term update is then given by
	\begin{equation}
		\label{eq:BlendedDG-FV}
		\dot{\boldsymbol u}_i \vert_\text{DG--FV} \coloneqq (1 - \beta_i) \dot{\boldsymbol u}_i \vert_\text{FD} + \beta_i \dot{\boldsymbol u}_i \vert_\text{FV} \, , \quad \beta_i \in [0, \beta_\text{max}] \, ,
	\end{equation}
	where $\dot{\boldsymbol u}_i \vert_\text{FV}$ is a standard finite--volume update on the subcells of the $[-1, 1]^d$ reference element.
	The \textit{blending parameter} $\beta_i$ determines the amount of dissipation introduced by the \ac{FV} update.
	For a $p$--degree \ac{DG} scheme, $(p+1)^d$ subcells are used per element in $d$ spatial dimensions.
	The natural size of these subcells is given by the quadrature weights $\omega_i$ of the Lobatto--Legendre basis \cite{hennemann2021provably}.
	Thus, the first subcell is given by $[-1, -1 + \omega_1]$, the second by $[-1 + \omega_1, -1 + \omega_1 + \omega_2]$, and so on.
	The scheme \cref{eq:BlendedDG-FV} is provably entropy stable and the main shock--capturing methodology used in \texttt{Trixi.jl} \cite{trixi1,trixi2, trixi3}.
	\subsection{Combining Shock--Capturing and Entropy--Production Indicators}
	Due to the different nature (a--priori vs a--posteriori) of the shock--capturing and $v$--adaptive indicators, they can be readily combined to give a nested \ac{VT} computation strategy.
	In particular, the shock--capturing indicator is computed first to decide whether any stabilization via blended \ac{DG}--\ac{FV} is required.
	If not, the \ac{VT} for this element is computed with the weak form, and the $v$--adaptive indicator is evaluated to decide whether the \ac{FD} volume term is required for this element to maintain e.g. entropy stability.

	In practice, however, we have for the \ac{DGSEM} (i.e., tensor--product elements) not found examples where this hybrid strategy is more efficient to maintain \ac{EC/ES} properties than using always the \ac{FD} volume term in non--stabilized regions.
	We believe that this hybrid strategy could be beneficial for non tensor--product elements where the cost of the \ac{FD} volume term is much higher than the \ac{WF} volume term, but this is left for future work.
	%
	\section{Verification}
	\label{sec:Verification}
	All numerical tests and examples are performed with \texttt{Trixi.jl} \cite{trixi1, trixi2, trixi3}, a high--order \ac{DG} code written in \texttt{Julia} \cite{bezanson2017julia} for the numerical simulation of transport--dominated hyperbolic--parabolic balance laws and related multiphysics problems.
	All computations performed in this section can be reproduced without the need for commercial and proprietary software.
	The files to reproduce the results are publicly available on \texttt{GitHub} and archived on \texttt{Zenodo} \cite{doehring2026VTA_ReproRepo}.
	\subsection{Convergence}
	\subsubsection{Isentropic Vortex Advection with Viscous Surface Flux}
\begin{figure}
		\centering
		\subfloat[{Domain--normalized $L^2$ error in density.}]{
			\label{fig:IVA_Convergence_L2}
			\resizebox{.46\textwidth}{!}{\includegraphics{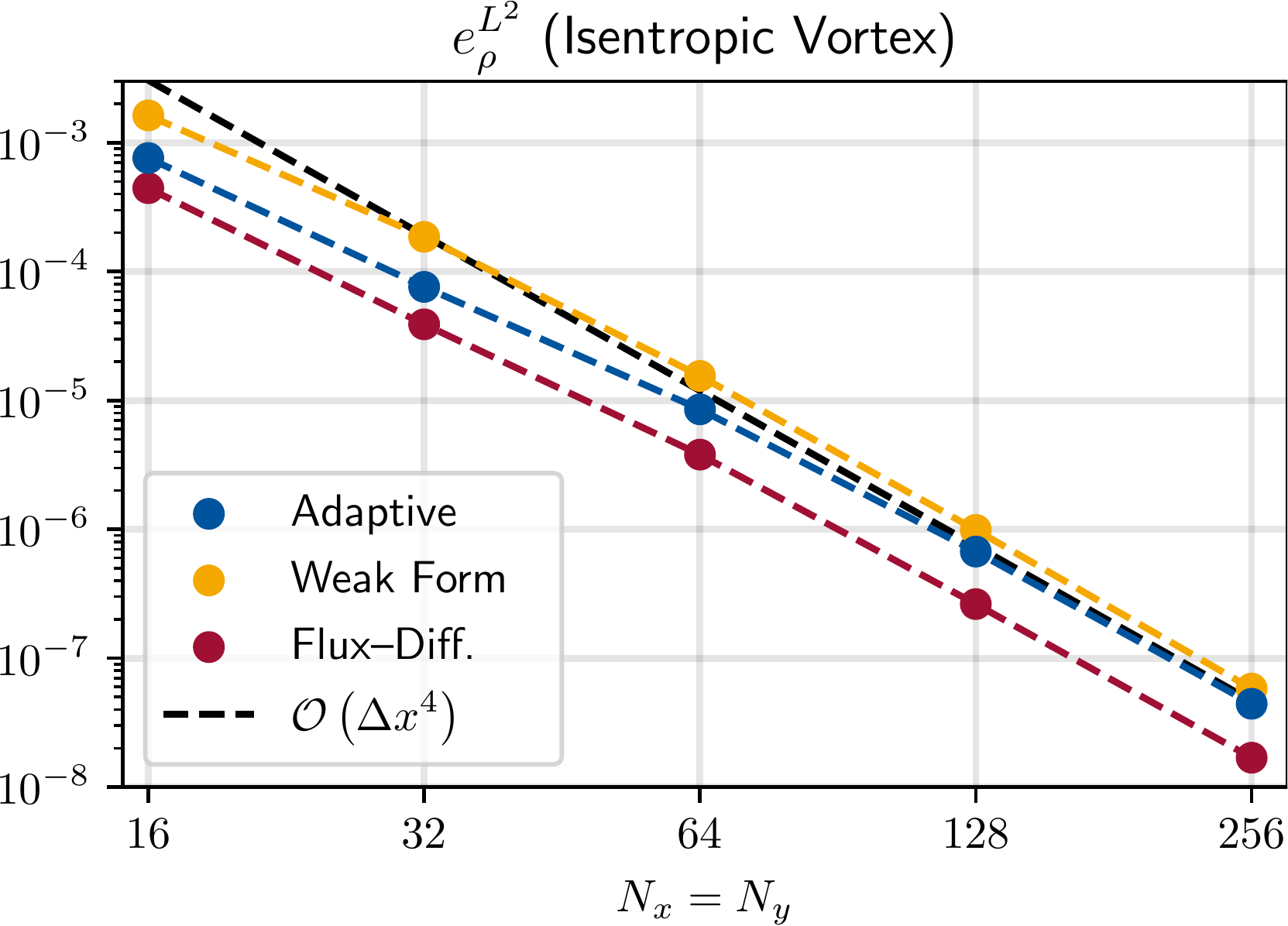}}
		}
		\hfill
		\subfloat[{Pointwise $L^\infty$ error in density.}]{
			\label{fig:IVA_Convergence_LInf}
			\resizebox{.46\textwidth}{!}{\includegraphics{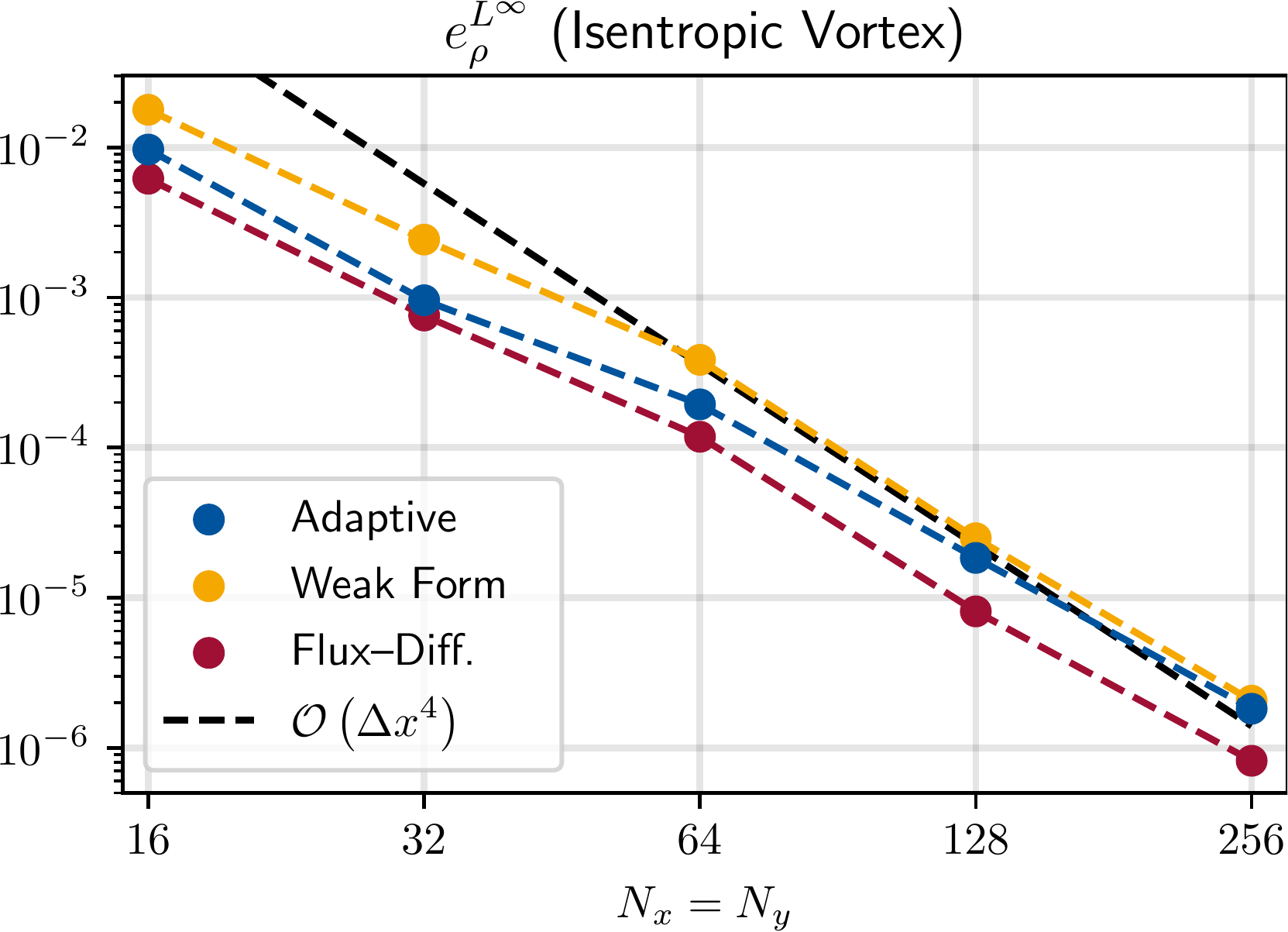}}
		}
	\caption[Errors of the isentropic vortex advection testcase with viscous surface flux.]
	{Errors of the isentropic vortex advection testcase with viscous surface flux.}
	\label{fig:IVA_Convergence_L2_LInf}
	\end{figure}

	We consider the isentropic vortex advection problem \cite{shu1988efficient, wang2013high} using the \ac{HLLC} interface flux \cite{toro1994restoration} and the entropy-conservative, kinetic energy-preserving volume flux described by Ranocha \cite{Ranocha2020Entropy}.
	The polynomial degree is set to $p = 3$ and the vortex setup follows the configuration from \cite{vermeire2019paired}.
	The simulation involves a single pass of the vortex through the computational domain.
	We consider four refinements of the computational domain, starting with $2^4 \times 2^4$ elements up to $2^8 \times 2^8$ elements.
	We consider the \ac{WF}, \ac{FD}, and the $v$--adaptive scheme with rigorous indicator from \cref{sec:indicator}.
	Domain--normalized $L^2$ error and pointwise $L^\infty$ error in density are depicted in \cref{fig:IVA_Convergence_L2_LInf}.
	The expected asymptotic convergence rates are observed for all schemes, while the lowest errors are obtained with the \ac{FD} volume term discretization.
	\subsubsection{Density Wave with Inviscid Surface Flux}
	As a first verification we consider the 1D density wave test case from \cite{gassner2022stability} given by
	\begin{subequations}
		\begin{align}
			\rho(x) &= 1 + 0.98 \sin \big(2 \pi x \big) \\
			v_x(x)  &= 0.1 \\
			p(x)    &= 20
		\end{align}
	\end{subequations}
	on $\Omega = [-1, 1]$ with periodic boundary conditions.
	The compressible Euler equations reduce for this setup to a linear advection equation for the density $\rho$ with velocity $v_x = 0.1$.
	For the purpose of this testcase, we employ the inviscid \ac{EC} two--point flux by Ranocha \cite{Ranocha2020Entropy} for both the surface flux and the volume flux in the \ac{FD} volume term discretization.
	The solution polynomials are discretized with polynomial degree $p = 3$ and the simulation is advanced up to $t_f = 2$ using the five--stage, fourth--order low--storage Runge--Kutta method from \cite{carpenter1994fourth}.

\begin{figure}
		\centering
		\subfloat[{Domain--normalized $L^2$ error in density.}]{
			\label{fig:DW_Convergence_L2}
			\resizebox{.46\textwidth}{!}{\includegraphics{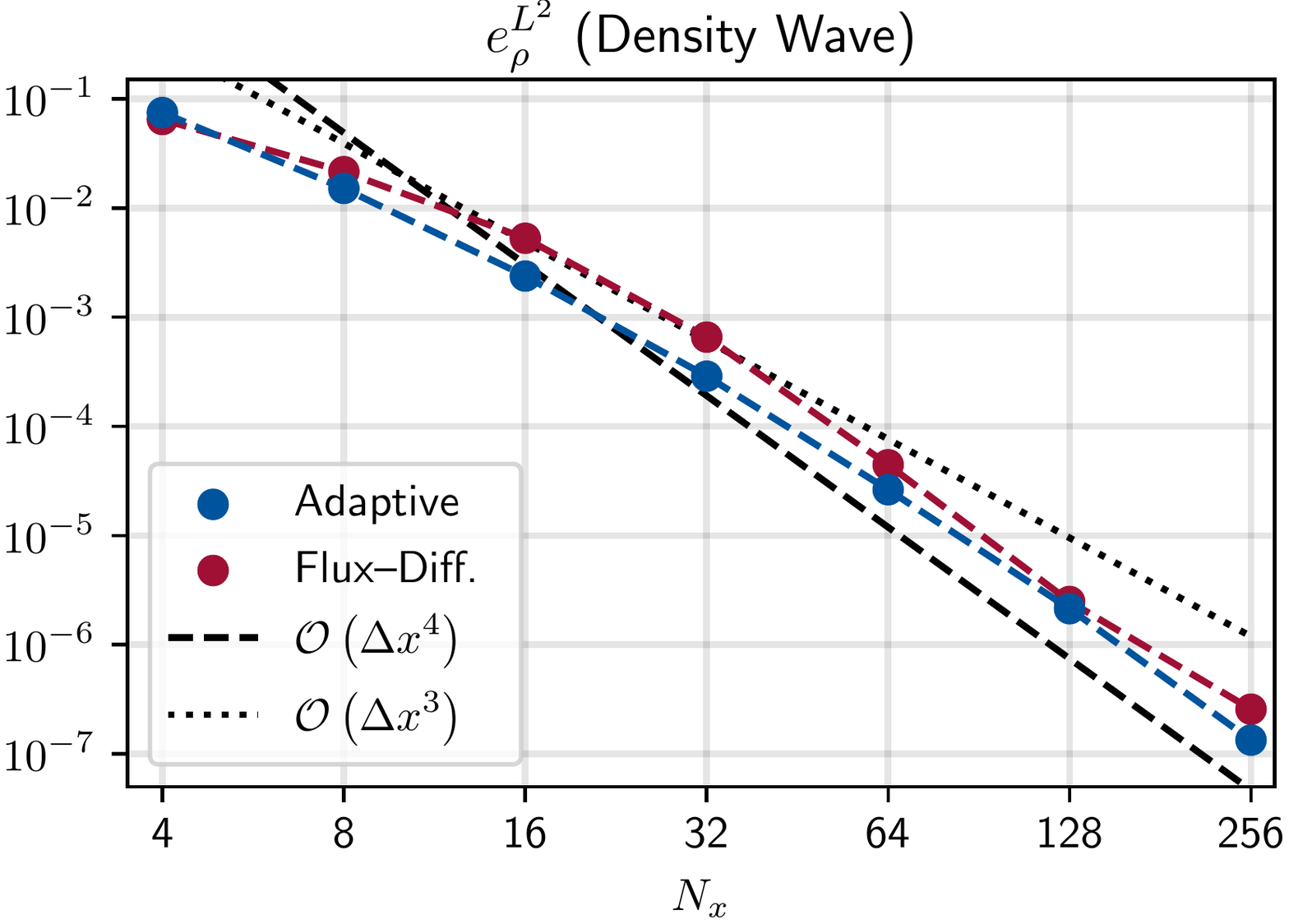}}
		}
		\hfill
		\subfloat[{Pointwise $L^\infty$ error in density.}]{
			\label{fig:DW_Convergence_LInf}
			\resizebox{.46\textwidth}{!}{\includegraphics{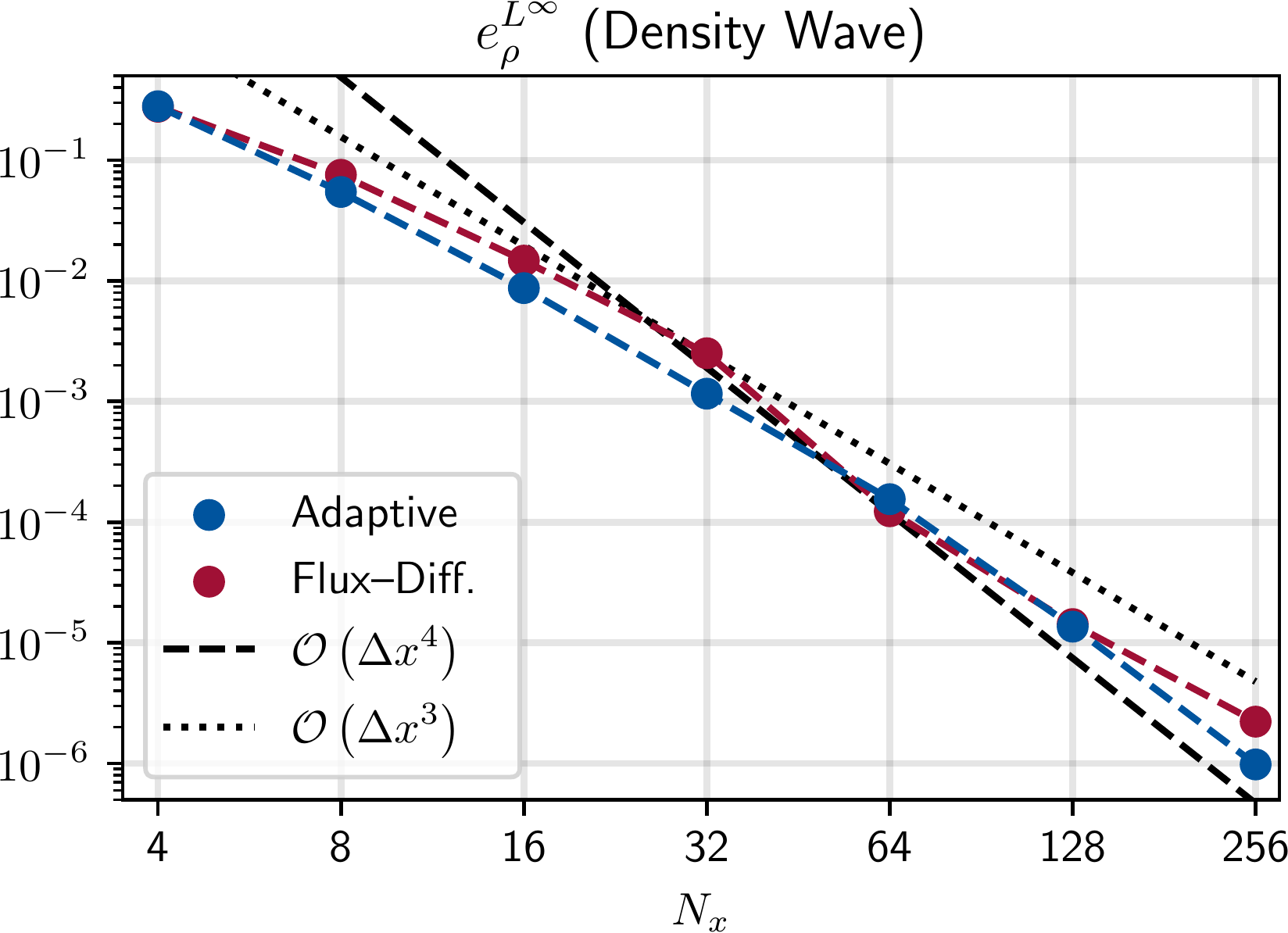}}
		}
	\caption[Errors for the density wave convergence testcase with inviscid surface flux.]
	{Errors for the density wave convergence testcase with inviscid surface flux.}
	\label{fig:DW_Convergence_L2_LInf}
	\end{figure}

	In contrast to the previous example, where the \ac{FD} discretization yielded the lowest errors, we observe here that the adaptive, entropy--dissipative discretization (i.e., with rigorous indicator) leads to lower errors for both $L^2$ and $L^\infty$ norms as depicted in \cref{fig:DW_Convergence_L2_LInf}.
	The \ac{WF} discretization is not depicted here as it is unstable for this configuration.
	Notably, only the adaptive scheme recovers fourth--order convergence in $L^\infty$ norm, while the pure \ac{EC} scheme only yields third--order convergence.
	\subsection{Linear Stability}
	In \cite{gassner2022stability} linear stability defects of entropy--conservative \ac{FD} volume term discretizations were analyzed.
	In particular, it was shown that for certain configurations of the \ac{DGSEM} the entropy--conservative \ac{FD} volume term leads to linearly unstable schemes, even when paired with dissipative interface fluxes.
	Here, we focus on the two--dimensional density wave example with configuration from \cite{gassner2022stability} given by
	\begin{subequations}
		\begin{align}
			\rho(x,y) &= 1 + 0.98 \sin \big(2 \pi (x + y) \big) \\
			v_x(x,y) &= 0.1 \\
			v_y(x,y) &= 0.2 \\
			p(x,y) &= 20
		\end{align}
	\end{subequations}
	on $\Omega = [-1, 1]^2$ with periodic boundary conditions.
	The compressible Euler equations reduce for this setup to a linear advection equation for the density $\rho$ with velocity $\boldsymbol v = (0.1, 0.2)^T$.
	The solution is discretized using $4 \times 4$ quadrilateral elements of polynomial degree $p = 5$ with the Rusanov/\ac{LLF} flux \cite{RUSANOV1962304} and the entropy--conservative volume flux from Chandrashekar \cite{Chandrashekar_2013}.
	This setup is deliberately chosen to be sensitive to perturbations in density, which pose the danger to violate positivity.
	The desired final time is $t_f = 5.0$, timestepping is performed with the explicit fourth--order, five--stage low--storage Runge--Kutta method from \cite{carpenter1994fourth}.
	The timestep is computed from the \ac{CFL} condition 
	\begin{equation}
		\Delta t = \text{CFL} \frac{1}{(p + 1)} \min_i \frac{\Delta x_i}{\lambda_i}
	\end{equation}
	where $\lambda$ denote the wave speeds $\lambda \in \sigma \left(\partial_{\boldsymbol u} \boldsymbol f \right)$ and $\Delta x$ is a suitable measure for the local element size.
	The minimum is taken over all elements and their interpolation points.

\begin{figure}
		\centering
		\subfloat[{Adaptive}]{
			\label{fig:LinearStab_t0_AV}
			\resizebox{.309\textwidth}{!}{\includegraphics{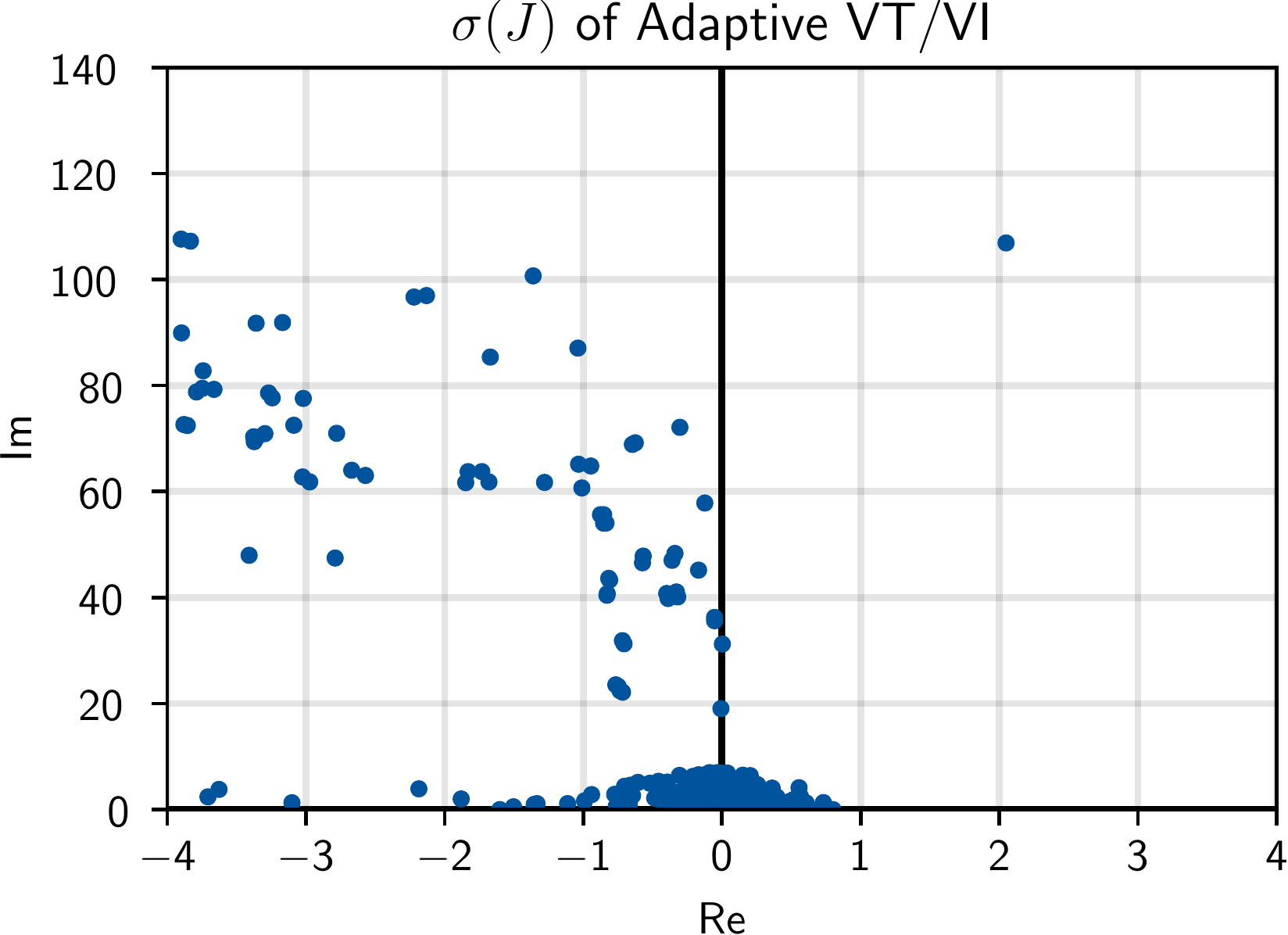}}
		}
		\subfloat[{\ac{WF}}]{
			\label{fig:LinearStab_t0_WF}
			\resizebox{.309\textwidth}{!}{\includegraphics{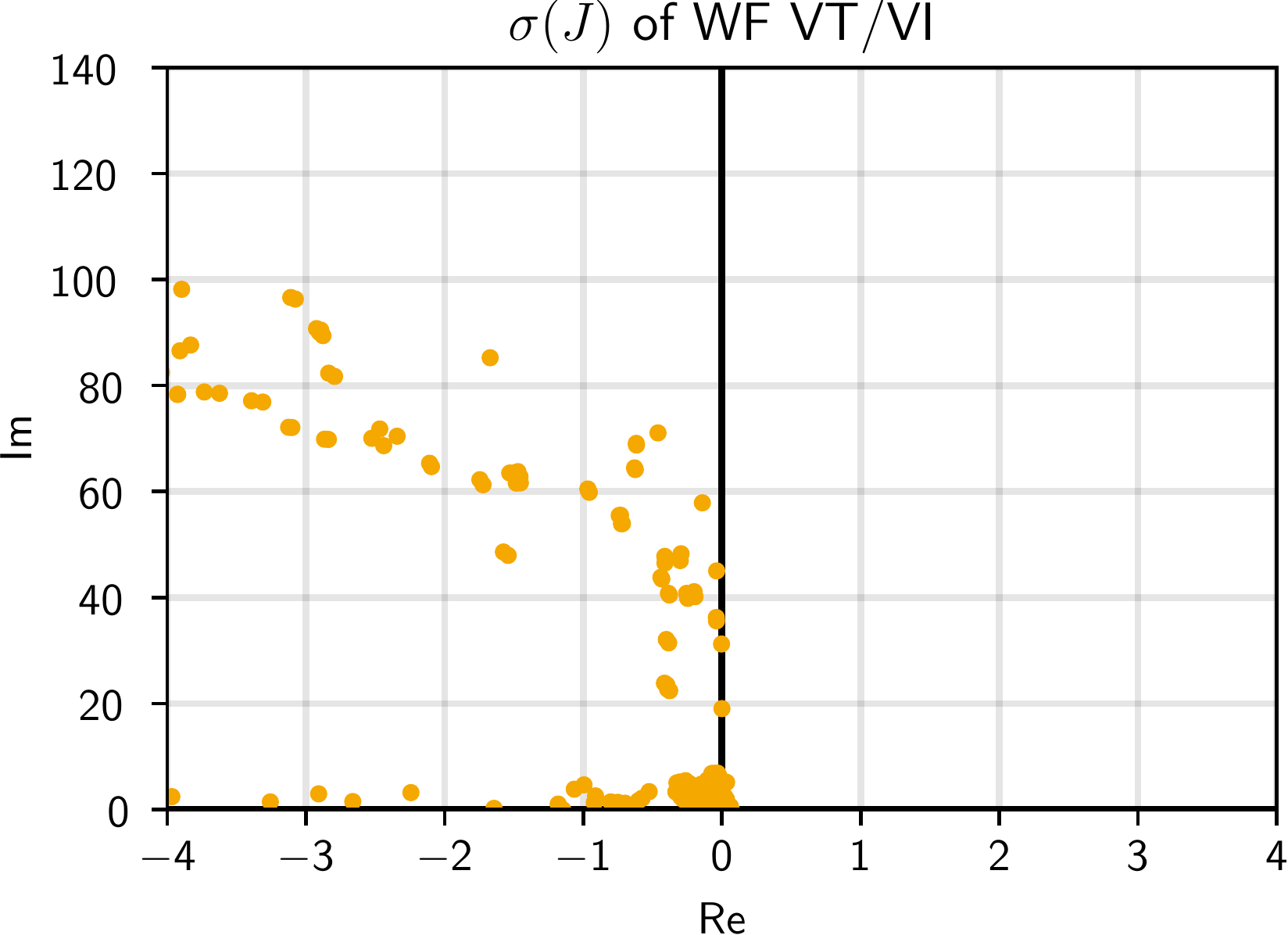}}
		}
		\subfloat[{\ac{FD}}]{
			\label{fig:LinearStab_t0_FD}
			\resizebox{.309\textwidth}{!}{\includegraphics{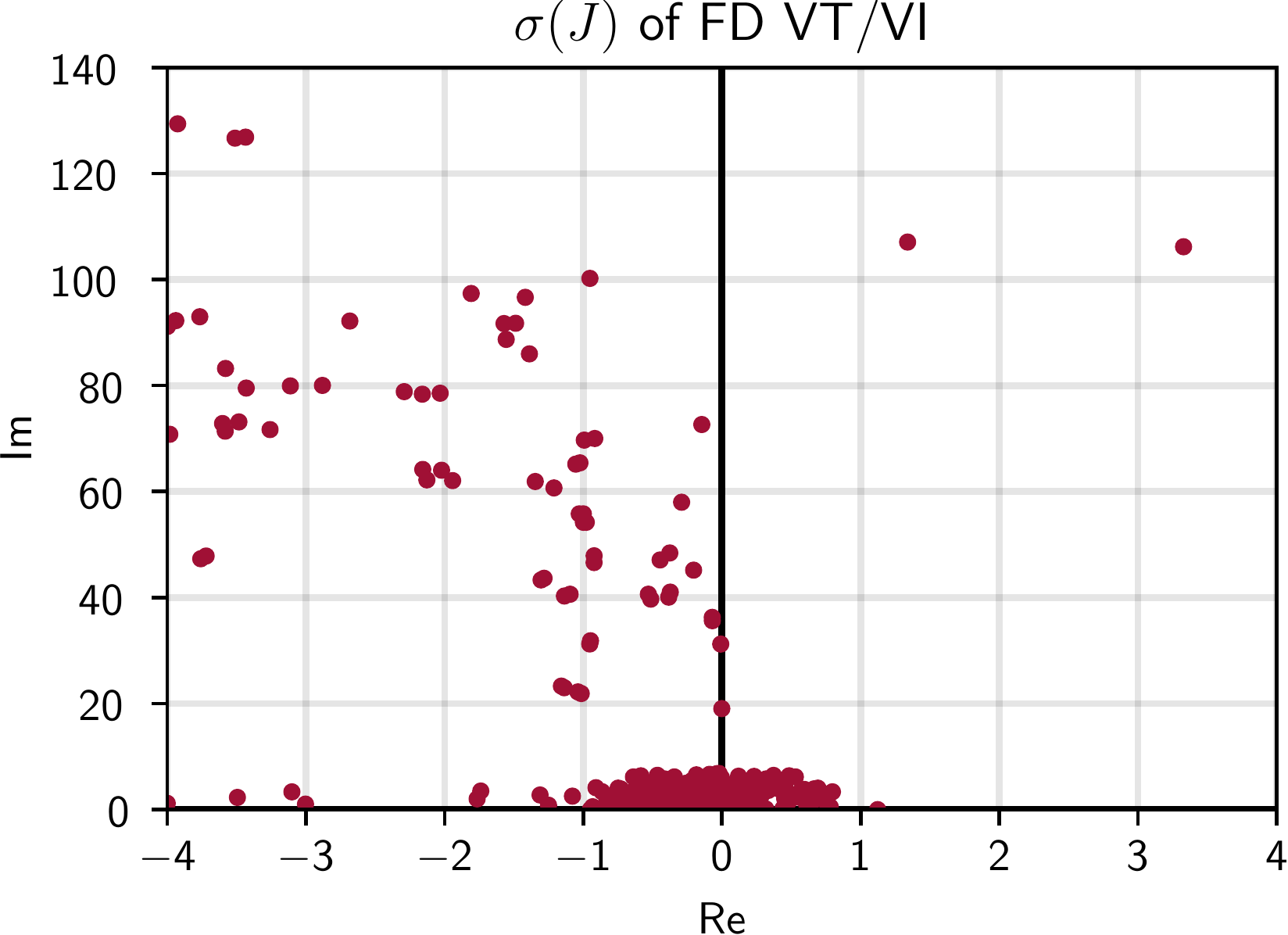}}
		}
	\caption[Spectra at initial time for the linear stability testcase.]
	{Spectra at initial time $t=0$ for the linear stability testcase.}
	\label{fig:Spectra_LinearStability_t0}
	\end{figure}

\begin{figure}
		\centering
		\resizebox{.46\textwidth}{!}{\includegraphics{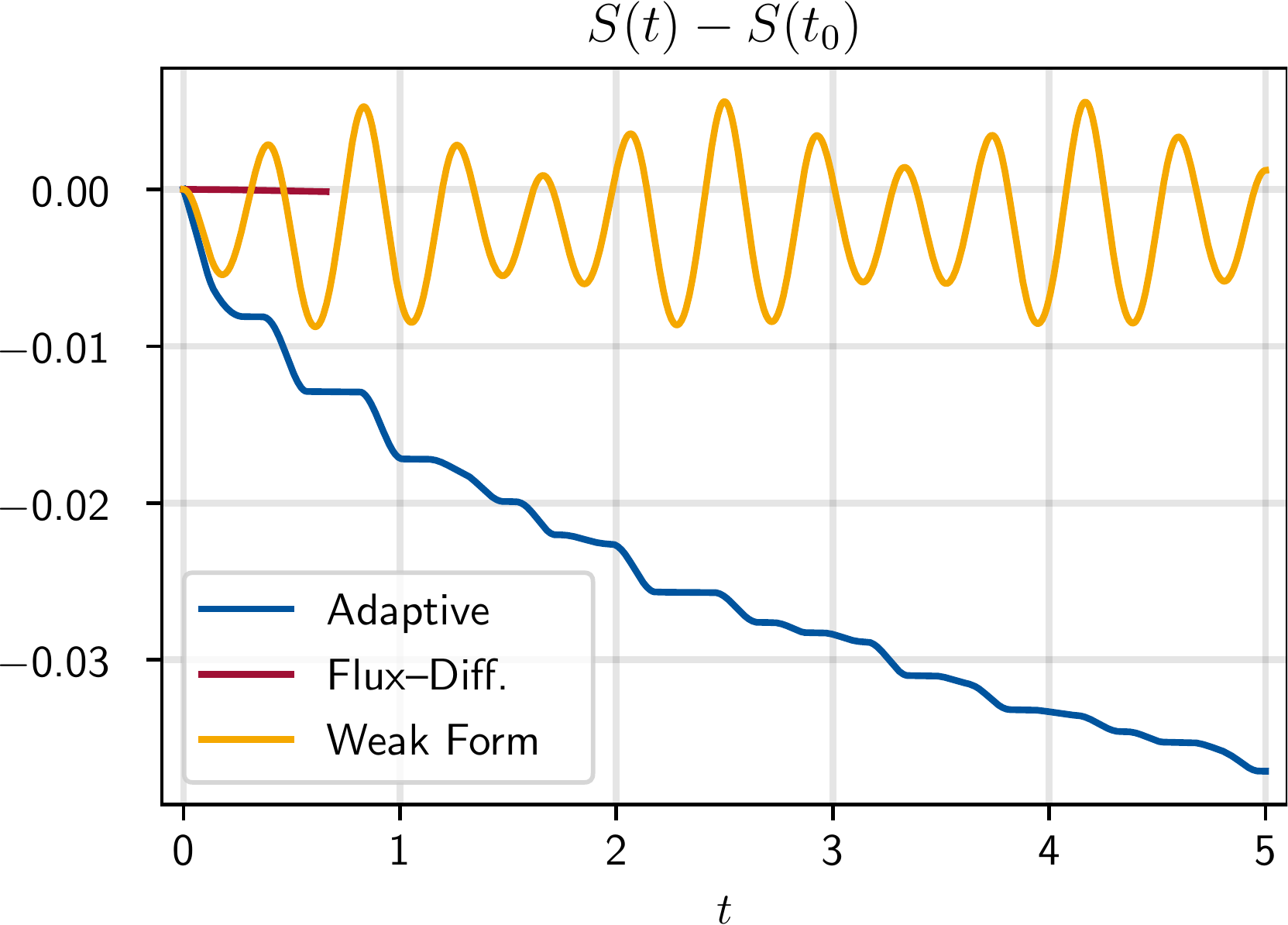}}
		\caption[]
		{Entropy difference over time for the density wave linear stability testcase with rigorous indicator for the adaptive volume term discretization.
		The simulations are performed with \ac{CFL} $=0.9$.}
		\label{fig:Entropy_LinearStability}
	\end{figure}

	As demonstrated in \cite{gassner2022stability}, the pure \ac{FD} volume term discretization is linearly unstable for this configuration, which has been confirmed both by spectral analysis and numerical experiments.
	We repeat this investigation here and include the \ac{WF} and the $v$--adaptive scheme with rigorous indicator from \cref{sec:indicator}.
	The required Jacobians are obtained exactly by employing algorithmic differentiation \cite{RevelsLubinPapamarkou2016} and the eigenvalue computation is performed using \texttt{Julia}'s built--in linear algebra package.
	The spectra at initial time $t=0$ are depicted in \cref{fig:Spectra_LinearStability_t0}.
	The most relevant metric for linear stability is the largest real part of the spectrum $r \coloneqq \max_{\lambda \in \sigma(J)} \text{Re}(\lambda)$.
	For the three considered schemes, we obtain the values $r_\text{AV} \approx 2.05$, $r_\text{WF} \approx 0.07$, and $r_\text{FD} \approx 3.32$ for the adaptive, weak form, and flux--differencing volume term discretizations, respectively.
	In other words, all discretizations are linearly unstable for this configuration, i.e., small perturbations will grow exponentially.
	In practice, however, we observe significantly different behaviour for the three schemes.
	As already reported in \cite{gassner2022stability}, the pure \ac{FD} volume term discretization crashes before reaching $t = 1$, while the \ac{WF} discretization is stable.
	Based on $r_\text{AV}$, one might infer that the adaptive scheme inherits the poor linear stability properties of the \ac{FD} volume term discretization.
	However, we are able to run the simulation up to the desired final time $t_f = 5.0$ (and way beyond) without loosing positivity in density.
	In fact, we conducted simulation runs up to $t_f = 5000$ corresponding to almost $1.2 \cdot 10^6$ timesteps with increased $\text{CFL} = 0.9$ without loss of positivity.
	Importantly, however, the \ac{WF} discretization does not imply entropy decay as depicted in \cref{fig:Entropy_LinearStability}.
	As evident from \cref{fig:Entropy_LinearStability}, the adaptive scheme dissipates entropy correctly, while maintaining stability.
	We recomputed the maximum real part of the spectrum for the \ac{WF} and adaptive scheme at $t_f = 5.0$ and obtain $r_\text{AV} \approx 0.55$, $r_\text{WF} \approx 0.07$, while we obtain for the \ac{FD} discretization shortly before the crash $r_\text{FD} \approx 1.42$.
	\subsection{Maintaining Entropy Solutions}
	\label{subsec:ConvEntropySolutions}
	A crucial necessity for the applicability of the $v$--adaptive scheme is that it maintains convergence to an entropy solution, even if the \ac{FD} volume term discretization is only applied on a subset of elements.
	We demonstrate this for the examples considered in \cite{lin2024high} and additionally the inviscid Burgers equation.
	\subsubsection{Burgers' Equation}
\begin{figure}
		\centering
		\subfloat[{\ac{WF} \ac{VT} discretization.}]{
			\label{fig:EntropySolution_Burgers_WF}
			\includegraphics[width=0.47\textwidth]{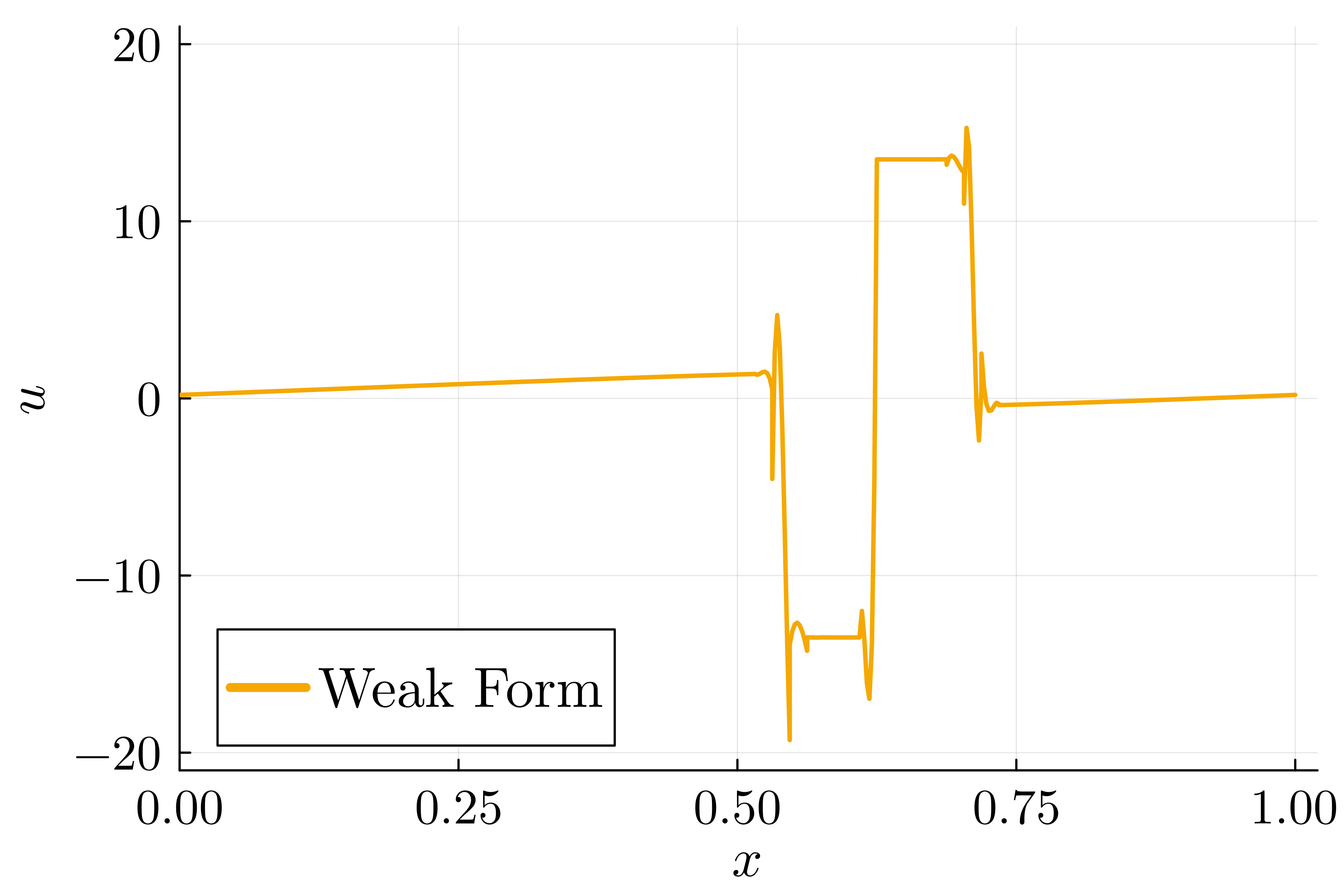}
		}
		\hfill
		\subfloat[{\ac{FD} and adaptive \ac{VT} discretization.}]{
			\label{fig:EntropySolution_Burgers_other}
			\includegraphics[width=0.47\textwidth]{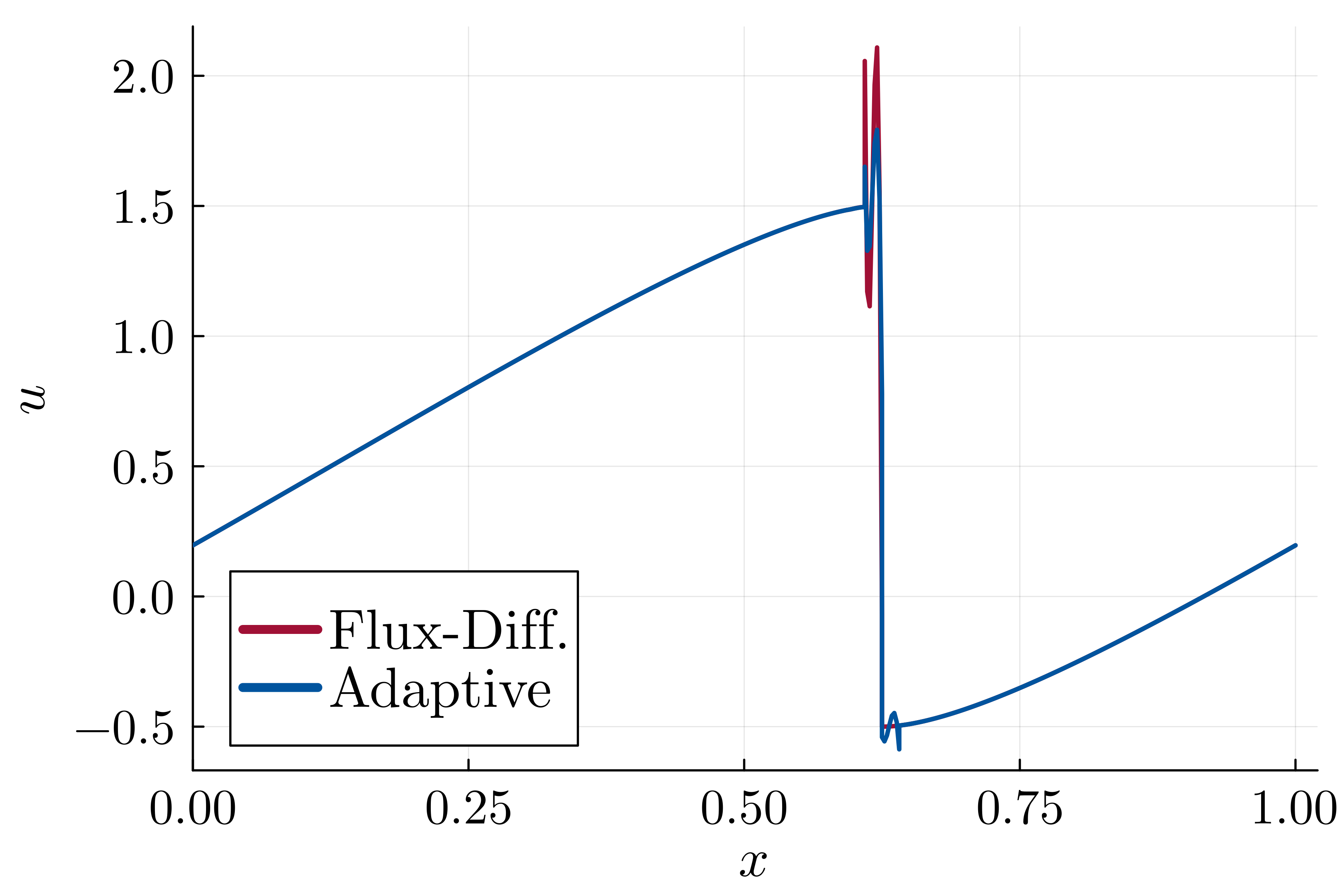}
		}
	\caption[Shock formation and propagation for Burgers' equation.]
	{Shock formation and propagation for Burgers' equation \cref{eq:BurgersEquation} at $t_f = 0.25$.
	Note the different scalings of the vertical axis.}
	\label{fig:EntropySolution_Burgers_WF_other}
	\end{figure}

	To begin, we consider the inviscid Burgers' equation
	\begin{equation}
		\label{eq:BurgersEquation}
		\partial_t u + \partial_x \left( \frac{u^2}{2} \right) = 0
	\end{equation}
	with entropy function, variable, flux and potential given by
	\begin{equation}
		S(u) = \frac{1}{2} u^2, \quad w = u, \quad F(u) = \frac{1}{3} u^3, \quad \psi(u) = \frac{1}{6} u^3 \: .
	\end{equation}
	We consider the initial condition
	\begin{equation}
		u_0(x) = \sin(2 \pi x) + 0.5
	\end{equation}
	on the periodic domain $\Omega = [0, 1]$.
	The solution develops a shock at time $t_s = \sfrac{1}{2 \pi}$ which then travels with constant speed $s = 0.5$.
	The solution is build from $p=3$ polynomials on $N = 64$ uniform elements.
	The surface fluxes are computed using the Godunov flux, i.e., an exact Riemann solver 
	is employed.
	Time integration is performed with the explicit fourth--order, five--stage low--storage \ac{SSP} method from \cite{ruuth2006global} implemented in the \texttt{Julia} \cite{bezanson2017julia} package \texttt{OrdinaryDiffEq.jl} \cite{DifferentialEquations.jl-2017} with CFL--based timestep selection.
	A plain weak form discretization will result in an unstable scheme due to uncontrolled oscillations.
	This is depicted in \cref{fig:EntropySolution_Burgers_WF} at $t_f = 0.25$.
	We see that the wild oscillations at the shock lead to formation of other completely artificial discontinuities, which eventually cause a blow up of the simulation.

	In contrast, both the pure \ac{FD} discretization with volume flux
	\begin{equation}
		f^\text{EC}(u_L, u_R) = \frac{1}{6} \left( u_L^2 + u_L u_R + u_R^2 \right)
	\end{equation}
	as well as the $v$--adaptive scheme with rigorous indicator from \cref{subsec:RigorousIndicator} are able to capture the correct shock propagation as depicted in \cref{fig:EntropySolution_Burgers_other}.
	As expected, oscillations arise also for these setups at the discontinuity.
	Importantly, both schemes maintain stability, even for long runs up to $t_f = 10$.
	Furthermore, the adaptive scheme produces slightly smaller oscillations than the pure \ac{FD} discretization, which is a direct consequence of the entropy--diffusive nature of the adaptive scheme.
	\subsubsection{1D Compressible Euler Equations}
	We consider the modified Sod shock tube problem proposed in \cite{toro2009riemann} with initial condition
	\begin{equation}
		\label{eq:ModifiedSodShockTube}
		\boldsymbol u(x) = 
		\begin{cases}
			\boldsymbol u_L, & x < 0.3 \\
			\boldsymbol u_R, & x \geq 0.3
		\end{cases}	
		, \quad
		\boldsymbol u_L =
		\begin{pmatrix}
			1.0 \\ 0.75 \\ 1.0
		\end{pmatrix}
		, \quad
		\boldsymbol u_R =
		\begin{pmatrix}
			0.125 \\ 0.0 \\ 0.1
		\end{pmatrix}
		\: .
	\end{equation}
	The domain is set to $\Omega = [0, 1]$ which we discretize with $N = 64$ elements.
	As in \cite{lin2024high} we employ $p=3$ polynomials and run the simulations until final time $t_f = 0.2$.
	The surface fluxes are computed with the \ac{HLLC} flux \cite{toro1994restoration} and the volume fluxes with the entropy--conservative and kinetic energy preserving flux from Ranocha \cite{Ranocha2020Entropy} for the \ac{FD} and adaptive \ac{VT} discretizations.
	Time integration is performed with the adaptive explicit four/three--stage embedded \ac{SSP} Runge--Kutta method from \cite{fekete2022embedded}, implemented in \texttt{OrdinaryDiffEq.jl} \cite{DifferentialEquations.jl-2017}.

\begin{figure}
		\centering
		\subfloat[{\ac{WF} \ac{VT} discretization.}]{
			\label{fig:EntropySolution_CEE_WF}
			\includegraphics[width=0.47\textwidth]{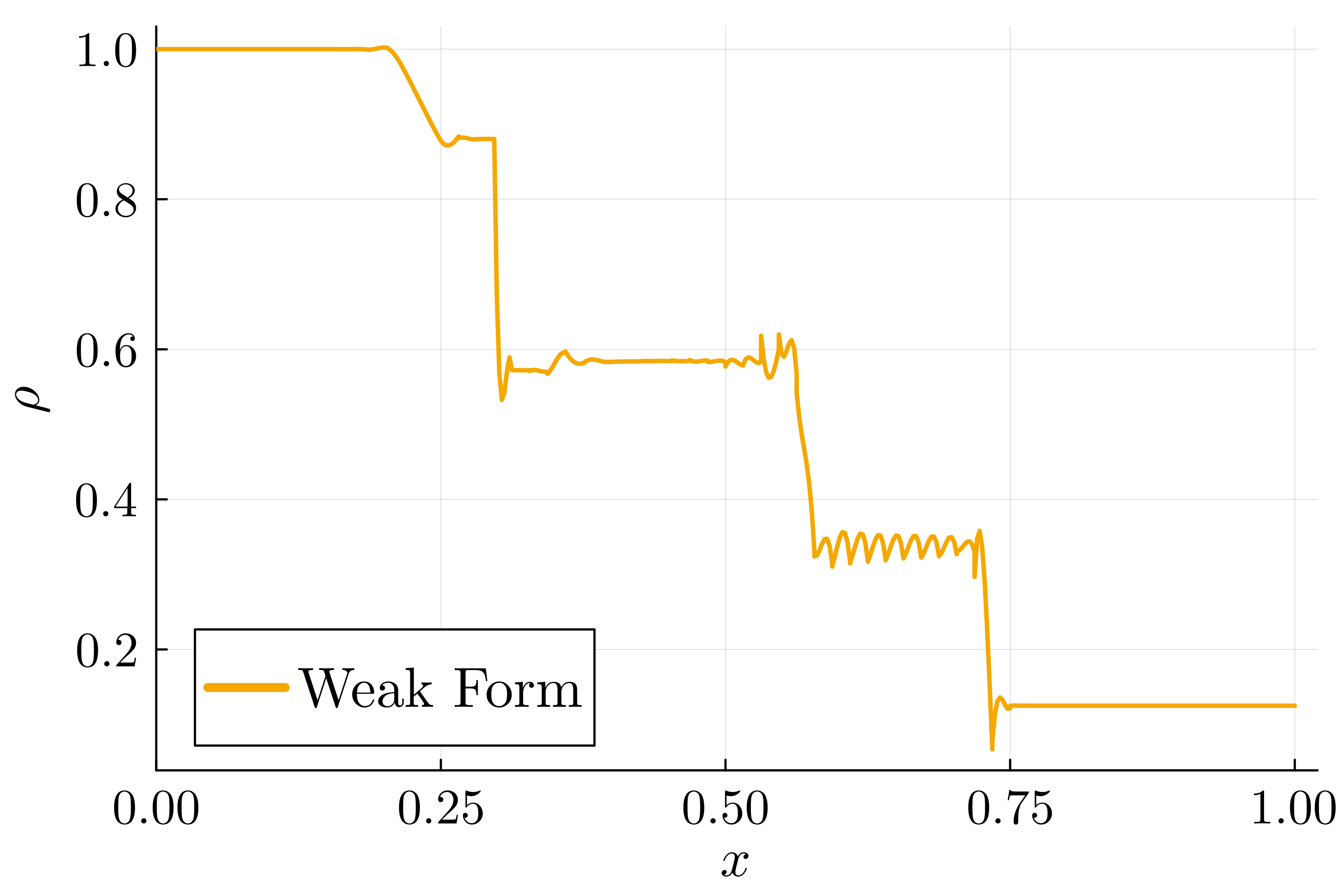}
		}
		\hfill
		\subfloat[{\ac{FD} and adaptive \ac{VT} discretization.}]{
			\label{fig:EntropySolution_CEE_other}
			\includegraphics[width=0.47\textwidth]{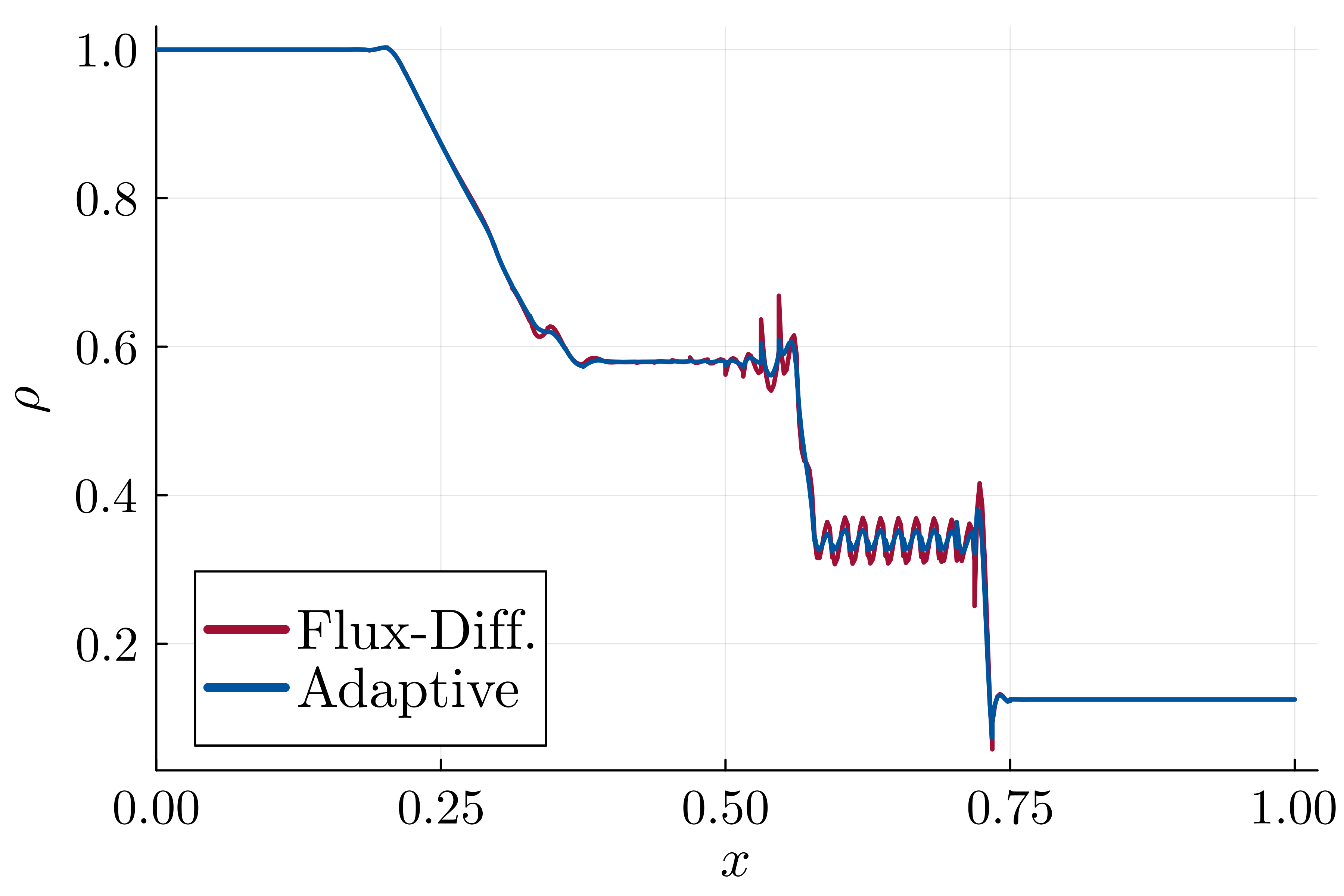}
		}
	\caption[Density $\rho$ at $t_f = 0.2$ for the modified Sod shock tube problem.]
	{Density $\rho$ at $t_f = 0.2$ for the modified Sod shock tube problem \cref{eq:ModifiedSodShockTube}.}
	\label{fig:EntropySolution_CEE_WF_other}
	\end{figure}

	As described also in \cite{lin2024high}, the pure \ac{WF} volume term discretization does not capture the smooth rarefaction wave correctly, cf. \cref{fig:EntropySolution_CEE_WF}.
	Furthermore, the \ac{WF} volume term will violate positivity of density and pressure if not explicitly enforced.
	To be able to run the simulation until final time $t_f = 0.2$, we therefore apply the positivity limiter from \cite{zhang2010positivity} with limiting threshold $5 \cdot 10^{-6}$ for both density and pressure.
	The positivity limiter can ensure positivity of density and pressure if the corresponding cell mean values are still positive.
	In contrast, both the pure \ac{FD} and the $v$--adaptive scheme with rigorous indicator from \cref{subsec:RigorousIndicator} are able to capture the correct rarefaction wave as depicted in \cref{fig:EntropySolution_CEE_other}.
	For both, the simulations preserve positivity of density and pressure without the need for a limiter.
	As for Burgers' equation, the adaptive scheme produces slightly smaller oscillations than the pure \ac{FD} discretization due to the selection of \ac{WF} integral if it introduces additional dissipation.
	\subsubsection{Kurganov--Petrova--Popov (KPP) Equation}
	\label{subsubsec:KPP}
\begin{figure}
		\centering
		\subfloat[{\ac{WF} discretization, realized in \ac{FD} form with central volume flux.}]{
			\label{fig:KPP_WF}
			\includegraphics[width=0.47\textwidth]{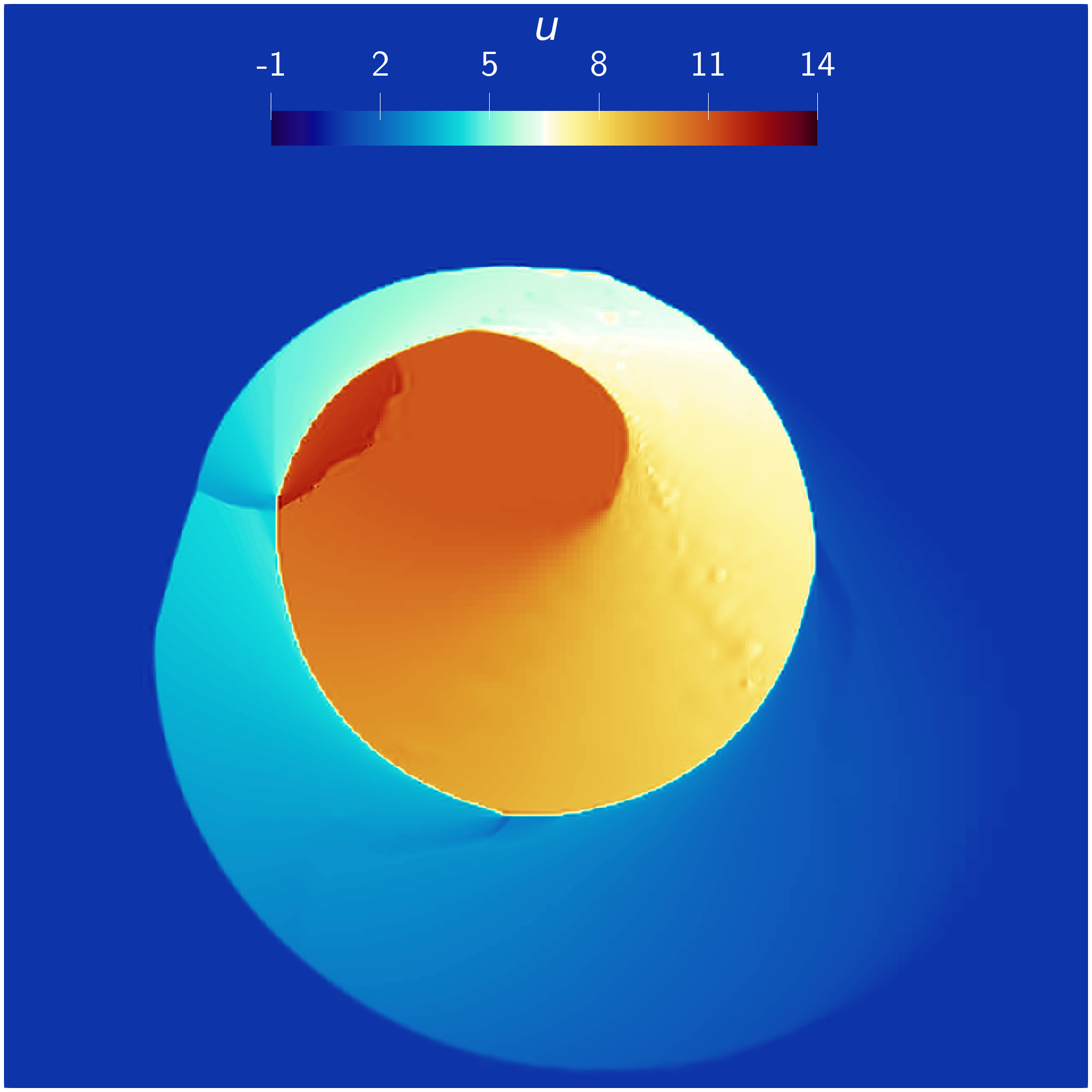}
		}
		\hfill
		\subfloat[{Blended \ac{FD}--\ac{FV} discretization with adaptive volume flux selection.}]{
			\label{fig:KPP_VTA}
			\includegraphics[width=0.47\textwidth]{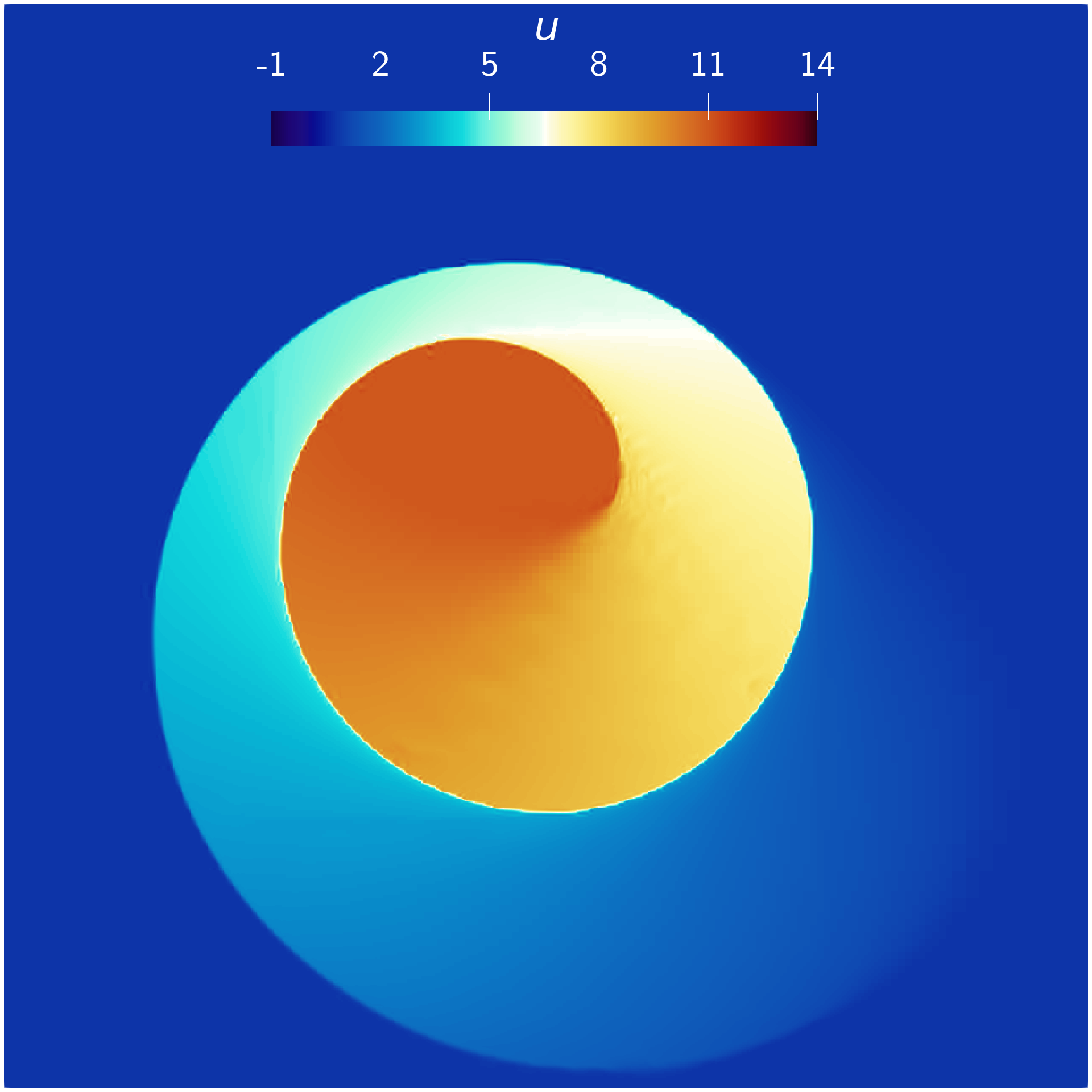}
		}
		\\
		\subfloat[{Blending parameter $\beta$.}]{
			\label{fig:KPP_beta}
			\includegraphics[width=0.47\textwidth]{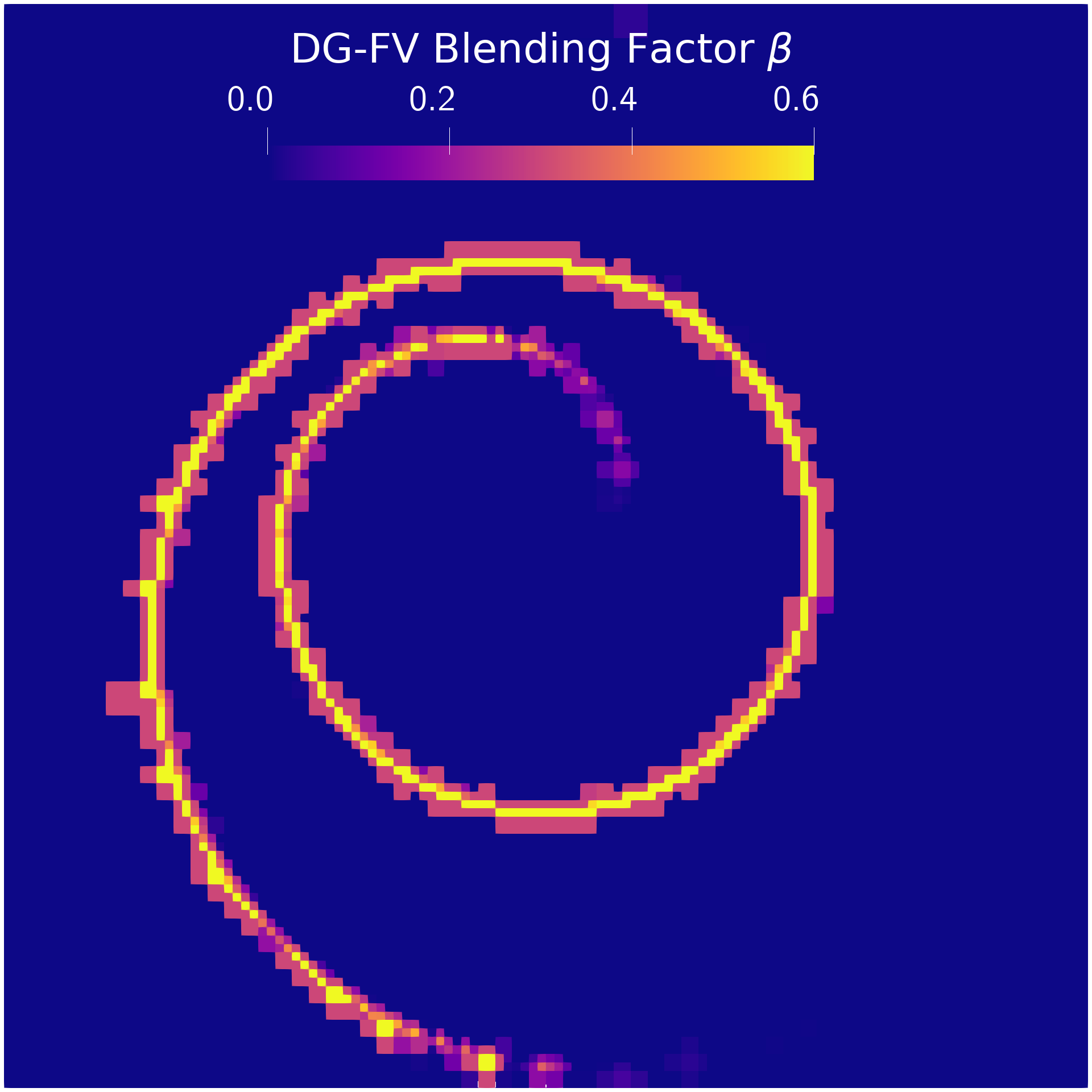}
		}
		\subfloat[{Adaptive mesh.}]{
			\label{fig:KPP_mesh}
			\includegraphics[width=0.47\textwidth]{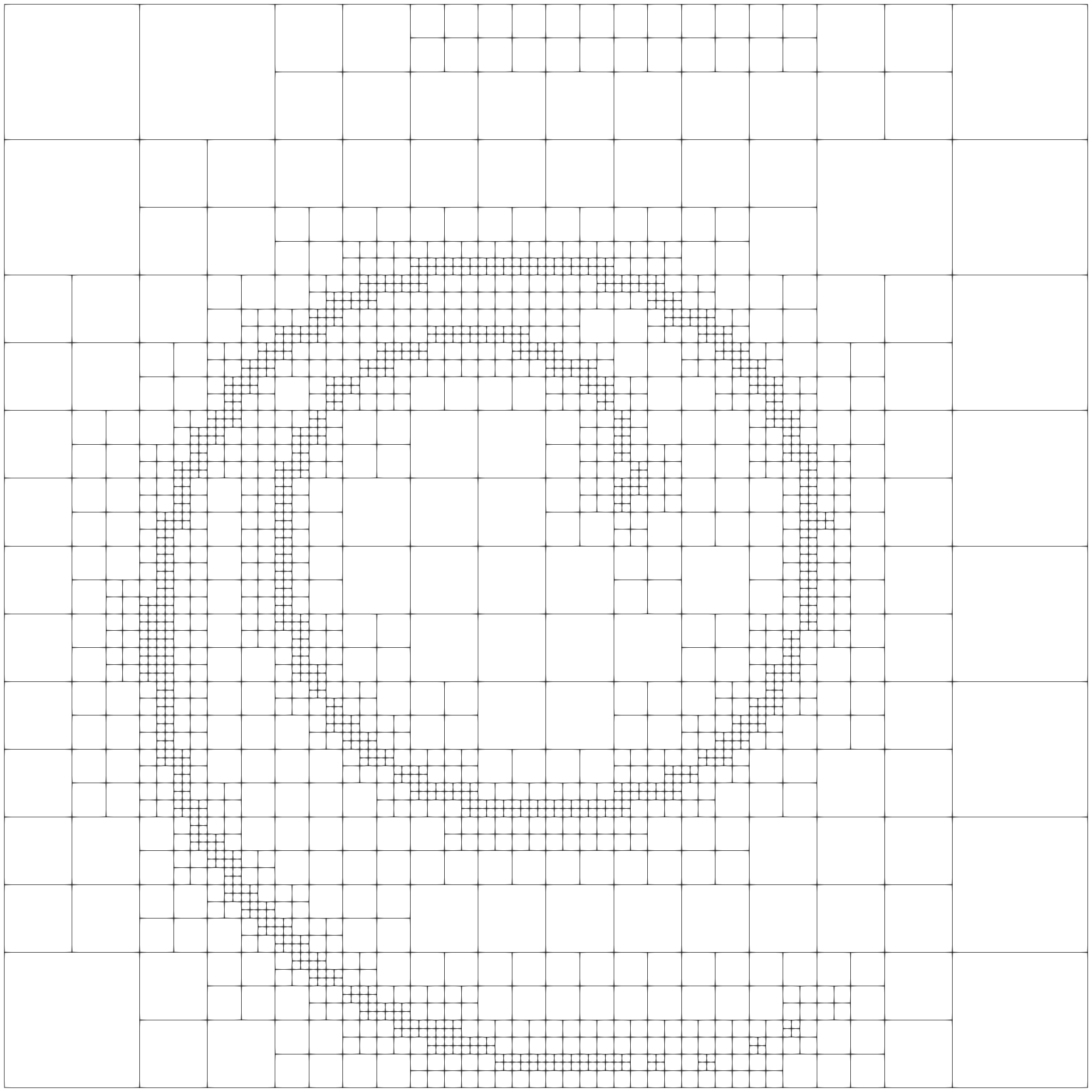}
		}
	\caption[KPP problem. Solution for naive central flux and entropy--conservative flux.
	Stabilization parameter $\beta$ and adaptive mesh.]
	{Solution for the KPP problem \cref{eq:KPP}, \cref{eq:KPP_IC} at $t_f = 1.0$ obtained with the blended \ac{FD}--\ac{FV} scheme.
	The solution depicted in \cref{fig:KPP_WF} is obtained with the naive central flux, while the solution in \cref{fig:KPP_VTA} is obtained by selecting the entropy--conserving flux \cref{eq:KPP_EC_Flux} in cells where the shock--capturing parameter $\beta_i$ is non--zero, cf. \cref{fig:KPP_beta}.
	The adaptive mesh is depicted in \cref{fig:KPP_mesh}.}
	\label{fig:KPP_WF_VTA_beta_mesh}
	\end{figure}

	Finally, we consider the two--dimensional Kurganov--Petrova--Popov (KPP) problem \cite{kurganov2007adaptive} which is a scalar conservation law
	\begin{equation}
		\label{eq:KPP}
		\partial_t u + \partial_x \sin(u) + \partial_y \cos(u) = 0
	\end{equation}
	with non--convex flux functions in two spatial dimensions.
	The discontinuous initial condition is given by
	\begin{equation}
		\label{eq:KPP_IC}
		u_0(x,y) = 
		\begin{cases}
			3.5 \pi & x^2 + y^2 < 1 \\
			0.25\pi & x^2 + y^2 \geq 1
		\end{cases}	
		\: .
	\end{equation}
	The entropy is given by $S(u) = 0.5 u^2$ and the entropy--conserving fluxes per direction are given by
	\begin{subequations}
		\label{eq:KPP_EC_Flux}
		\begin{align}
			f^\text{EC}_x(u_L, u_R) &= \frac{\cos(u_L) - \cos(u_R)}{u_R - u_L} \\
			f^\text{EC}_y(u_L, u_R) &= \frac{\sin(u_R) - \sin(u_L)}{u_R - u_L} \: .
		\end{align}
	\end{subequations}
	If $\vert u_R - u_L \vert \leq 10^{-12}$ we fall back to the central flux $0.5 \left( f(u_L) + f(u_R) \right)$.
	As pointed out in \cite{kurganov2007adaptive} high--order methods may fail to converge to the correct entropy solution.
	We use this example to test now the combination of the \ac{WF} volume integral with the subcell--stabilized \ac{FD}--\ac{FV} volume integral from \cite{hennemann2021provably, rueda2021entropy}.
	The KPP problem is simulated on $\Omega = [-2, 2]^2$ with periodic boundary conditions.
	The solution is discretized using $p=4$ polynomials, and we run the simulations until final time $t_f = 1.0$ using the five--stage, fourth--order \ac{SSP} Runge--Kutta method from \cite{ruuth2006global} with constant timestep $\Delta t = 10^{-3}$.
	We employ an adaptive mesh with $6$ levels, starting from base resolution $4 \times 4$, i.e., $\Delta x_\text{max} = 1.0$ down to $\Delta x_\text{min} = 0.03125$.
	The surface fluxes are discretized with the Rusanov/\ac{LLF} flux \cite{RUSANOV1962304} and the volume fluxes with the entropy--conservative fluxes \cref{eq:KPP_EC_Flux} provided above.
	To be able to use the shock--capturing framework from \cite{hennemann2021provably, rueda2021entropy}, we realize the \ac{WF} volume integral in the \ac{DG}--\ac{FV} blending in \ac{FD} form with central flux.
	This leads to a wrong solution, as depicted in \cref{fig:KPP_WF}.
	In contrast, the adaptive scheme which selects the entropy--conserving flux \cref{eq:KPP_EC_Flux} in cells where the shock--capturing parameter $\beta_i$ is non--zero is able to capture the correct solution, cf. \cref{fig:KPP_VTA}.
	Importantly, the usage of the \ac{WF}/central flux \ac{FD} outside the shock--capturing regions does not lead to a wrong solution.
	\section{Numerical Examples}
	\label{sec:Examples}
	We present a range of applications of the $v$--adaptive scheme to compressible flows governed by the Euler equations and the Navier--Stokes equations.
	As in the previous section, we begin by coupling the weak form \ac{VT/VI} with the flux--differencing \ac{VT/VI}.
	Then, we consider the combination of the weak form volume integral with the blended \ac{DG}--\ac{FV} scheme from \cite{hennemann2021provably, rueda2021entropy}.
	The numerical experiments can be reproduced without any commercial or proprietary software, with files available on \texttt{GitHub} and archived on \texttt{Zenodo} \cite{doehring2026VTA_ReproRepo}.
	\subsection{Weak Form --- Flux--Differencing Combination}
	For the first numerical examples we consider the combination of the weak form and flux--differencing volume term discretizations with the entropy production indicator introduced in \cref{sec:indicator}.
	Thus, we strictly maintain order of accuracy in all cells, but obtain a potentially more stable and efficient scheme.
	\subsubsection{Kelvin--Helmholtz Instability with Entropy--Diffusive Scheme}
	\label{subsubsec:KHI_DGSEM_Rigorous}
	We consider the Kelvin--Helmholtz instability with parametrization from \cite{rueda2021subcell} which is governed by the compressible Euler equations with $\gamma = 1.4$ in two spatial dimensions.
	The initial condition is given by
	\begin{equation}
		\renewcommand\arraystretch{1.15}
		\begin{pmatrix}
			\rho \\ v_x \\ v_y \\ p
		\end{pmatrix}
		= 
		\renewcommand\arraystretch{1.15}
		\begin{pmatrix}
			\frac{1}{2} + \frac{3}{4} b \\
			\frac{1}{2} (b - 1) \\
			\frac{1}{10} \sin(2 \pi x) \\
			1
		\end{pmatrix}, 
		\quad b \coloneqq \tanh \left( 15y + 7.5 \right) - \tanh \left(15y - 7.5 \right)
	\end{equation}
	on periodic domain $\Omega = [-1, 1]^2$.
	This simulation features an always under--resolved flow due to the fractal--like formation of small--scale flow features.
	The initial condition provides a Mach number $\text{Ma} \leq 0.6$ which renders the flow compressible, but does not cause shocks to develop in the domain.
	The interface fluxes are computed with the \ac{HLLC} \cite{toro1994restoration} flux and the volume fluxes with the entropy--conservative and kinetic energy preserving flux from Ranocha \cite{Ranocha2020Entropy}.
	Time integration is performed with the explicit fourth--order, five--stage low--storage Runge--Kutta method from Carpenter and Kennedy \cite{carpenter1994fourth} with \ac{CFL} based timestepping implemented in \texttt{OrdinaryDiffEq.jl} \cite{DifferentialEquations.jl-2017}.

\begin{figure}
		\centering
		\resizebox{.46\textwidth}{!}{\includegraphics{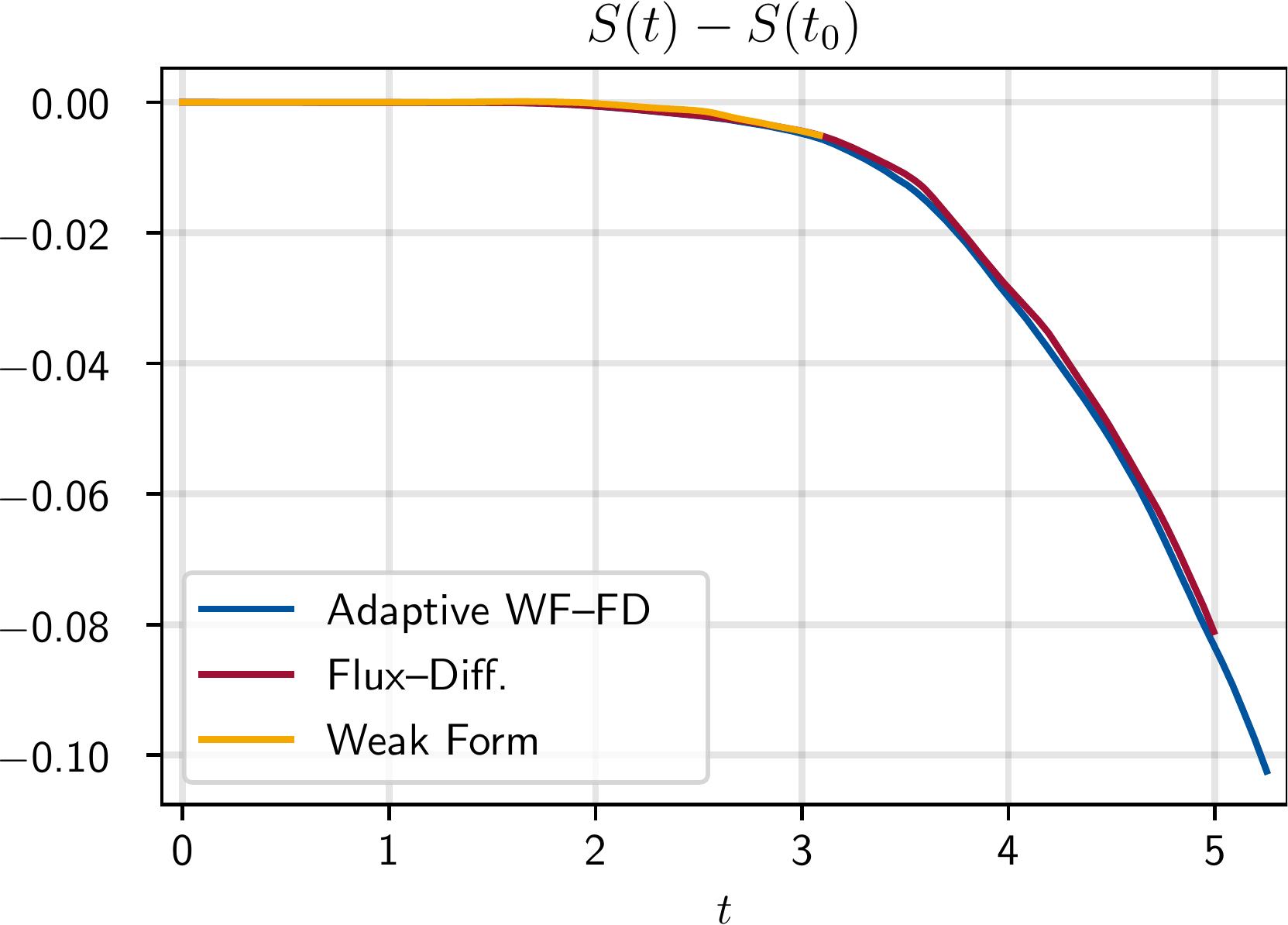}}
		\caption[]
		{Entropy difference over time for different \ac{VT} discretizations for the Kelvin--Helmholtz instability testcase on a $64 \times 64$ quadrilateral mesh with $p = 3$.}
		\label{fig:KHI_Rigorous_DGSEM_Entropy}
	\end{figure}

\begin{figure}
		\centering
		\resizebox{.55\textwidth}{!}{\includegraphics{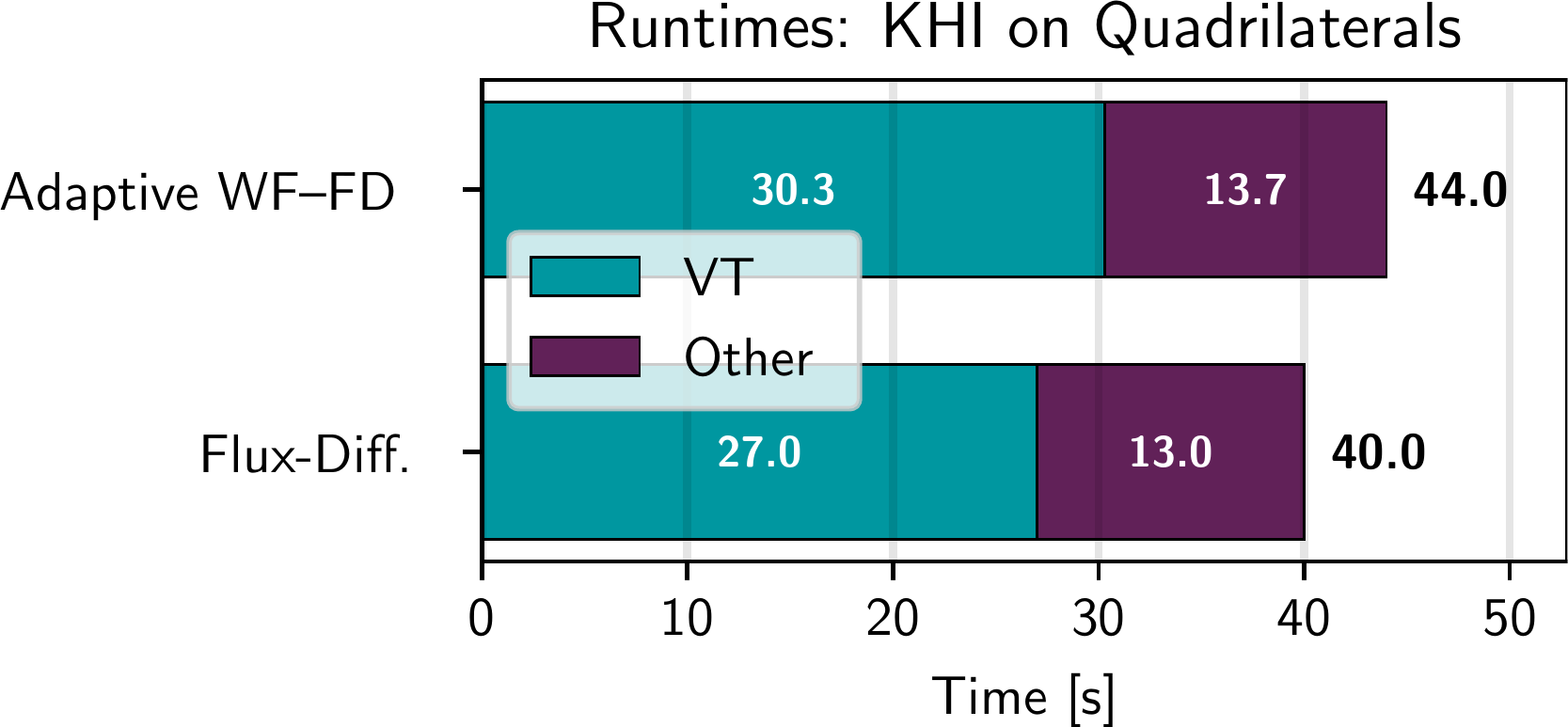}}
		\caption[Total RHS and volume term/integral computation runtime for the Kelvin--Helmholtz instability on quadrilaterals.]
		{Total \ac{RHS} and \ac{VT} computation runtime for the Kelvin--Helmholtz instability for $t_f = 4.95$.
		Volume term adaptivity is governed by the rigorous indicator resulting in an entropy--diffusive scheme.}
		\label{fig:RunTimes_KHI_Rigorous_DGSEM}
	\end{figure}

	To begin, we consider the nodal \ac{DGSEM} on $p = 3$ uniform $64 \times 64$ quadrilateral elements.
	We investigate how long different \ac{VT/VI} discretizations are able to maintain stability, i.e., positivity of density and pressure.
	The pure \ac{WF} volume term discretization is stable until around $t \approx 3.1$, the pure \ac{FD} volume term discretization is stable until around $t \approx 5.0$, and the entropy--diffusive $v$--adaptive scheme with rigorous indicator is stable until around $t \approx 5.33$.
	The difference in mathematical entropy is depicted in \cref{fig:KHI_Rigorous_DGSEM_Entropy} for the three simulations.
	The density at $t_f = 5.25$ is plotted in \cref{fig:KHI_Rigorous_DGSEM_rho} (see the Appendix) for the simulation with the entropy--diffusive $v$--adaptive scheme with rigorous indicator.
	We compare the computational cost of the different \ac{VT} discretizations by recording the runtime spent on the volume term computation for the pure \ac{FD} and the $v$--adaptive scheme with rigorous indicator on the $t \in [0, 4.95]$ interval, i.e., while the pure \ac{FD} scheme is still stable.
	We observe that the adaptive scheme is about $11\%$ more expensive than the pure \ac{FD} scheme, cf. \cref{fig:RunTimes_KHI_Rigorous_DGSEM}.
	This is due the fact that at later simulation times the \ac{FD} \ac{VT} is applied on a significant portion of the elements, rendering the \ac{WF} \ac{VT} computation performed for the indicator in principle redundant.
	Nevertheless, even at later simulation times some elements dissipate entropy through the \ac{WF} \ac{VT} which allows the simulation to advance further in time than the pure \ac{FD} scheme.
	We observe small discrepancies in runtime for the non \ac{VT} ("other") computations, most notably interface flux and surface integral computations which we cannot fully explain.
	The values reported here are averaged over five simulation runs.
	\subsubsection{Kelvin--Helmholtz Instability with Tolerated Entropy Increase}
	\label{subsubsec:KHI_DGMulti_Heuristic}
\begin{figure}
		\centering
		\resizebox{.46\textwidth}{!}{\includegraphics{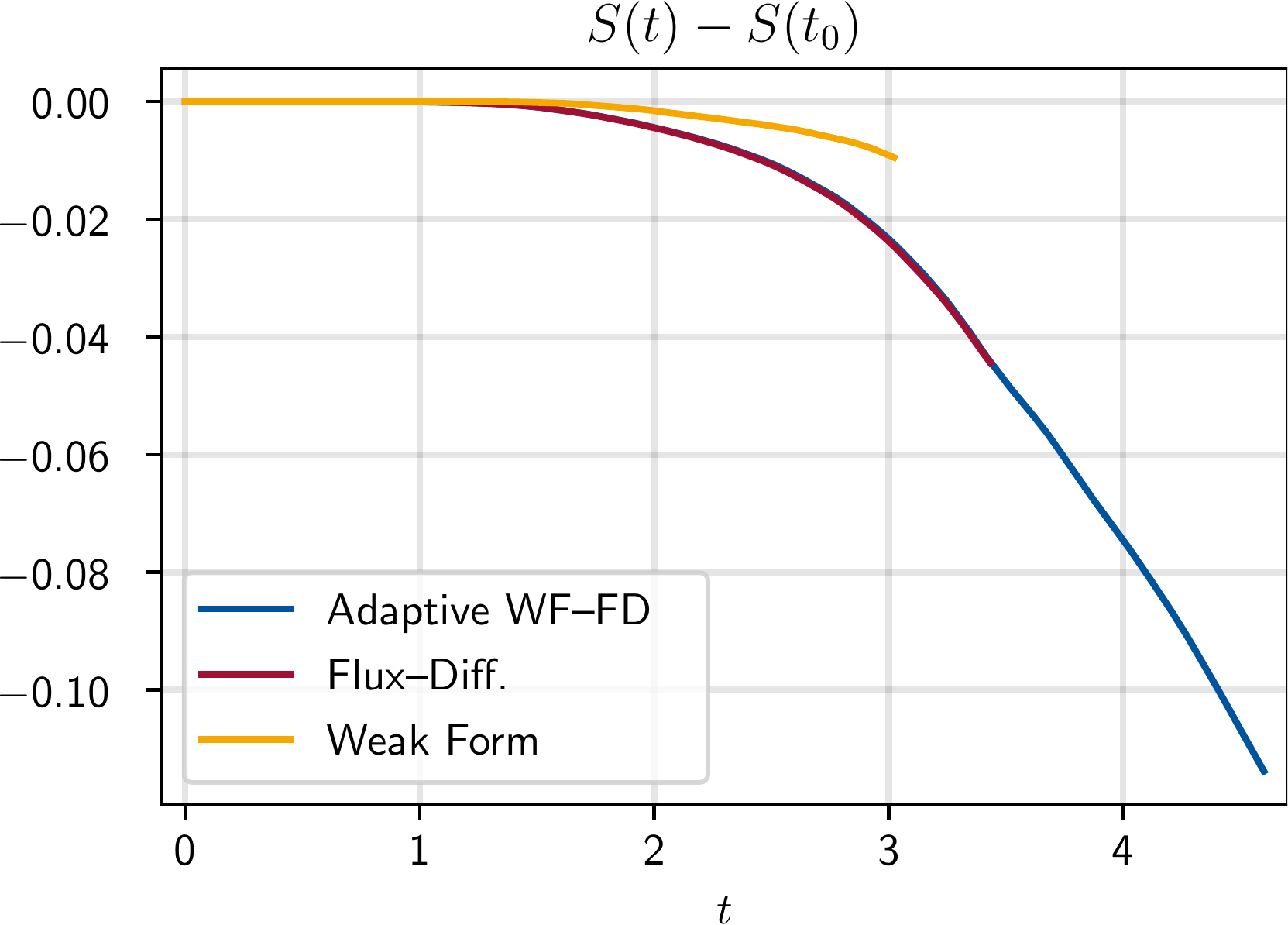}}
		\caption[]
		{Entropy difference over time for different \ac{VT} discretizations for the Kelvin--Helmholtz instability testcase on a $32 \times 32$ triangular mesh with $p = 3$. For the $v$--adaptive scheme, the heuristic indicator with tolerated entropy increase $\sigma_i = \sigma = 5 \cdot 10^{-3}$ is employed.}
		\label{fig:KHI_Heuristic_DGMulti_Entropy}
	\end{figure}

\begin{figure}
		\centering
		\resizebox{.55\textwidth}{!}{\includegraphics{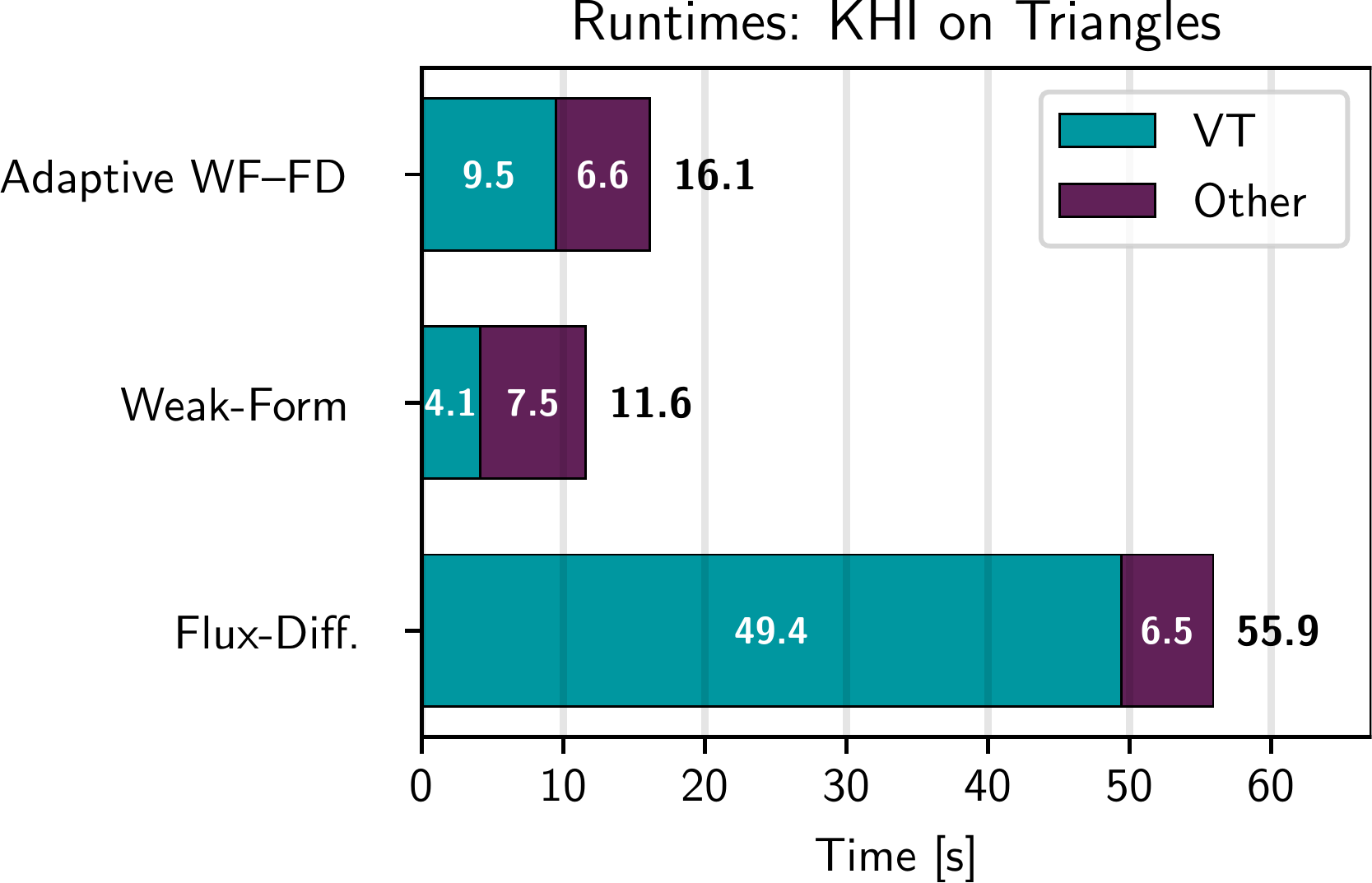}}
		\caption[Total RHS and volume term/integral computation runtime for the Kelvin--Helmholtz instability on triangles.]
		{Total \ac{RHS} and \ac{VT} computation runtime for the Kelvin--Helmholtz instability for $t_f = 3.0$.
		Volume term adaptivity is governed by the heuristic indicator with tolerated entropy increase $\sigma_i = \sigma = 5 \cdot 10^{-3}$.}
		\label{fig:RunTimes_KHI_Heuristic_DGMulti}
	\end{figure}

	We repeat the Kelvin--Helmholtz instability simulation on triangles with the heuristic indicator, recall \cref{eq:EntropyTargetDecay} with tolerated entropy increase $\sigma_i = \sigma = 5 \cdot 10^{-3}$.
	We discretize the domain with $32 \times 32$ triangles with symmetric Gauss--Jacobi quadrature \cite{wandzurat2003symmetric, xiao2010numerical} of order $2p$ and solution polynomial degree $p = 3$.
	These are implemented in the \texttt{Julia} packages \texttt{StartUpDG.jl} \cite{StartUpDG} and \texttt{NodesandModes.jl} \cite{NodesandModes}.
	The surface fluxes are computed with the \ac{HLLC} flux \cite{toro1994restoration} and the volume fluxes with the entropy--conservative and kinetic energy preserving flux from Ranocha \cite{Ranocha2020Entropy}.
	Time integration is performed with the explicit fourth--order, five--stage low--storage Runge--Kutta method from Carpenter and Kennedy \cite{carpenter1994fourth} with \ac{CFL} based timestepping implemented in \texttt{OrdinaryDiffEq.jl} \cite{DifferentialEquations.jl-2017}.
	We run the simulation until $t_f = 4.6$ (a plot of the density is provided in the Appendix in \cref{fig:KHI_Heuristic_DGMulti_rho}) with the $v$--adaptive scheme using the entropy production indicator with tolerated entropy production $\sigma_i = \sigma = 5 \cdot 10^{-3}$.
	Both the pure weak form and pure flux--differencing form simulations crash significantly earlier at around $t \approx 3.0$ and $t \approx 3.4$, respectively.
	The difference in mathematical entropy is displayed in \cref{fig:KHI_Heuristic_DGMulti_Entropy} for the three simulations.
	As expected, the adaptive scheme is dramatically cheaper than the pure flux--differencing form simulation as shown in \cref{fig:RunTimes_KHI_Heuristic_DGMulti}.
	In particular, the adaptive scheme is almost seven times faster than the pure flux--differencing form simulation.
	In comparison to the pure weak form simulation, the adaptive scheme is about $67\%$ more expensive.
	Surprisingly, the interface flux computation is for the \ac{WF} only scheme slightly more expensive than for the adaptive and pure \ac{FD} scheme, which results in the differences in "other" runtime.
	\subsubsection{Taylor--Green Vortex with Tolerated Entropy Increase}
	The three-dimensional Taylor--Green vortex is a well-known reference problem showcasing emerging turbulence from a simple initial condition which is followed by the decay of the turbulent structures accompanied by kinetic energy dissipation \cite{debonis2013solutions, bull2014simulation}.
	The initial condition in primitive variables is given by
	\begin{equation}
		\label{eq:IC_TaylorGreen}
		\boldsymbol u_\text{prim}(t_0=0, x, y, z) = \begin{pmatrix} \rho \\ v_x \\ v_y \\ v_z \\ p \end{pmatrix} = \begin{pmatrix} 1 \\ \sin(x) \cos(y) \cos(z) \\ -\cos(x) \sin(y) \cos(z) \\ 0 \\ p_0 + \frac{1}{16} \rho \Big(\big[\cos(2x) + \cos(2y)\big] \big[ 2+ \cos(2z) \big] \Big) \end{pmatrix} ,
	\end{equation}
	where $p_0 = \frac{\rho}{M^2 \gamma}$ with Mach number $M = 0.1$.
	The Prandtl and Reynolds number are $\text{Pr} = 0.72$ and $\text{Re} = 1600$, and the compressible Navier-Stokes equations are simulated on $\Omega = [0, 2 \pi]^3$ equipped with periodic boundaries until $t_f = 20.0$.
	We discretize the domain with a uniform grid of $32^3$ hexahedral elements and solution polynomials of degree $p = 3$, which is severe under--resolution for this testcase \cite{wang2013high, debonis2013solutions, bull2014simulation}.
	The interface fluxes are computed with the \ac{HLLC} flux \cite{toro1994restoration} and for the volume flux we employ the flux by Ranocha \cite{Ranocha2020Entropy}.
	As tolerated entropy increase (cf. \cref{eq:EntropyTargetDecay}) we prescribe $\sigma_i = \sigma = 4 \cdot 10^{-4}$.
	Time integration is performed with the nine--stage, fourth--order low--storage Runge--Kutta method from \cite{ranocha2022optimized} with a fixed timestep of $5.5 \cdot 10^{-3}$ which is optimized for compressible flows, implemented in \texttt{OrdinaryDiffEq.jl} \cite{DifferentialEquations.jl-2017}.

\begin{figure}
		\centering
		\subfloat[{Integrated kinetic energy \cref{eq:IntegratedKineticEnergy}.}]{
			\label{fig:TGV_EKin}
			\resizebox{.46\textwidth}{!}{\includegraphics{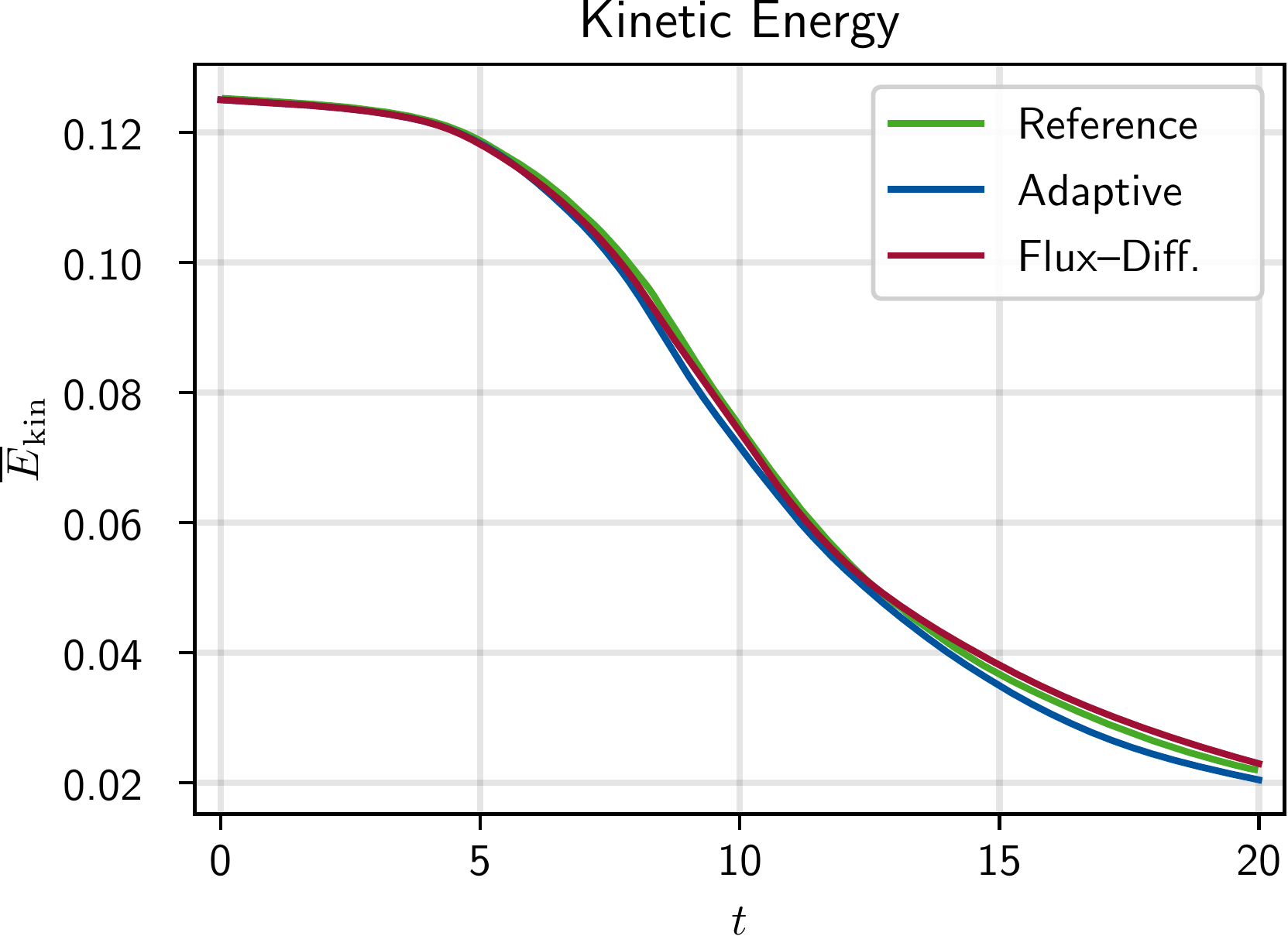}}
		}
		\hfill
		\subfloat[{Integrated enstrophy \cref{eq:IntegratedEnstrophy}.}]{
			\label{fig:TGV_Enstrophy}
			\resizebox{.46\textwidth}{!}{\includegraphics{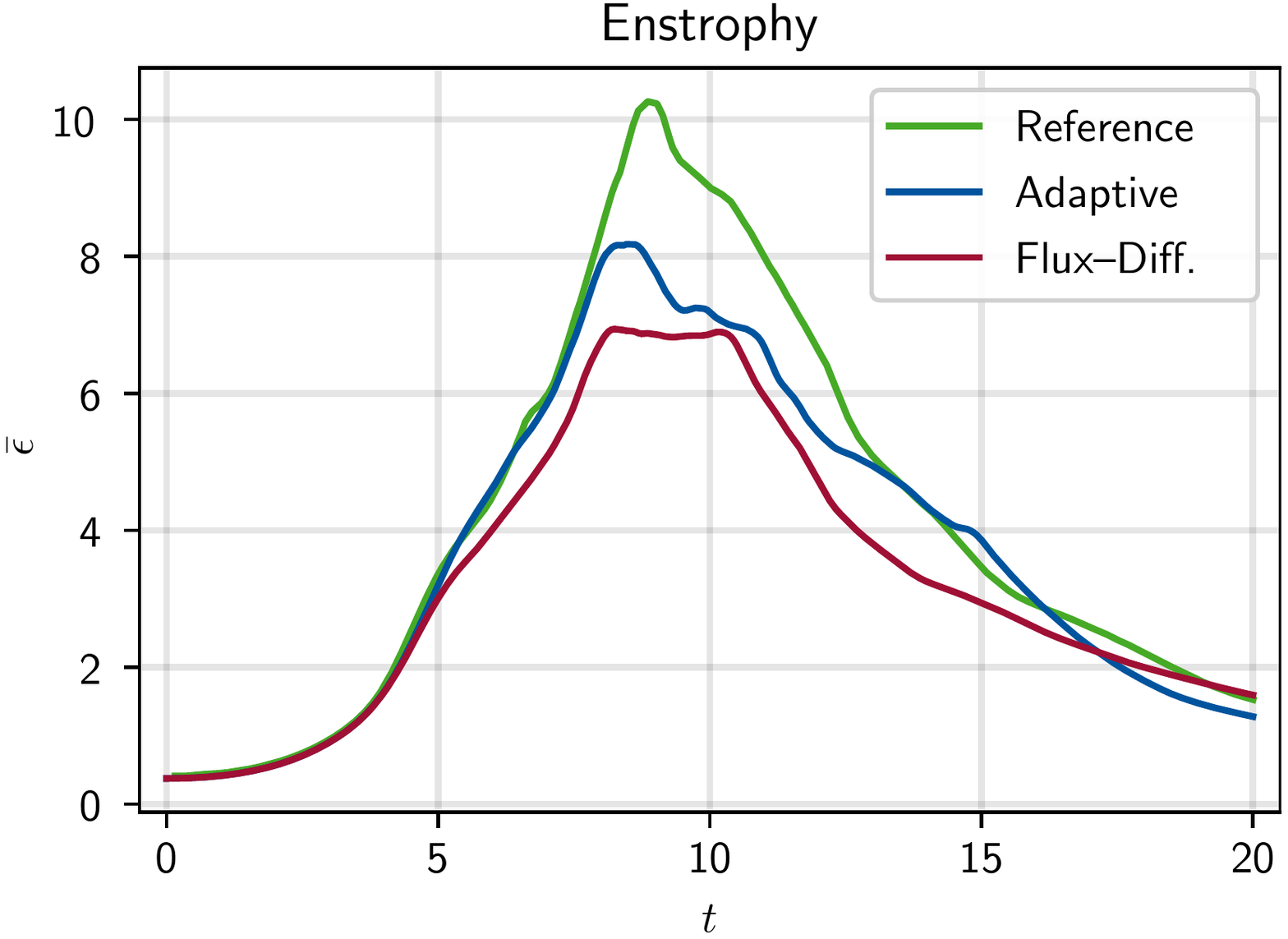}}
		}
	\caption[]
	{Kinetic energy (\cref{fig:TGV_EKin}) and enstrophy (\cref{fig:TGV_Enstrophy}) for the Taylor--Green vortex with tolerated element--wise entropy production $\sigma = 4 \cdot 10^{-4}$.
	Reference data taken from \cite{wang2013high}.}
	\label{fig:TGV_EKin_Enstropy}
	\end{figure}

\begin{figure}
		\centering
		\resizebox{.55\textwidth}{!}{\includegraphics{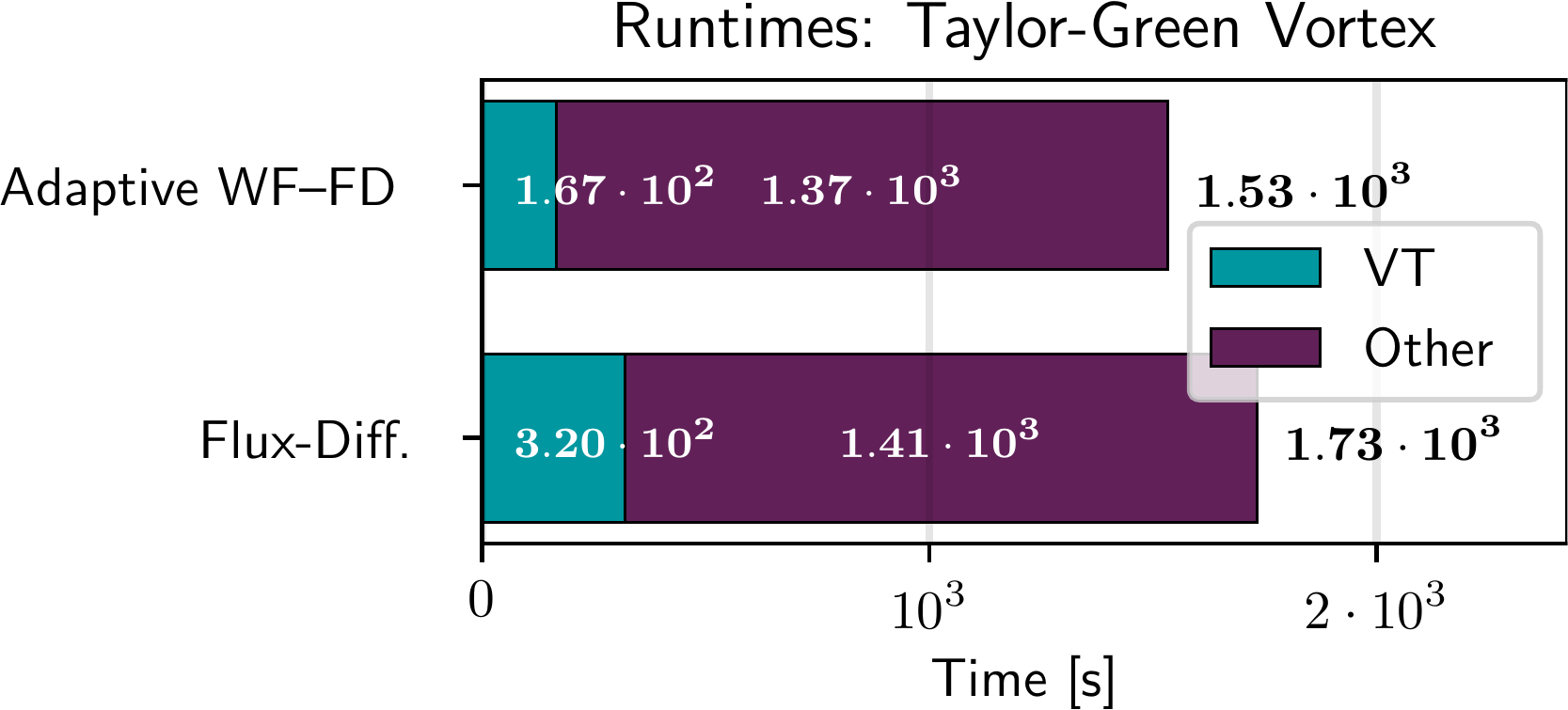}}
		\caption[Recorded runtime spent on volume integral/term computation for the Taylor--Green vortex for different volume term discretizations.]
		{Total \ac{RHS} and \ac{VT} computation runtime for the Taylor--Green vortex for $t_f = 20.0$.
		Volume term adaptivity is governed by the heuristic indicator with tolerated entropy increase $\sigma_i = \sigma = 4 \cdot 10^{-4}$.}
		\label{fig:RunTimes_TGV}
	\end{figure}

	The integrated enstrophy
	\begin{equation}
		\label{eq:IntegratedEnstrophy}
		\bar{\epsilon} \coloneqq \frac{1}{\rho_0 \vert \Omega \vert }  \int_\Omega \epsilon \, \text{d} \Omega,
		\quad \epsilon \coloneqq 0.5 \, \rho \, \boldsymbol{\omega} \cdot \boldsymbol{\omega}, \quad  \boldsymbol{\omega} \coloneqq \nabla \times \boldsymbol{v} \, ,
	\end{equation}
	and integrated kinetic energy 
	\begin{equation}
		\label{eq:IntegratedKineticEnergy}
		\overline{E}_\text{kin} = \frac{1}{\rho_0 \vert \Omega \vert } \int_\Omega \frac{1}{2} \rho \boldsymbol{v} \cdot \boldsymbol{v} \, \text{d} \Omega 
	\end{equation}
	showcase for the Taylor--Green vortex a characteristic behavior which allows for a meaningful comparison of different simulations.
	These are plotted in \cref{fig:TGV_EKin_Enstropy}.
	We observe slightly increased dissipation of kinetic energy for the adaptive scheme when compared to the pure \ac{FD} scheme, compare \cref{fig:TGV_EKin}.
	For enstrophy, however, we see that the adaptive scheme is able to follow the reference solution much better than the pure \ac{FD} scheme, see \cref{fig:TGV_Enstrophy}.
	In particular, up to about $t \approx 7.6$ the adaptive scheme is very close to the reference solution and produces also a distinct peak, albeit at a slightly earlier time.
	Tolerating some entropy increase results in a super--resolution alike behaviour for the enstrophy.
	We remark that this choice of $\sigma$ is close to the stability limit for this testcase, as the simulation crashes for $\sigma = 5 \cdot 10^{-4}$.
	In terms of performance, the $v$--adaptive scheme cuts the hyperbolic \ac{VT} roughly in half.
	For the Navier--Stokes equations, which introduce $12$ additional variables due to the $9$ velocity and $3$ temperature gradients, the overall runtime savings are only about $13\%$, cf. \cref{fig:RunTimes_TGV}.
	\subsubsection{Inviscid ONERA M6 Wing with Shock--Capturing Indicator}
	\label{subsubsec:ONERA_M6}
\begin{figure}
		\centering
		\subfloat[{Rescaled ONERA M6 wing geometry following \cite{onera_m6_original_report, slater_study_1}.
		}]{
			\label{fig:ONERA_M6_geometry}
			\centering
			\resizebox{.47\textwidth}{!}{\includegraphics{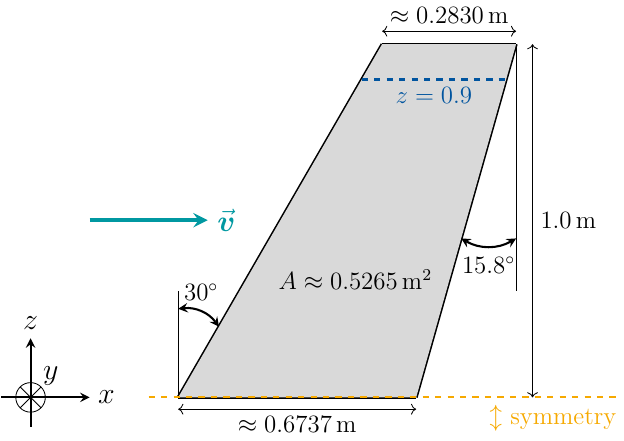}}
		}
		\hfill
		\subfloat[{Computational mesh (2D/surface element edges).
		}]{
			\label{fig:ONERA_M6_Mesh}
			\centering
			\includegraphics[width=0.45\textwidth]{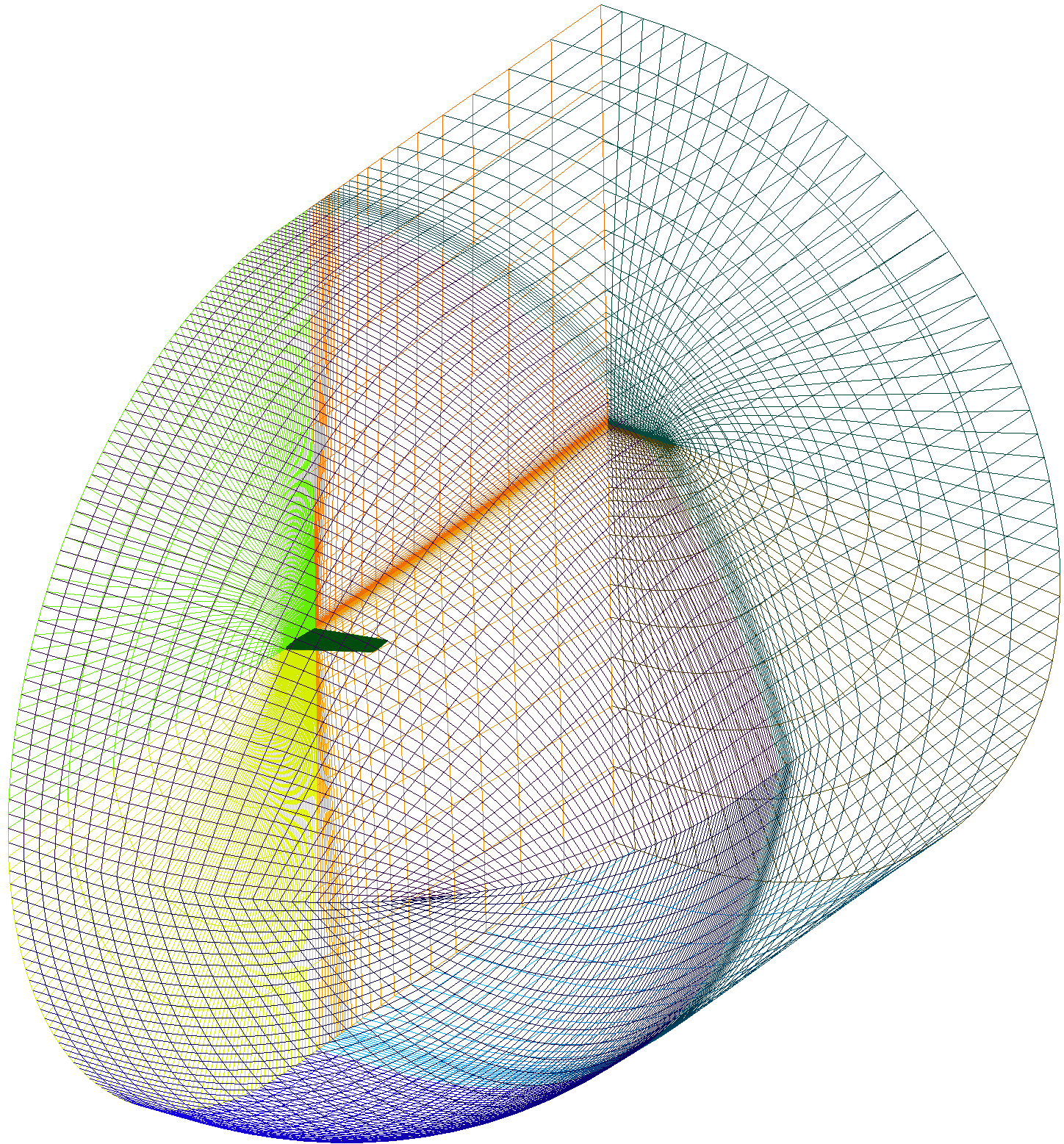}
		}
		\caption[ONERA M6 wing geometry and used computational mesh.]
		{ONERA M6 wing geometry \cref{fig:ONERA_M6_geometry} and used computational mesh \cref{fig:ONERA_M6_Mesh}.
		The displayed mesh is the sanitized version of the mesh from \cite{slater_study_1}.
		Note that the $y$-axis points into the plane in \cref{fig:ONERA_M6_geometry} and correspondingly upwards in \cref{fig:ONERA_M6_Mesh}.
		}
		\label{fig:ONERA_M6_geometry_mesh}
	\end{figure}

	We consider the inviscid flow around the ONERA M6 wing \cite{onera_m6_original_report, slater_study_1} at Mach number $\text{Ma} = 0.84$ and angle of attack $\alpha = 3.06^\circ$.
	In our study we follow the geometry as presented in \cite{onera_m6_original_report}, although we employ a rescaled variant with wingspan set to $1$.
	A sketch of the wing geometry in the $x-z$ plane is given by \cref{fig:ONERA_M6_geometry}.
	Details on the chamber of the wing are also provided in the cited reference.

	For the geometry proposed in \cite{onera_m6_original_report} a purely hexahedral mesh is available for download from publicly available NASA webpages as part of the NPARC Alliance Verification Archive \cite{slater_study_1}.
	The mesh as published is split into four different blocks (cf. edge coloring in \cref{fig:ONERA_M6_Mesh}) where a united version in \texttt{gmsh} format is available from the HiSA repository \cite{heyns2014modelling}.
	Some mesh sanitizing was necessary (see \cite{doehring2025paired} for details) to make this grid compatible with the \ac{DG} methodology as implemented in \texttt{Trixi.jl}.
	The final $294838$--element sanitized mesh is publicly available as part of the reproducibility repository which we provide publicly available on \texttt{GitHub} and archived on \texttt{Zenodo} \cite{doehring2026VTA_ReproRepo}.
	The mesh extends in $x$--direction from $-6.373$ to $ 7.411$, in $y$--direction from $-6.375$ to $6.375$ and in $z$--direction from $0$ to $7.411$.
	The 2D element edges are depicted for illustration in \cref{fig:ONERA_M6_Mesh}.

	At $z=0$ plane symmetry is enforced by rotating the symmetry-plane crossing velocity $v_z$ to zero and copying density and pressure from the domain.
	At all other outer boundaries we weakly enforce \cite{mengaldo2014guide} the freestream flow.
	The wing surface is modeled as a slip-wall \cite{Vegt2002SlipFB} with an analytic solution of the pressure Riemann problem \cite{toro2009riemann}.

	We represent the solution using $p=3$ degree element polynomials which results in about $94$ million total \ac{DoF}, i.e., almost $19$ million unknowns per solution field.
	The surface fluxes are computed with the Rusanov/\ac{LLF} flux \cite{RUSANOV1962304}.
	To run the simulation it suffices in this case to use the flux--differencing approach, i.e., no additional low--order stabilization is required.
	To select the cells which require stabilization, however, we employ the shock--capturing indicator from \cite{hennemann2021provably, rueda2021entropy} with indicator variable $\rho \cdot p$ and parameters $\beta_{\text{min}} = 0.01$, $\beta_\text{max} = 1.0$, thus only for $\beta_i > 0.01$ the more expensive flux--differencing volume term discretization is applied.
	We remark that this strategy is identical to the \textit{DG--ES} scheme from \cite{bilocq2025comparison} where the shock--capturing indicator is also used to distinguish between standard \ac{DG} and entropy--stable \ac{DG} cell-wise discretizations.
	The arising volume fluxes are discretized using the entropy--stable and kinetic--energy preserving flux by Ranocha \cite{Ranocha2020Entropy}.
	For time integration we employ fourth--order multirate Paired Explicit Runge--Kutta schemes \cite{vermeire2019paired, doehring2025fourth} with 15 different methods corresponding to $5$ to $17$ active stage evaluations.
	For the restarted simulation, only about 1140 out of the total 294838 elements require the more expensive flux--differencing volume term discretization, i.e., only roughly $0.39\%$ of the elements require the more expensive volume term discretization.

\begin{figure}
		\centering
		\resizebox{.55\textwidth}{!}{\includegraphics{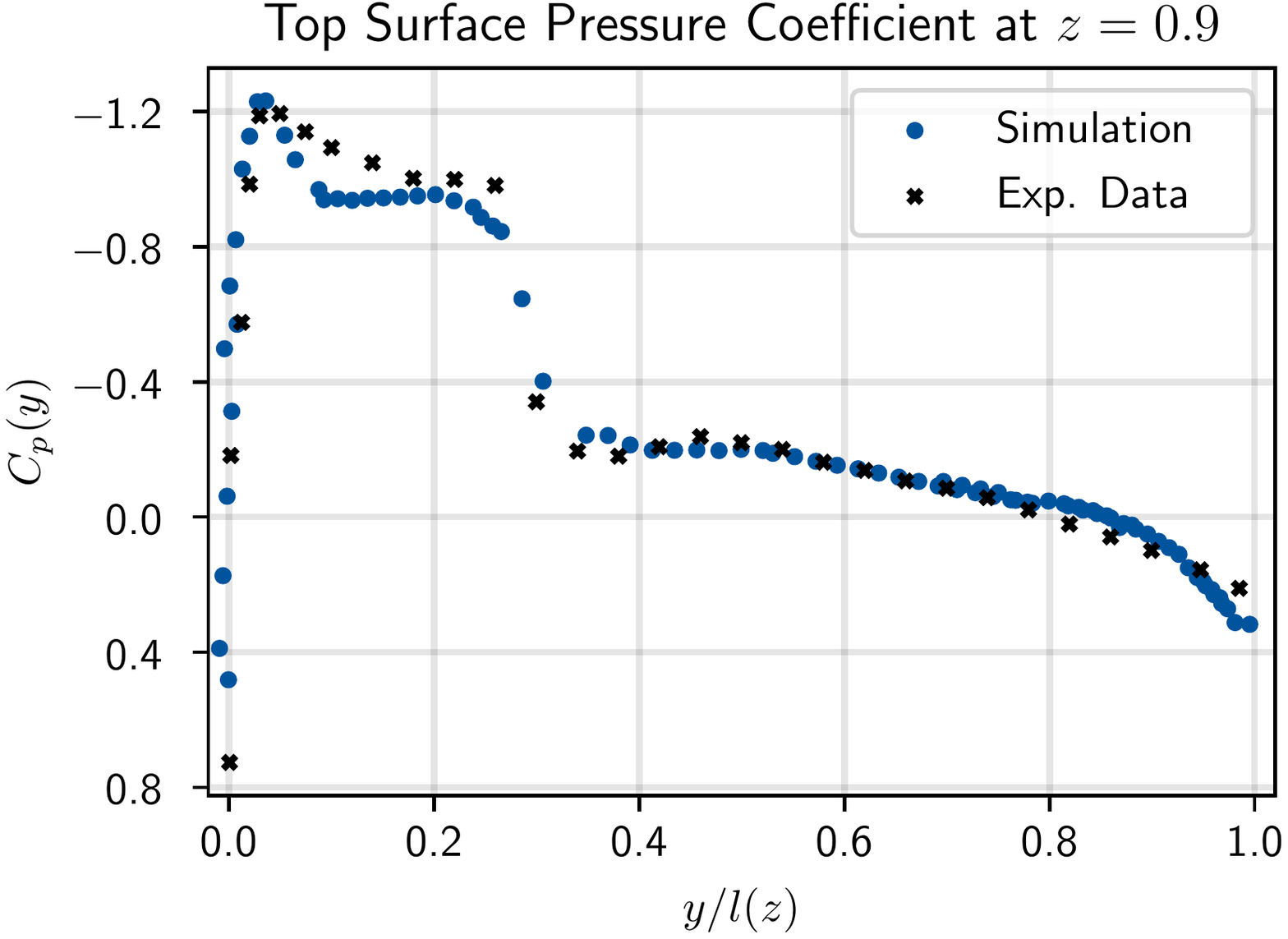}}
		\caption[Pressure coefficient recorded at $z = 0.9$ on the top surface of the ONERA M6 wing.]
		{Pressure coefficient \eqref{eq:PressureCoefficient} recorded at $z = 0.9$ (cf. \cref{fig:ONERA_M6_geometry}) on the top surface of the ONERA M6 wing.
		The numerical solution is compared to experimental data \cite{onera_m6_original_report}.}
		\label{fig:PressureCoefficientONERAM6_z090}
	\end{figure}

\begin{figure}
		\centering
		\includegraphics[width=0.65\textwidth]{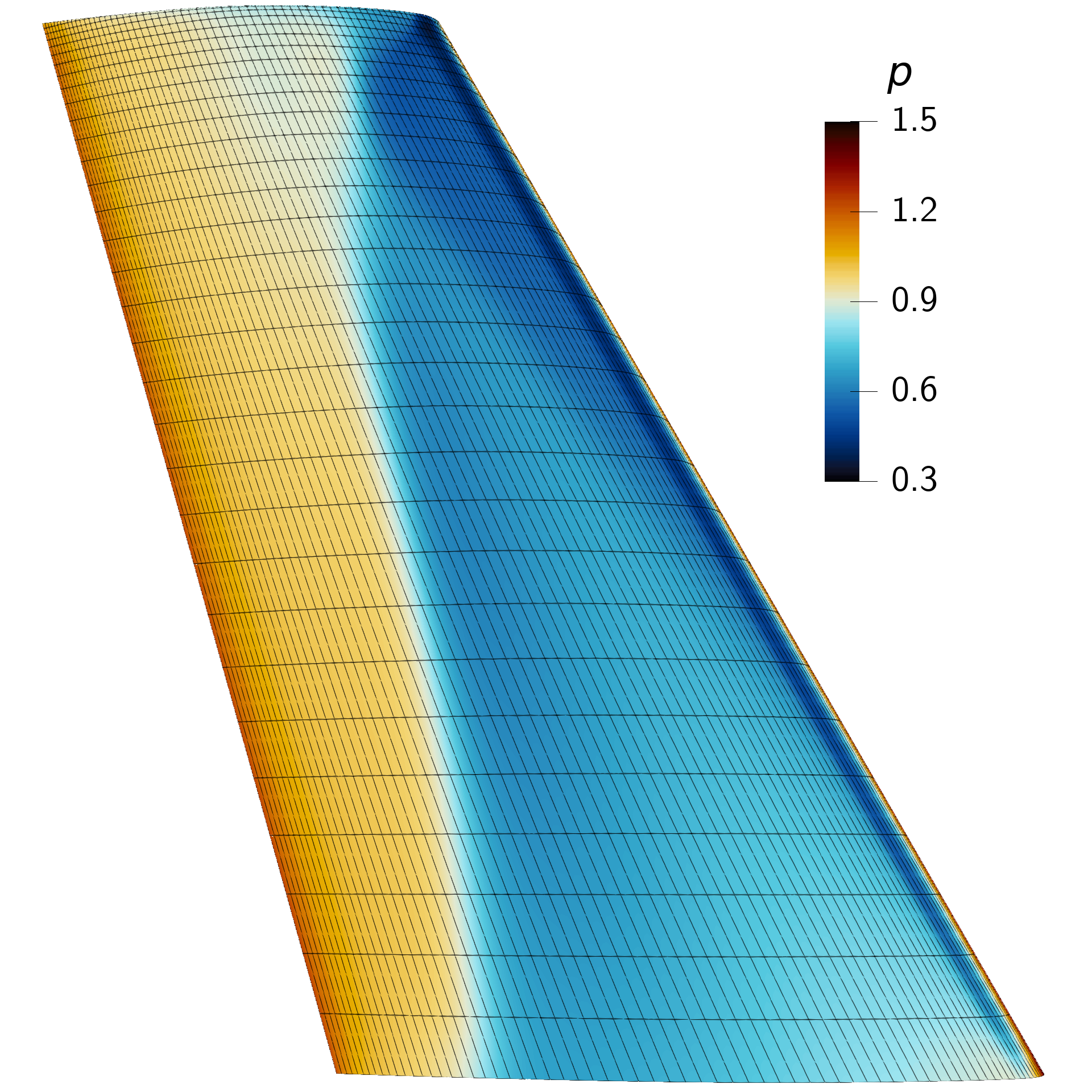}
		\caption[Pressure $p$ on the ONERA M6 wing's upper surface for the inviscid steady--state case.]
		{Pressure $p$ on the ONERA M6 wing's upper surface at $t_f = 6.05$ for the inviscid steady--state case.
		Note that in contrast to \cref{fig:ONERA_M6_geometry} the $y$--axis points in this figure out of the view plane.
		}
		\label{fig:ONERA_M6_Pressure_UpperWing}
	\end{figure}

\begin{figure}
		\centering
		\subfloat[{Runtimes for the simulation with multirate time integration.
		}]{
			\label{fig:RunTimes_OneraM6_Multirate}
			\centering
			\resizebox{.47\textwidth}{!}{\includegraphics{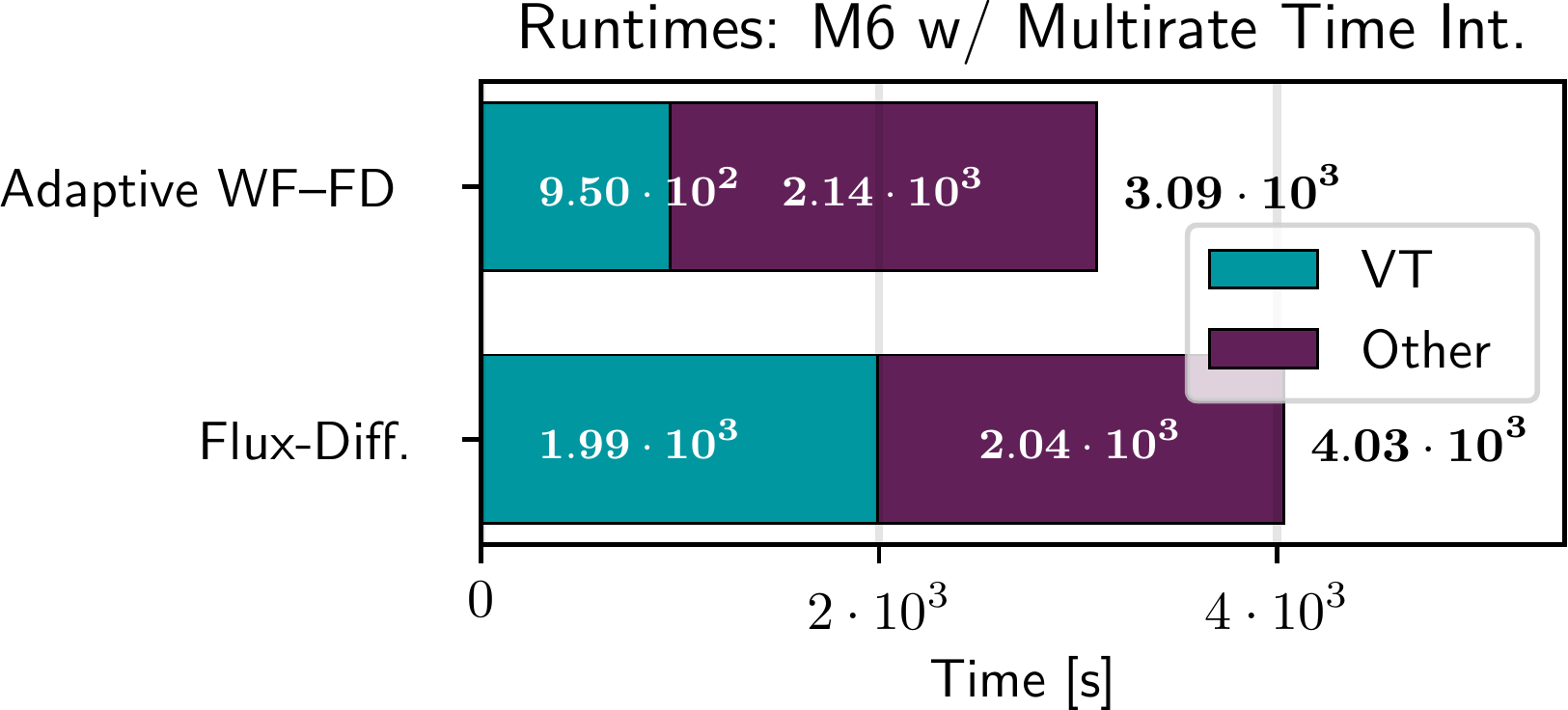}}
		}
		\hfill
		\subfloat[{Runtimes for the simulation with standard time integration.
		}]{
			\label{fig:RunTimes_OneraM6_Standard}
			\centering
			\resizebox{.47\textwidth}{!}{\includegraphics{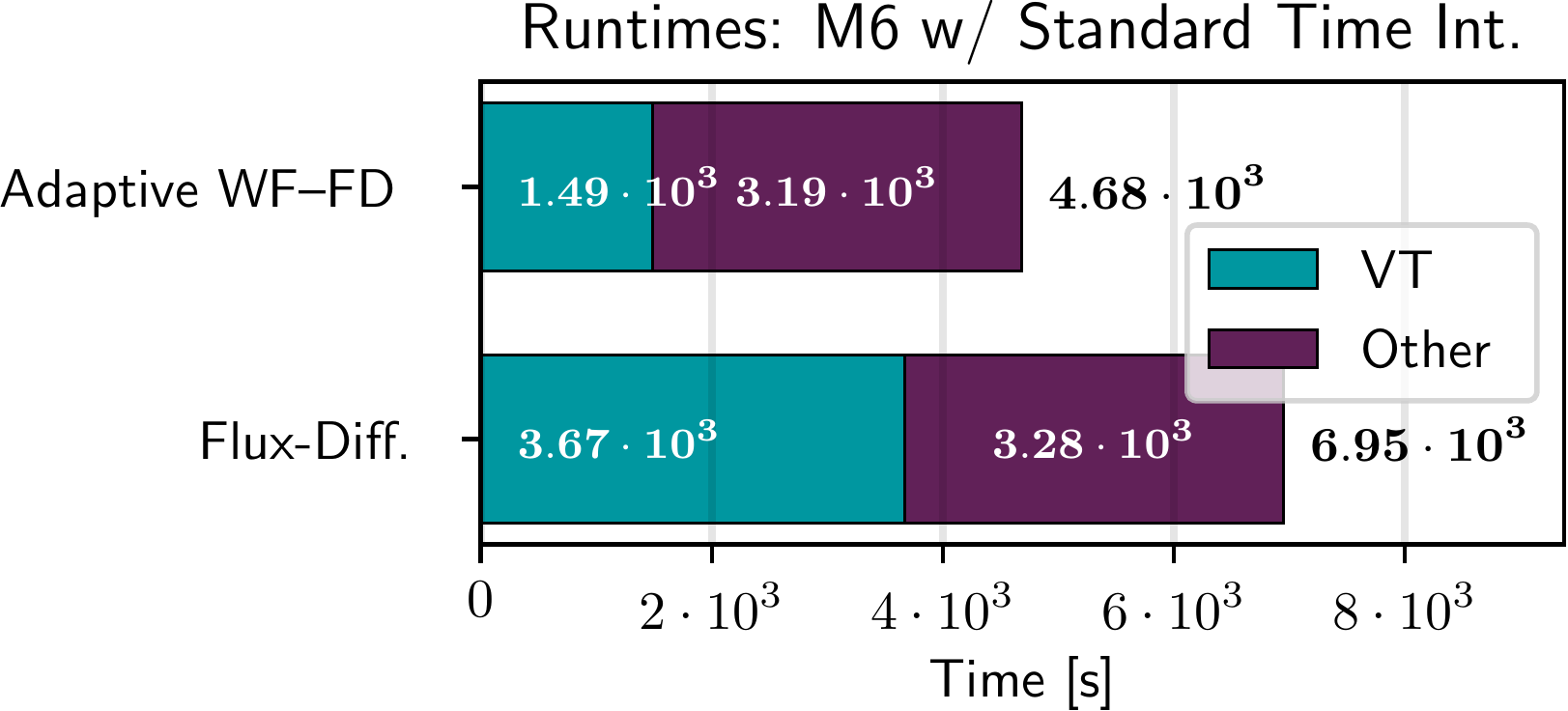}}
		}
		\caption[Recorded runtime spent on volume integral/term computation for the ONERA M6 wing for different volume term discretizations and time integration strategies.]
		{Total \ac{RHS} and \ac{VT} computation runtime for the ONERA M6 wing with multirate time integration (\cref{fig:RunTimes_OneraM6_Multirate}) and standard time integration (\cref{fig:RunTimes_OneraM6_Standard}).
		Volume term adaptivity is governed by the shock--capturing indicator from \cite{persson2006sub, hennemann2021provably}.}
		\label{fig:RunTimes_OneraM6}
	\end{figure}

	We compute the lift coefficient 
	\begin{equation}
		\label{eq:LiftCoeff}
		C_L = C_{L,p} = \oint_{\partial \Omega} \frac{ p \, \boldsymbol n \cdot \boldsymbol t^\perp}{0.5 \rho_{\infty} U_{\infty}^2 L_{\infty}} \, \text{d} S
	\end{equation}
	for $\alpha = 3.06\degree$ based on the reference area $A_\infty \approx 0.5265$.
	We obtain $C_L \approx 0.2951$ which is about \qty{3}{\percent} larger of what is reported in one of the tutorials of the \texttt{SU2} code \cite{economon2016su2}.
	This might be due to the fact that in the latter case the actual area of the wing is used which is necessarily larger than the projected referenced area $A_\infty$ (which we employ) due to the curvature of the wing.
	In that tutorial a lift coefficient of $C_L \approx 0.2865$ is reported, obtained from a second--order in space, first--order in (pseudo) time simulation on a $\sim 5.8 \cdot 10^5$ element tetrahedral mesh.
	Furthermore, we also provide a plot of the pressure coefficient
	\begin{align}
		\label{eq:PressureCoefficient}
		C_p(x) \coloneqq \frac{p(x) - p_{\infty}}
		{0.5 \, \rho_{\infty} U_{\infty}^2}
	\end{align}
	at $z=0.9$ on the upper surface of the wing (cf. \cref{fig:ONERA_M6_geometry}) in \cref{fig:PressureCoefficientONERAM6_z090}.
	The numerical results agree well with the experimental data \cite{onera_m6_original_report} apart from deviations at the leading, especially given that this is an inviscid simulation without viscosity and thus any turbulence modeling.
	For illustration a plot of the pressure on the upper wing surface at $t_f = 6.05$ is provided in \cref{fig:ONERA_M6_Pressure_UpperWing}.
	To compare performance, we record the runtime spent on volume integral/term computation $\tau_\text{VI}$ for both multirate time integration using the fourth--order P-ERK schemes \cite{vermeire2019paired, doehring2025fourth} and standalone, optimized explicit Runge--Kutta time integration over the $t \in[6.049, 6.05]$ interval.
	The results are provided in \cref{fig:RunTimes_OneraM6}.
	As for the previous examples, we observe significant computational savings when employing the \ac{WF}---\ac{FD} combination in comparison to the pure \ac{FD} volume term discretization.
	For the simulation with multirate time integration speedup of about $2.09$ is observed, while for standard explicit Runge--Kutta time integration speedup of more than $2.46$ is observed, cf. \cref{fig:RunTimes_OneraM6}.
	This manifests in overall about $32 \%$ and $35 \%$ faster \ac{RHS} computations, respectively.
	This discrepancy is due to the fact that the multirate time integration already reduces the number of right--hand side evaluations significantly in regions with large cells, i.e., far away from the wing.
	The small cells close to the wing, which demand more stage evaluations, are also more prone to use the more expensive \ac{FD} volume term discretization.
	Thus, the relative advantage of the \ac{WF} volume term discretization is slightly less pronounced for multirate time integration.
	\subsection{Weak Form --- Blended DG--FV Combination}
	We now consider examples where additional low--order stabilization is required in critical regions of the flow.
	To this end, we employ the blended \ac{DG}--\ac{FV} scheme from \cite{hennemann2021provably, rueda2021entropy} with the troubled--cell indicator from \cite{persson2006sub, hennemann2021provably} which we introduced in \cref{sec:shock_capturing_indicator}.

	Now, however, we employ the \ac{WF} volume integral as the default, unlimited ($\beta = 0$) volume term instead of the \ac{FD} volume term discretization.
	This will be more efficient if only a small number of elements is troubled and requires limiting, and furthermore the more frequent use of the \ac{WF} volume term does not result in more troubled and thus limited cells.
	\subsubsection{Rayleigh--Taylor Instability}
	\label{subsubsec:RTI}
	The Rayleigh--Taylor instability occurs at sharp interfaces of fluids of different density in a gravitational field.
	The heavy fluid tends to fall down while the light fluid rises, leading to the formation of characteristic mushroom cap structures.
	We consider a setup of the Rayleigh--Taylor instability based on \cite{shi2003resolution, remacle2003adaptive}.
	The domain is set to $\Omega = \left[0, \sfrac{1}{4}\right] \times [0, 1]$ with the interface placed at $y = \sfrac{1}{2}$.
	In $x$--direction periodic boundary conditions are imposed, while in $y$--direction weakly imposed Dirichlet boundary conditions are used \cite{carlson2011inflow, mengaldo2014guide}.
	We set the ratio of specific heats $\gamma = \sfrac{5}{3}$ and the initial state to
	\begin{equation}
		\label{eq:RTI_IC}
		\boldsymbol u_\text{prim}(t_0=0, x, y) = 
		\begin{pmatrix} \rho \\ v_x \\ v_y \\ p \end{pmatrix} 
		= 
		\begin{pmatrix} 
			\begin{cases} 2 & y \geq 0.5 \\ 1 & y < 0.5 \end{cases} \\ 0.0 \\ c k \cos(8 \pi x) \sin^6(\pi y)\\ \begin{cases} -2 y + 3 & y \geq 0.5 \\ -y + 2.5 & y < 0.5 \end{cases}  
		\end{pmatrix} 
	\end{equation}
	where $ c = \sqrt{\frac{\gamma p}{\rho}}$ denotes the speed of sound and $k= -0.025$ quantifies the strength of the initial perturbation.
	Gravity is included via the source term
	\begin{equation}
		\boldsymbol s(\boldsymbol u) = \begin{pmatrix} 0 \\ 0 \\ g \rho \\ g \rho v_y \end{pmatrix}, \quad g = -1 \: .
	\end{equation}
	%
	
\begin{figure}
		\centering
		\subfloat[{Density at $t_f = 3.0$.}]{
			\label{fig:RTI_Solution}
			\includegraphics[height=0.5\textheight]{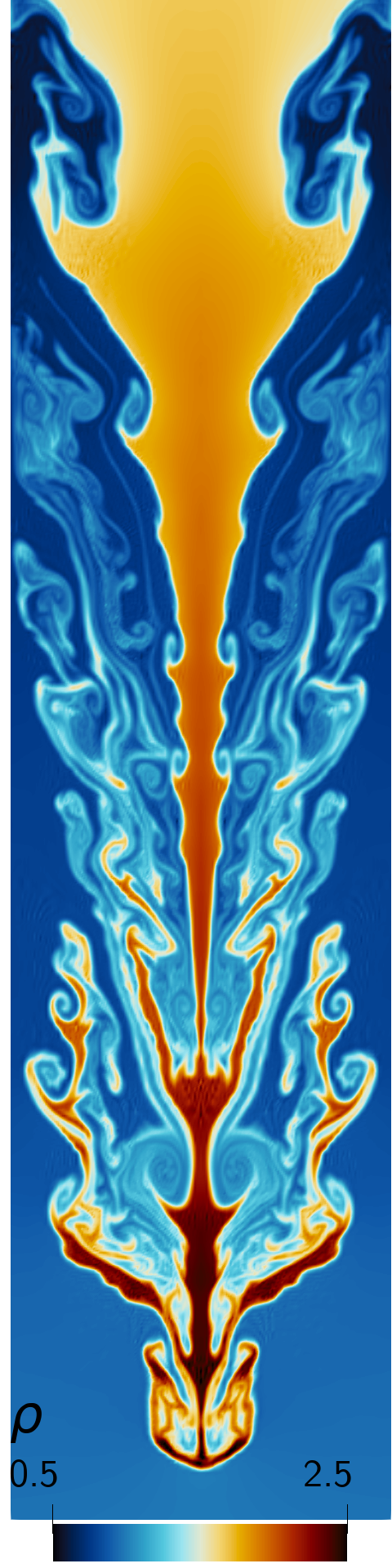}
		}
		\hspace{2cm}
		\subfloat[{Mesh at $t_f = 3.0$.}]{
			\label{fig:RTI_Mesh}
			\includegraphics[height=0.5\textheight]{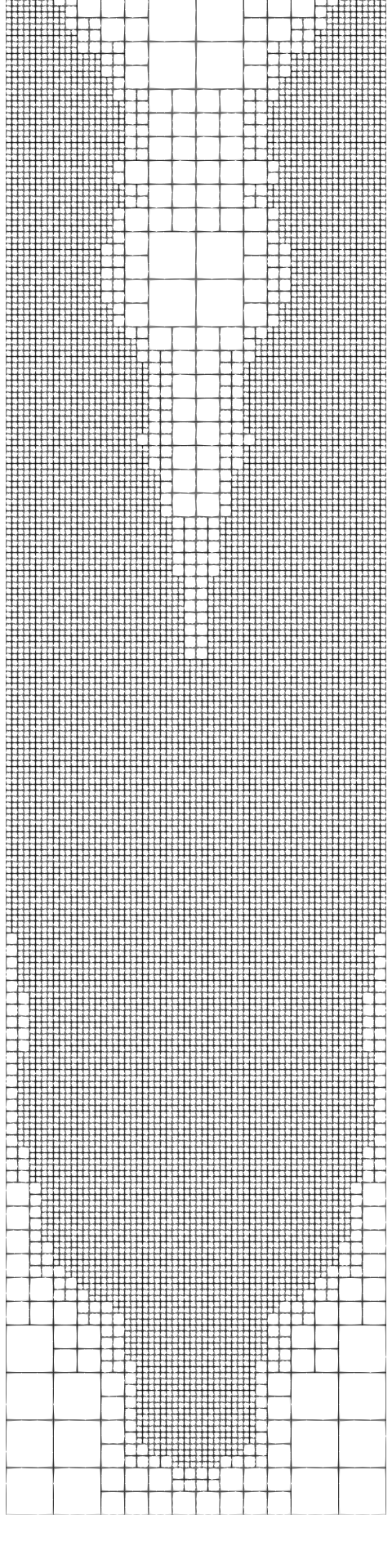}
		}
		\caption[]
		{Density and adaptive mesh at final time $t_f = 3.0$ for the Rayleigh--Taylor instability with \ac{AMR} and \ac{WF} --- \ac{FD}--\ac{FV} volume term combination.
		Note that symmetry across the horizontal centerline is preserved.}
		\label{fig:RTI_Solution_Mesh}
	\end{figure}

	We employ solution polynomials of degree $p = 3$ and \ac{HLL} surface flux \cite{harten1983upstream} with Davis--type wave speed estimates \cite{doi:10.1137/0909030}.
	For stabilizing, we use the blended \ac{DG} second--order \ac{FV} scheme from \cite{rueda2021entropy} with monotonized central limiter \cite{van1977towards} and the entropy--conservative and kinetic energy preserving volume flux from Ranocha \cite{Ranocha2020Entropy}.
	The shock--capturing indicator from \cite{persson2006sub, hennemann2021provably} is based on the density $\rho$ and is parametrized through $\beta_\text{min} = 0.001$, $\beta_\text{max} = 0.5$.
	We conduct two simulations up to final time $t_f = 3.0$.
	First, we consider a uniformly discretized domain with $64 \times 256$ quadrilateral elements.
	Second, we employ \ac{AMR} starting with a very coarse mesh with $2 \times 8$ elements which can be refined 5 times, i.e., allowing the same smallest cell size as for the uniform mesh.
	The solution alongside the mesh at final time is shown in \cref{fig:RTI_Solution_Mesh}.
	\ac{AMR} is based on the indicator by Löhner \cite{lohner1987adaptive} as implemented in the \texttt{FLASH} code \cite{fryxell2000flash} with density as refinement variable.
	Time integration for the \ac{AMR} case is performed with an adaptive nine--stage, fourth--order Runge--Kutta method from \cite{ranocha2022optimized} implemented in \texttt{OrdinaryDiffEq.jl} \cite{DifferentialEquations.jl-2017}.

\begin{figure}
		\centering
		\subfloat[{Runtimes for the simulation with \ac{AMR}.
		}]{
			\label{fig:RunTimes_RTI_AMR}
			\centering
			\resizebox{.47\textwidth}{!}{\includegraphics{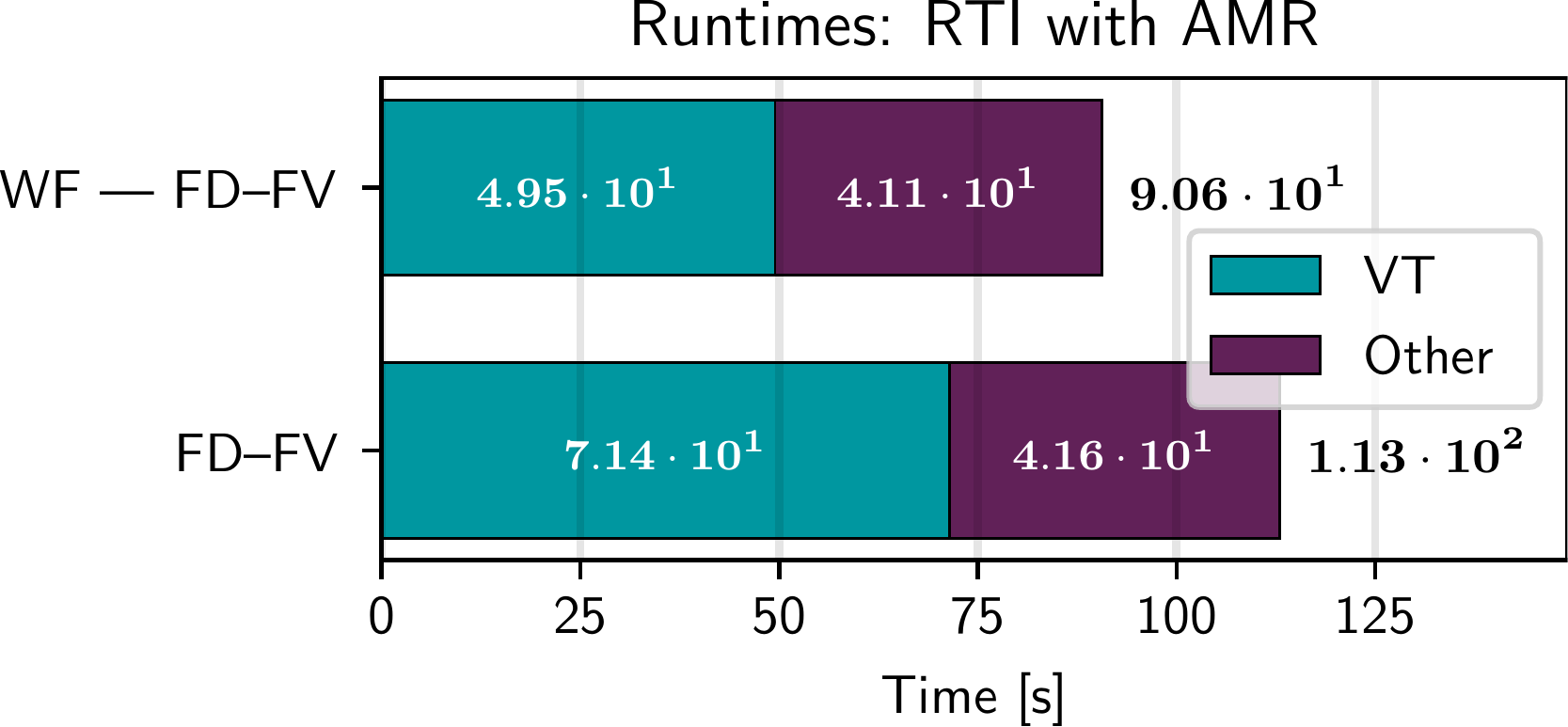}}
		}
		\hfill
		\subfloat[{Runtimes for the simulation on a constant uniform mesh.
		}]{
			\label{fig:RunTimes_RTI_Uniform}
			\centering
			\resizebox{.47\textwidth}{!}{\includegraphics{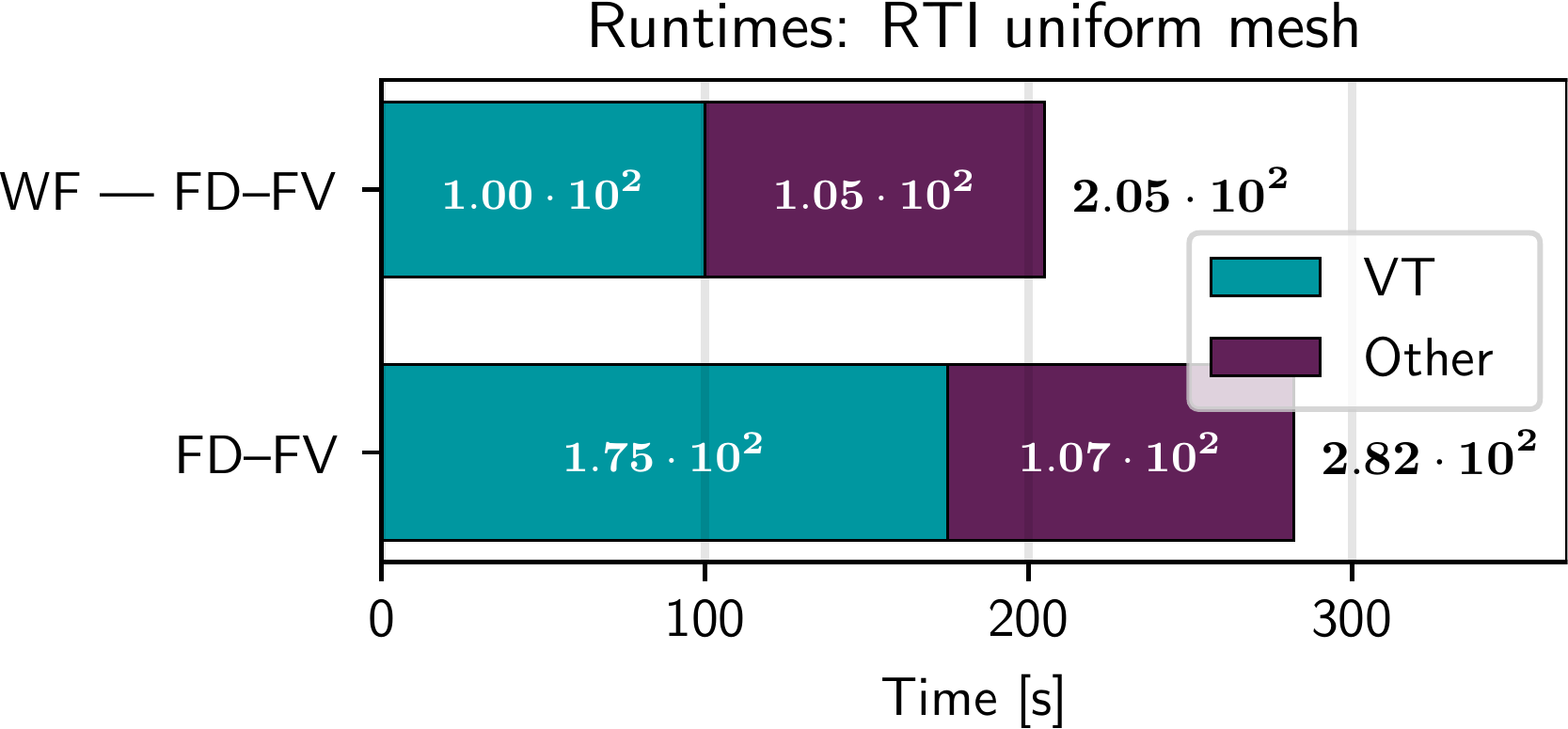}}
		}
		\caption[Recorded runtime spent on volume integral/term computation for the Rayleigh--Taylor instability for different volume term discretizations.]
		{Total \ac{RHS} and \ac{VT} computation runtime for the Rayleigh--Taylor instability with \ac{AMR} (\cref{fig:RunTimes_RTI_AMR}) and uniform constant mesh (\cref{fig:RunTimes_RTI_Uniform}).}
		\label{fig:RunTimes_RTI}
	\end{figure}

	For this setup, we observe that the volume term computation for the standard \ac{FD} discretization is about $35\%$ more expensive than the \ac{WF} discretization as shown in \cref{fig:RunTimes_RTI_AMR}.
	Due to the \ac{AMR} setup the solver automatically coarsens regions of the flow where no stabilization is required.
	Thus, for runs with \ac{AMR} we expect moderate speedup, as cells on which the \ac{WF} volume term discretization could be employed are reduced to a minimum.
	To confirm this, we also provide runtimes for the run on the uniform $64 \times 256$ mesh in \cref{fig:RunTimes_RTI_Uniform}.
	In this case, we observe that the \ac{FD} volume term discretization is $1.53$ times as expensive as the \ac{WF} discretization.
	For the uniform mesh case, \ac{CFL} based timestepping with the five--stage, fourth--order \ac{SSP} Runge--Kutta from \cite{ruuth2006global} is more efficient and thus employed.
	\subsubsection{3D Sedov Blast Wave}
	\label{subsubsec:SedovBlastWave}
\begin{figure}
		\centering
		\subfloat[{Runtimes for the simulation with \ac{AMR}.
		}]{
			\label{fig:RunTimes_SedovBlast_AMR}
			\centering
			\resizebox{.47\textwidth}{!}{\includegraphics{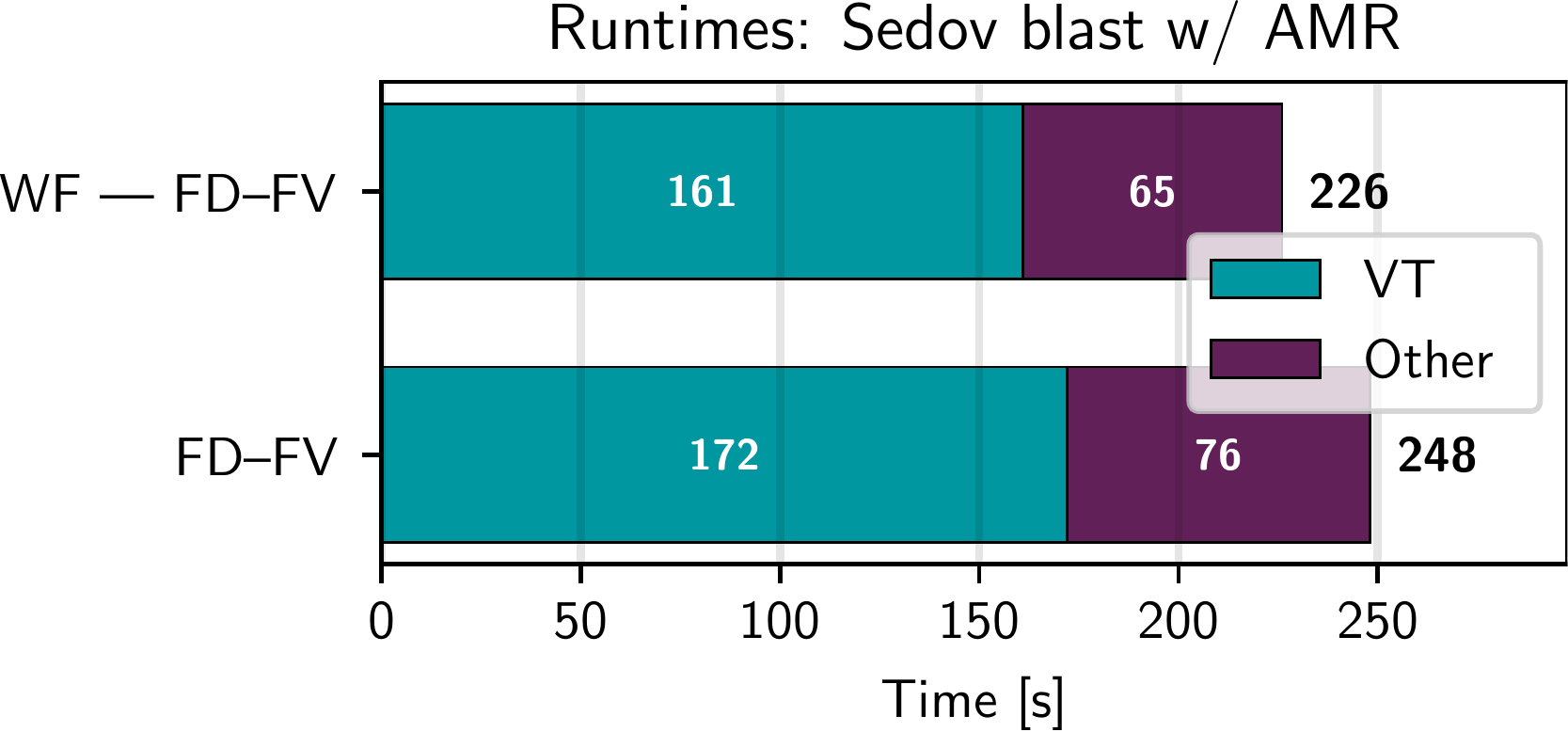}}
		}
		\hfill
		\subfloat[{Runtimes for the simulation on a constant uniform mesh.
		}]{
			\label{fig:RunTimes_SedovBlast_Uniform}
			\centering
			\resizebox{.47\textwidth}{!}{\includegraphics{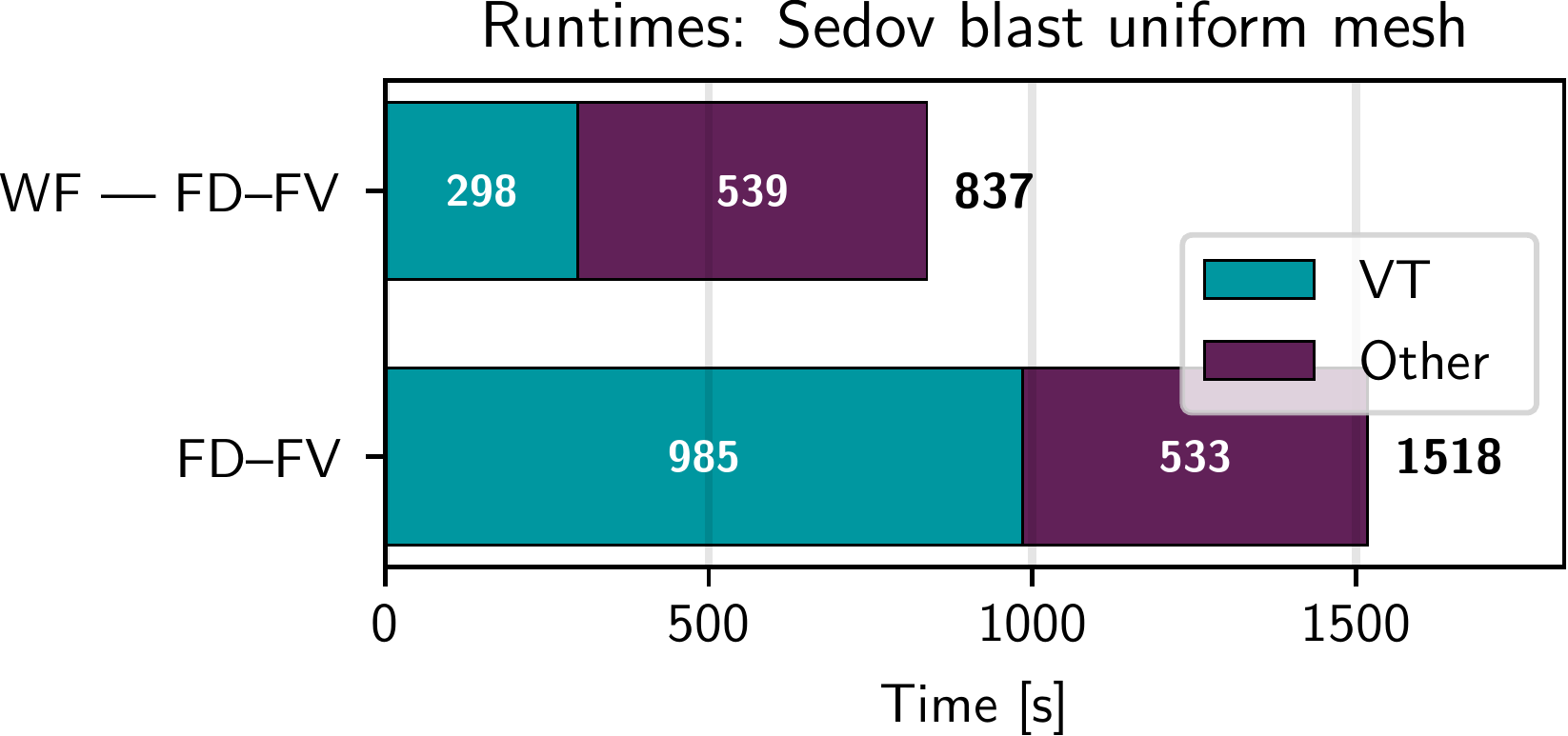}}
		}
		\caption[Recorded runtime spent on volume integral/term computation for the Sedov blast wave for different volume term discretizations.]
		{Total \ac{RHS} and \ac{VT} computation runtime for the Sedov blast wave with \ac{AMR} (\cref{fig:RunTimes_SedovBlast_AMR}) and uniform constant mesh (\cref{fig:RunTimes_SedovBlast_Uniform}).}
		\label{fig:RunTimes_SedovBlast}
	\end{figure}

\begin{figure}[ht]
		\includegraphics[width=0.46\textwidth]{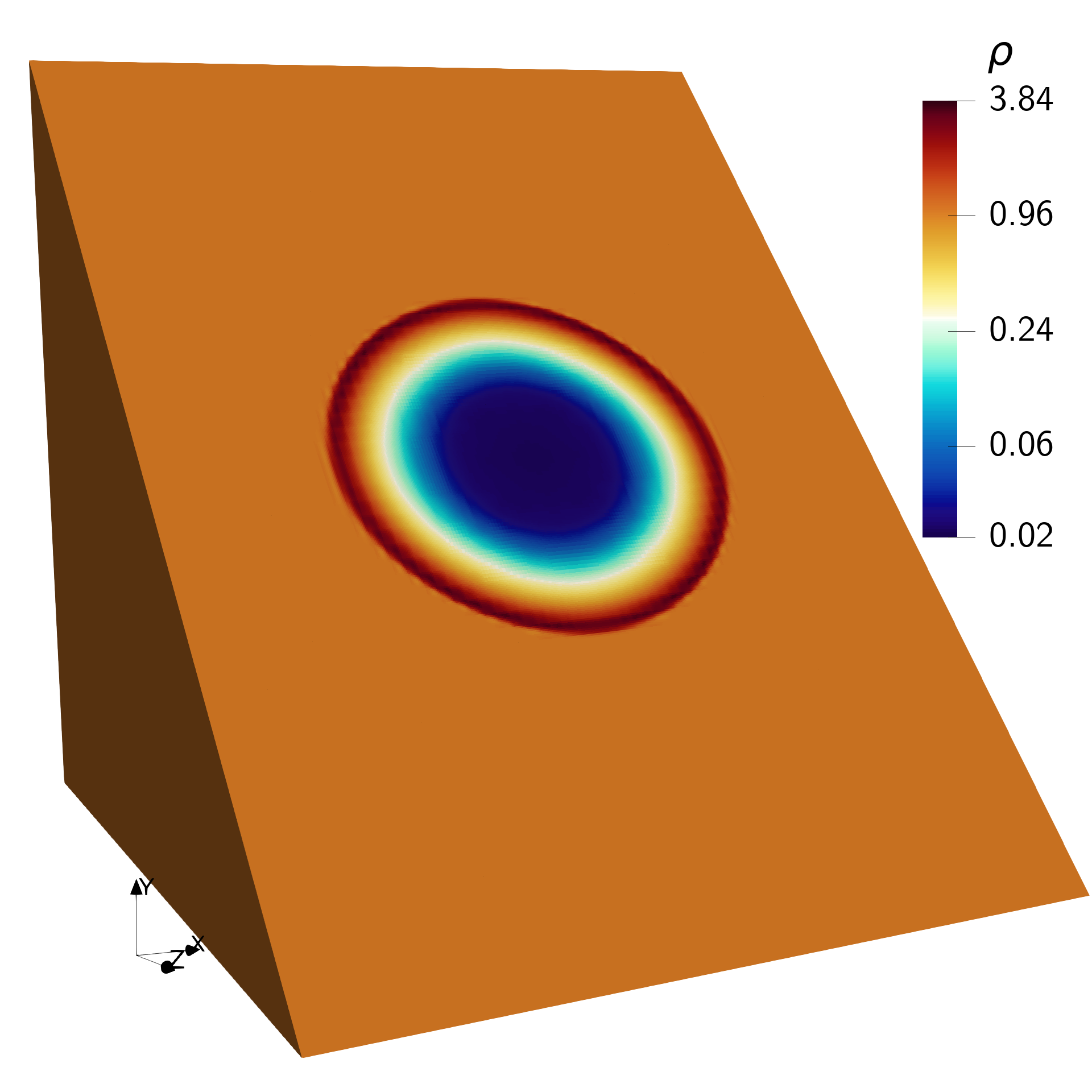}
		\caption[]{Density at final time $t_f = 1.0$ for the Sedov blast wave from \cref{subsubsec:SedovBlastWave} obtained with the \ac{WF}---\ac{FD}--\ac{FV} adaptive scheme.
		Note that the coloring is based on the logarithm of the density.}
		\label{fig:SedovBlastWave_rho}
	\end{figure}

	We consider the classic Sedov blast wave problem \cite{sedov1946propagation} with parametrization from \texttt{FLASH} code \cite{flashuserguide, fryxell2000flash}.
	The blast energy $E = 1$ is deposited into a small spherical region of radius $r = 3.5 \Delta x_\text{min}$ at the center of the domain $\Omega = [-2, 2]^3$.
	The blast pressure is set to 
	\begin{equation}
		p(r) = \frac{3(\gamma - 1) E}{4 \pi r^3}
	\end{equation}
	and the ambient pressure to $p_\text{ambient} = 10^{-5}$.
	Density and velocity are initialized to $\rho = 1$ and $\boldsymbol{v} = \boldsymbol{0}$, respectively.
	The surface fluxes are computed with the Rusanov/\ac{LLF} flux \cite{RUSANOV1962304} and the volume fluxes with the entropy--conservative and entropy conservative flux from Chandrashekar \cite{Chandrashekar_2013}.
	To stabilize the $p=3$ scheme, we use the blended \ac{DG}-- first order \ac{FV} scheme from \cite{hennemann2021provably}.
	The blending parameter $\beta$ is governed by the shock--capturing indicator from \cite{persson2006sub, hennemann2021provably} with indicator variable $\rho \cdot p$ and parameters $\beta_{\text{min}} = 0.001$, $\beta_\text{max} = 0.5$.
	The simulation is performed on an adaptive five--level mesh (minimum edge length $\Delta x = \sfrac{1}{16}$) where the indicator from the shock--capturing scheme is also employed for refinement and coarsening.
	Additionally, we also perform a simulation on a uniform mesh with $64^3$ elements which corresponds to the smallest cell size of the adaptive mesh.
	The simulation is advanced until final time $t_f = 1.0$ with the five--stage, fourth--order \ac{SSP} Runge--Kutta from \cite{ruuth2006global} with \ac{CFL} based timestepping.
	The recorded runtimes for the \ac{WF}---\ac{DG}--\ac{FV} and \ac{FD}---\ac{DG}--\ac{FV} volume term discretizations are provided in \cref{fig:RunTimes_SedovBlast}.
	Similar to the Rayleigh--Taylor instability from \cref{subsubsec:RTI}, we observe a moderate speedup of about $1.07$ for the \ac{AMR} case and a very pronounced speedup of almost $3.3$ for the uniform mesh case.
	A plot of the density at final time $t_f = 1.0$ is provided in \cref{fig:SedovBlastWave_rho}.
	\subsubsection{Transonic Viscous Flow over RAE2822 Airfoil}
\begin{figure}
		\centering
		\subfloat[{Horizontal velocity $v_x$. Note the shock--induced flow separation on the upper side of the airfoil.}]{
			\label{fig:RAE2822_vx}
			\includegraphics[width=0.95\textwidth]{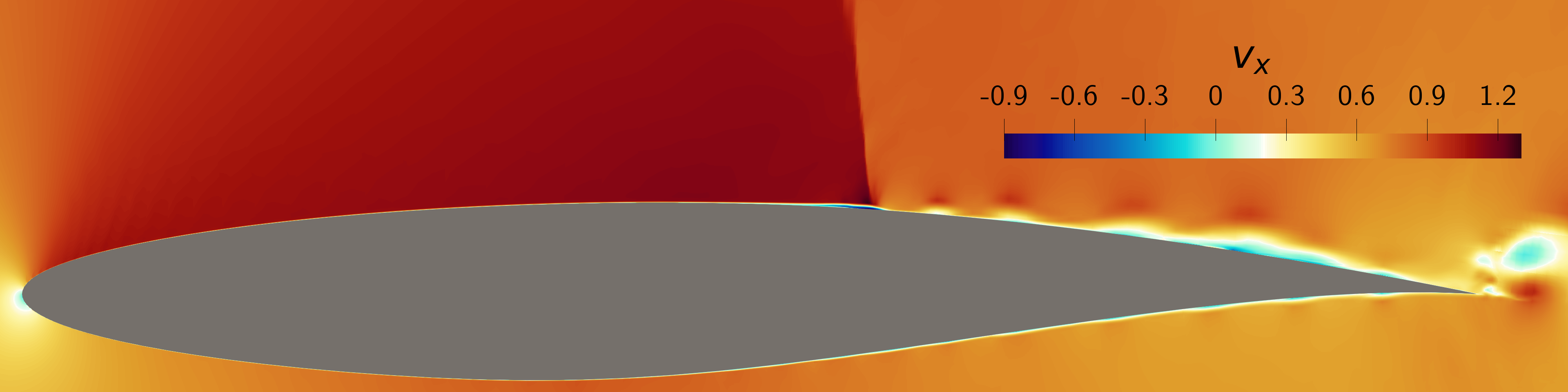}
		}
		\\
		\subfloat[{Blending parameter $\beta$ (recall \cref{eq:BlendedDG-FV}) governing the \ac{DG}--\ac{FV} combination.}]{
			\label{fig:RAE2822_alpha}
			\includegraphics[width=0.95\textwidth]{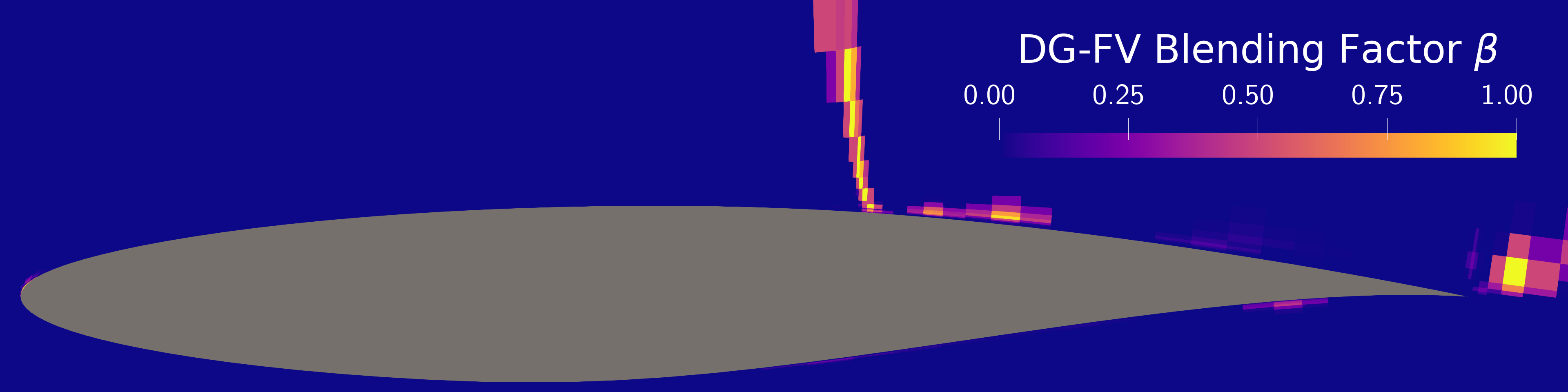}
		}
		\caption[Near airfoil plot of the horizontal velocity $v_x$ and blending value for the transonic viscous flow over the RAE2822 airfoil at $t_f = 26 t_c$.]
		{Near airfoil plot of the horizontal velocity $v_x$ (\cref{fig:RAE2822_vx}) and \ac{DG}--\ac{FV} blending value (\cref{fig:RAE2822_alpha}) at final time $t_f = 26 t_c$ for the transonic viscous flow over the RAE2822 airfoil.
		}
		\label{fig:RAE2822_vx_alpha}
	\end{figure}

	As the final example we consider the transonic, viscous flow over the RAE2822 airfoil \cite{cook1979aerofoil}.
	Our setup follows case 2.2 from the first International Workshop on High-Order CFD Methods \cite{case22_1stHiCFD}.
	Correspondingly, the Mach number is set to $\text{Ma}_\infty = 0.734$, the Reynolds number is $\text{Re}_\infty = 6.5 \cdot 10^6$, the Prandtl number is $\text{Pr}_\infty = 0.71$, and the angle of attack to $\alpha = 2.79^\circ$ \cite{case22_1stHiCFD, AGARD_702}.
	This testcase is particularly interesting due to shock--triggered flow separation.
	The suggested setup for this testcase is a steady \ac{RANS} simulation with turbulence modeling.
	Here, we consider an unsteady simulation without subgrid modeling.
	Nevertheless, this simulation setup recovers the expected separation as depicted in \cref{fig:RAE2822_vx}, paired with unsteady shock position oscillations.

\begin{figure}
		\centering
		\resizebox{.55\textwidth}{!}{\includegraphics{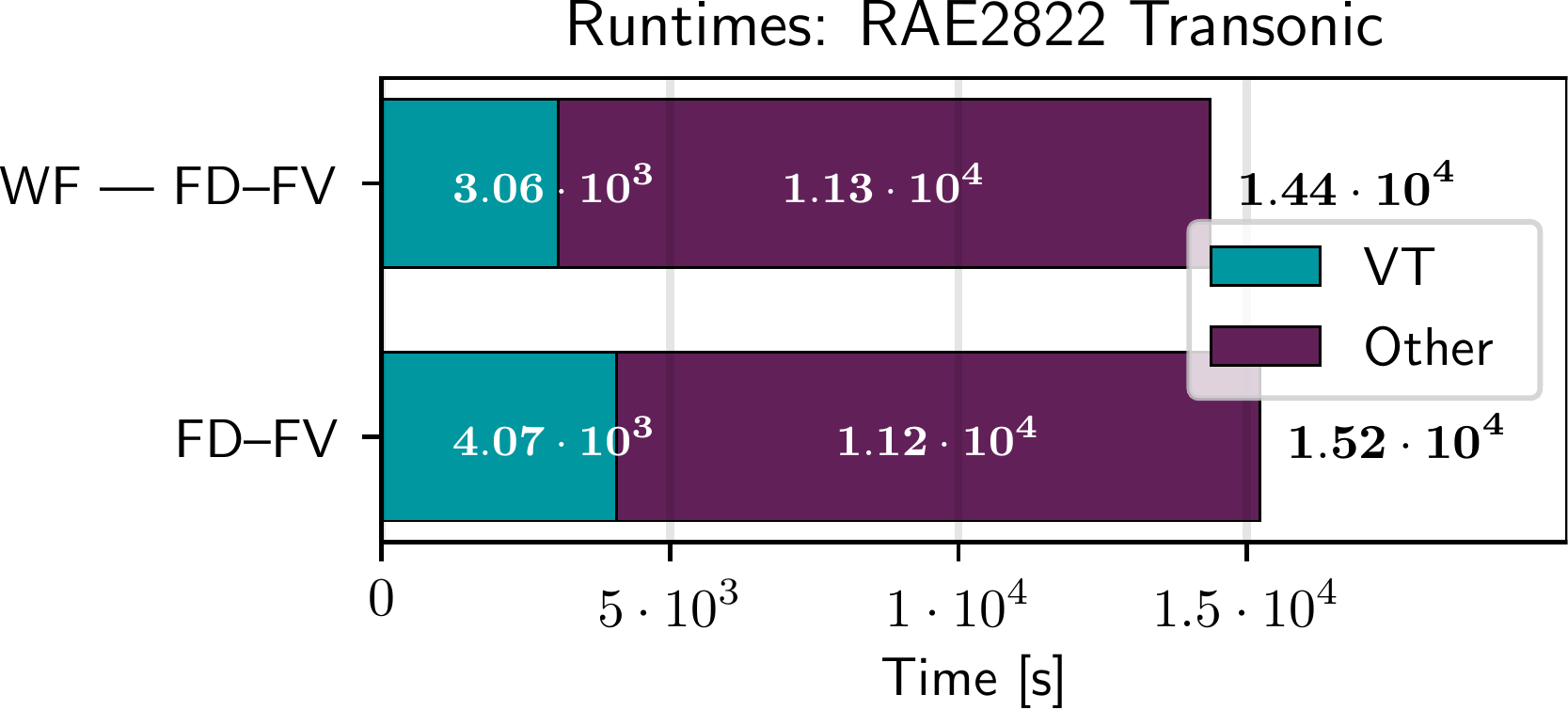}}
		\caption[Recorded runtime spent on volume integral/term computation for the transonic viscous flow over the RAE2822 airfoil.]
		{Total \ac{RHS} and \ac{VT} computation runtime for the transonic viscous flow over the RAE2822 airfoil for $t \in [25 t_c, 26 t_c]$.}
		\label{fig:RunTimes_RAE2822}
	\end{figure}

	The simulation is performed on the $8096$ element mesh provided online \cite{case22_1stHiCFD}.
	As for the ONERA M6 wing from \cref{subsubsec:ONERA_M6} a substantial amount of grid cells is only required to avoid influence of the farfield boundary conditions on the airfoil.
	Thus, we expect performance gains from employing the \ac{WF} volume term discretization in large parts of the domain, where the solution is almost constant and no stabilization is required.
	The solution is discretized using $p = 3$ polynomials, the surface fluxes are computed with the \ac{HLLC} flux \cite{toro1994restoration}, and the volume fluxes with the entropy--conservative and kinetic energy preserving flux from Ranocha \cite{Ranocha2020Entropy}.
	The simulation is first run until $t_f = 25 t_c$ with convective time $t_c \coloneqq t \frac{L}{U_\infty}$ with $L = 1$ and $U_\infty = \text{Ma}_\infty$.
	Explicit time stepping is performed with the adaptive third--order, four/three--stage embedded \ac{SSP} Runge--Kutta method from \cite{fekete2022embedded}, implemented in \texttt{OrdinaryDiffEq.jl} \cite{DifferentialEquations.jl-2017}.
	The timings for the \ac{WF}---\ac{DG}--\ac{FV} and \ac{FD}---\ac{DG}--\ac{FV}$\mathcal{O}2$ \cite{rueda2021entropy} volume term combinations are then recorded over the interval $t \in [25 t_c, 26 t_c]$.
	The required slope limiter for the second--order subcell finite volume scheme is chosen here as the minmod limiter. 
	The blending parameter $\beta$ is governed by the shock--capturing indicator from \cite{persson2006sub, hennemann2021provably} with indicator variable $\rho \cdot p$ and parameters $\beta_{\text{min}} = 0.001$, $\beta_\text{max} = 1.0$.
	This leads to about $14.1 \%$ of elements requiring stabilization with the blended \ac{DG}--\ac{FV} scheme, i.e., $\beta_i > 0$.
	As expected, the \ac{WF}---\ac{DG}--\ac{FV} combination is noticeably faster than the pure \ac{FD}---\ac{DG}--\ac{FV} combination, see \cref{fig:RunTimes_RAE2822}.
	This is due to the fact that large parts of the domain can be treated with the cheaper \ac{WF} volume term discretization as no stabilization is required there, cf. \cref{fig:RAE2822_alpha}.
	In particular, only in regions where the indicator from \cite{hennemann2021provably} indicates troubled cells which should be stabilized using the low--order \ac{FV} scheme (i.e., $\beta_i > 0$, recall \cref{eq:BlendedDG-FV}) the more expensive \ac{FD} volume term discretization is employed.
	The overall runtime savings, however, are comparably little which is due to the high cost of the parabolic part of the right--hand side which amounts to 50\% of the total runtime and is not affected by the choice of the volume term discretization for the convective part of the right--hand side.
	The inviscid volume integral contributes in this example only with $27.1\%$ or $36.3\%$ to the total runtime for the \ac{WF}---\ac{DG}--\ac{FV} and \ac{FD}---\ac{DG}--\ac{FV} combination, respectively.
	%
	%
	\section{Conclusions}
	\label{sec:Conclusions}
	We have presented a general framework for adaptively exchanging the volume term/integral discretization in every cell and Runge--Kutta stage of high--order \ac{DG} schemes.
	Based on a--priori and a--posteriori indicators, the cheap weak form volume integral discretization can be employed in non--critical regions, while the expensive and robust volume term discretizations like flux--differencing and finite--volume subcell limiting are used in cells where stabilization is required.
	Additionally, we have demonstrated that usage of entropy--dissipative \ac{WF} volume integrals increases the robustness of a discretization.
	An important property of this framework is that it strictly preserves the order of accuracy of the scheme.
	We also verified that the $v$--adaptive schemes maintain convergence to entropy solutions and can in practice remedy linear stability issues of entropy--conservative schemes.
	Furthermore, we have shown that tolerating some limited entropy increase due to the weak form volume integral can still give high--quality results.

	We applied the proposed framework to a variety of test cases governed by the compressible Euler and Navier--Stokes equations in two three space dimensions.
	In every case, the adaptive volume term discretization resulted in computational savings while maintaining stability.
	Precisely, we used both the a--posteriori entropy--production based indicator to combine weak form and flux--differencing integrals, as well as the a--priori shock--capturing indicator to combine the weak form and blended \ac{DG}--\ac{FV} scheme.
	Depending on the test case, the computational cost of the volume term computation can be reduced by factors up to $3$.
	In general, performance gains for the Navier--Stokes equations are less pronounced due to the fact that there are more gradient variables than conserved variables, where only the latter require the more expensive \ac{FD} volume term discretization.
	For the examples considered here, we showed that the volume integral cost can be reduced by about $25$ to $50\%$.
	
	Although we focused on the \ac{DGSEM} with collocated, tensor--product operators on quadrilateral and hexahedral meshes, the proposed framework is not limited to this particular discretization and can be applied to more general \ac{DG} schemes on unstructured meshes as well.
	In fact, the $v$--adaptive schemes are most beneficial for discretizations where the \ac{VT/VI} operators admit no tensor--product structure, e.g., for \ac{DG} schemes on simplices with non--collocated nodes.

	Consequently, we will focus on constructing efficient $v$--adaptive schemes with shock--capturing on simplices in future work.
	In particular, nested, i.e., combined weak form, flux--differencing, and blended \ac{DG}--\ac{FV} schemes on triangles and tetrahedra will be of interest.
	Furthermore, we will explore the adaptive combination of \ac{WF}, \ac{FD} and invariant domain preserving 
	\cite{pazner2021sparse, rueda2022subcell}
	and subcell monolithic convex limiting schemes 
	\cite{hajduk2021monolithic, rueda2024monolithic}.
	Finally, it is interesting to explore adaptivity in quadrature rules for the involved volume integrals as well.
	Potentially, through using e.g. Gauss--Legendre instead of Gauss--Lobatto--Legendre quadrature, even more cells could be treated with the \ac{WF} volume integral instead of the \ac{FD} volume term discretization.
	\section*{Data Availability}
	All data generated or analyzed during this study are included in this published article and its supplementary information files.
	\section*{Code Availability \& Reproducibility}
	We provide a reproducibility repository publicly available on GitHub \cite{doehring2026VTA_ReproRepo}.
	\section*{Acknowledgments}
	The authors thank Andr{\'e}s Rueda--Ram{\'\i}rez for fruitful discussions and exchange of ideas on the topic of this work.\\
	Daniel Doehring, Michael Schlottke--Lakemper, Manuel Torrilhon, and Gregor Gassner acknowledge funding by German Research Foundation (DFG) under Research Unit FOR5409: \\
	"Structure-Preserving Numerical Methods for Bulk- and Interface-Coupling of Heterogeneous Models ~(SNuBIC)" (Grants \#463312734; \#528753982).\\
	Hendrik Ranocha and Manuel Torrilhon acknowledge funding by German Research Foundation (DFG) under Research Unit SPP2410:\\
	"Hyperbolic Balance Laws in Fluid Mechanics: Complexity, Scales, Randomness (CoScaRa)" (Grants \#526031774; \#525660607).\\
	Gregor Gassner acknowledges funding by the German Reserach Foundation (DFG) under Germany´s Excellence Strategy - EXC 3037 - \#533607693.
	Gregor Gassner further acknowledges funding through the German Federal Ministry for Education and Research (BMFTR) project "ICON-DG" (\#01LK2315B) of the "WarmWorld Smarter" program.
	%
	%
	\section*{Declaration of competing interest}
	The authors declare the following financial interests/personal relationships which may be considered as potential competing interests:
	Daniel Doehring's financial support was provided by German Research Foundation.
	\section*{CRediT authorship contribution statement}
	\noindent
	\textbf{Daniel Doehring}: Conceptualization, Methodology, Investigation, Software, Validation, Writing - original draft. \\
	\textbf{Jesse Chan}: Conceptualization, Methodology, Software, Writing - review \& editing. \\
	\textbf{Hendrik Ranocha}: Conceptualization, Methodology, Software, Writing - review \& editing. \\
	\textbf{Michael Schlottke--Lakemper}: Methodology, Software, Investigation, Funding acquisition, Writing - review \& editing. \\
	\textbf{Manuel Torrilhon}: Funding acquisition, Supervision, Writing - review \& editing. \\
	\textbf{Gregor Gassner}: Conceptualization, Methodology, Investigation, Funding acquisition, Writing - review \& editing.
	\newpage
	\section*{Supplementary Figures}
	\begin{figure}[ht]
		\includegraphics[width=0.46\textwidth]{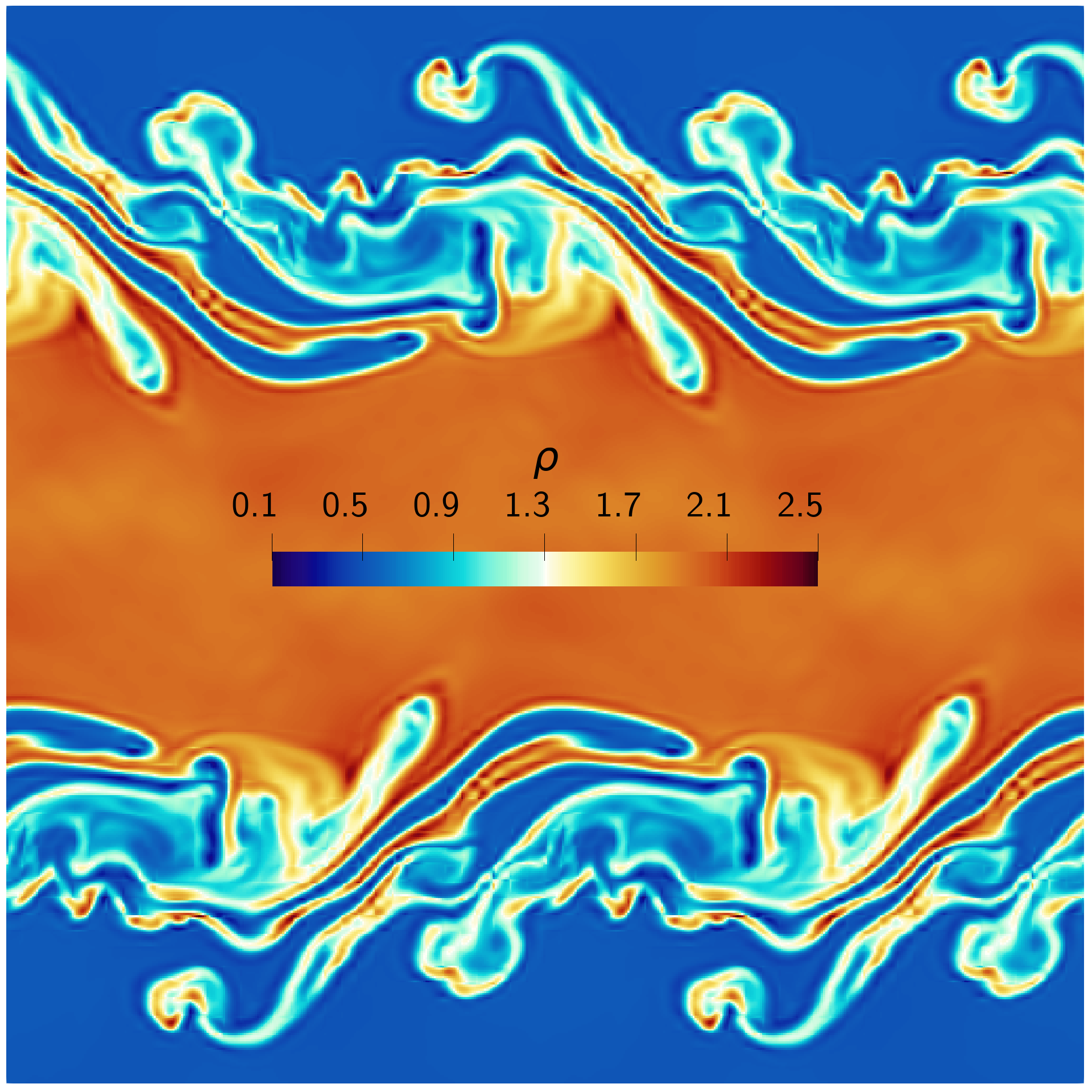}
		\caption[]{Density at final time $t_f = 5.25$ for the Kelvin--Helmholtz instability from \cref{subsubsec:KHI_DGSEM_Rigorous} obtained with the entropy--diffusive adaptive scheme.}
		\label{fig:KHI_Rigorous_DGSEM_rho}
	\end{figure}
	\begin{figure}[ht]
		\includegraphics[width=0.46\textwidth]{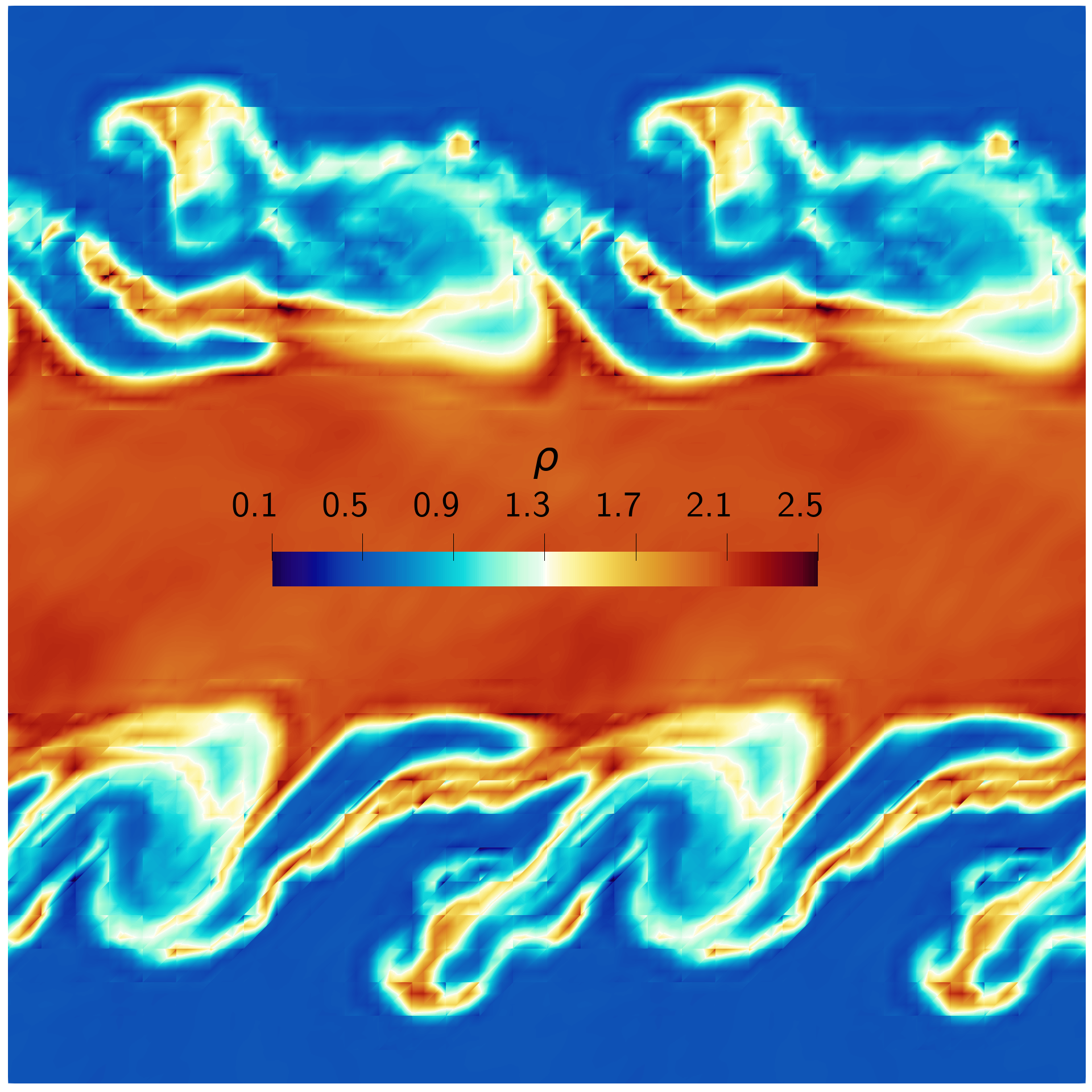}
		\caption[]{Density at final time $t_f = 4.6$ for the Kelvin--Helmholtz instability from \cref{subsubsec:KHI_DGMulti_Heuristic} obtained with the heuristic adaptive scheme.}
		\label{fig:KHI_Heuristic_DGMulti_rho}
	\end{figure}
	\bibliographystyle{elsarticle-num-names} 
	\newpage 
	\bibliography{references.bib}
	\end{document}